\documentclass[11pt]{amsart}
\usepackage{amssymb}
\usepackage{amsmath}
\usepackage{fancyhdr}
\usepackage[british]{babel}
\usepackage{geometry}
\usepackage{enumitem}
\usepackage{algpseudocode}
\usepackage{dsfont}
\usepackage{centernot}
\usepackage{xstring}
\usepackage{colortbl}
\usepackage[section]{algorithm}
\usepackage{graphicx}
\usepackage[font={footnotesize}]{caption}
\usepackage[usenames,dvipsnames,table]{xcolor}
\usepackage{pstricks,tikz}
\usepackage[h]{esvect}
\usepackage[unicode,
		bookmarksopen=true,
		bookmarksopenlevel=1,
		colorlinks=true,
		linkcolor=darkblue,
        linktoc=page,
		citecolor=darkblue,
]{hyperref}

\title[]{Decompositions into isomorphic rainbow spanning trees}
\date{\today}
\author[S.~Glock, D. K\"uhn, R. Montgomery, D. Osthus]{Stefan Glock, Daniela K\"uhn, Richard Montgomery and Deryk Osthus}
\address{School of Mathematics, University of Birmingham,
Edgbaston, Birmingham, B15 2TT, United Kingdom}
\email{[s.glock,d.kuhn,r.h.montgomery,d.osthus]@bham.ac.uk}
\thanks{This project has received partial funding from the European Research
Council (ERC) under the European Union's Horizon 2020 research and innovation programme (grant agreement no. 786198, D.~K\"uhn and D.~Osthus).
The research leading to these results was also partially supported by the EPSRC, grant nos. EP/N019504/1 (S.~Glock and D.~K\"uhn) and EP/S00100X/1  (D.~Osthus),
as well as the Royal Society and the Wolfson Foundation (D.~K\"uhn).}

\newtheorem{theorem}{Theorem}

\newtheorem{lemma}[theorem]{Lemma}
\newtheorem{cor}[theorem]{Corollary}
\newtheorem{fact}[theorem]{Fact}
\newtheorem{conj}[theorem]{Conjecture}

\theoremstyle{definition}

\newtheorem{defin}[theorem]{Definition}

\newtheoremstyle{claimstyle}{5pt}{5pt}{\em}{5pt}{\em}{:}{5pt}{}
\theoremstyle{claimstyle}
\newtheorem{claim}{Claim}

\newtheoremstyle{stepstyle}{10pt}{5pt}{\em}{0pt}{\em}{:}{5pt}{}
\theoremstyle{stepstyle}
\newtheorem{step}{Step}

\definecolor{darkblue}{rgb}{0,0,0.5}

\def\noproof{{\unskip\nobreak\hfill\penalty50\hskip2em\hbox{}\nobreak\hfill%
       $\square$\parfillskip=0pt\finalhyphendemerits=0\par}\goodbreak}
\def\endproof{\noproof\bigskip}

\def\noclaimproof{{\unskip\nobreak\hfill\penalty50\hskip2em\hbox{}\nobreak\hfill%
       $-$\parfillskip=0pt\finalhyphendemerits=0\par}\goodbreak}

\def\endclaimproof{\noclaimproof\medskip}

\newdimen\margin
\def\textno#1&#2\par{
   \margin=\hsize
   \advance\margin by -4\parindent
          \setbox1=\hbox{\sl#1}
   \ifdim\wd1 < \margin
      $$\box1\eqno#2$$
   \else
      \bigbreak
      \hbox to \hsize{\indent$\vcenter{\advance\hsize by -3\parindent
      \it\noindent#1}\hfil#2$}
      \bigbreak
   \fi}

\def\proof{\removelastskip\penalty55\medskip\noindent\setcounter{claim}{0}\setcounter{step}{0}{\bf Proof. }} 

\def\lateproof#1{\removelastskip\penalty55\medskip\noindent\setcounter{claim}{0}\setcounter{step}{0}{\bf Proof of #1. }} 

\DeclareMathOperator{\Ima}{Im}
\def\claimproof{\removelastskip\penalty55\medskip\noindent{\em Proof of claim: }}

\begin{document}

\def\COMMENT#1{}
\def\TASK#1{}

\def\eps{{\varepsilon}}
\newcommand{\ex}{\mathbb{E}}
\newcommand{\pr}{\mathbb{P}}
\newcommand{\cB}{\mathcal{B}}
\newcommand{\cA}{\mathcal{A}}
\newcommand{\cE}{\mathcal{E}}
\newcommand{\cS}{\mathcal{S}}
\newcommand{\cF}{\mathcal{F}}
\newcommand{\cG}{\mathcal{G}}
\newcommand{\bL}{\mathbb{L}}
\newcommand{\bF}{\mathbb{F}}
\newcommand{\bZ}{\mathbb{Z}}
\newcommand{\cH}{\mathcal{H}}
\newcommand{\cC}{\mathcal{C}}
\newcommand{\cM}{\mathcal{M}}
\newcommand{\bN}{\mathbb{N}}
\newcommand{\bR}{\mathbb{R}}
\newcommand{\rhosub}{\eta}
\def\O{\mathcal{R}}
\def\cO{\mathcal{R}}
\newcommand{\cP}{\mathcal{P}}
\newcommand{\cQ}{\mathcal{Q}}
\newcommand{\cR}{\mathcal{R}}
\newcommand{\cJ}{\mathcal{J}}
\newcommand{\cL}{\mathcal{L}}
\newcommand{\cK}{\mathcal{K}}
\newcommand{\cD}{\mathcal{D}}
\newcommand{\cI}{\mathcal{I}}
\newcommand{\cV}{\mathcal{V}}
\newcommand{\cT}{\mathcal{T}}
\newcommand{\cU}{\mathcal{U}}
\newcommand{\cX}{\mathcal{X}}
\newcommand{\cW}{\mathcal{W}}
\newcommand{\cZ}{\mathcal{Z}}
\newcommand{\1}{{\bf 1}_{n\not\equiv \xi}}
\newcommand{\eul}{{\rm e}}
\newcommand{\Erd}{Erd\H{o}s}
\newcommand{\cupdot}{\mathbin{\mathaccent\cdot\cup}}
\newcommand{\whp}{whp }

\newcommand{\ind}[1]{\mathds{1}_{\Set{#1}}}

\newcommand{\doublesquig}{%
  \mathrel{%
    \vcenter{\offinterlineskip
      \ialign{##\cr$\rightsquigarrow$\cr\noalign{\kern-1.5pt}$\rightsquigarrow$\cr}%
    }%
  }%
}

\newcommand{\defn}{\emph}

\newcommand\restrict[1]{\raisebox{-.5ex}{$|$}_{#1}}

\newcommand{\prob}[1]{\mathrm{\mathbb{P}}\left(#1\right)}
\newcommand{\cprob}[2]{\mathrm{\mathbb{P}}_{#1}\left(#2\right)}
\newcommand{\expn}[1]{\mathrm{\mathbb{E}}\left(#1\right)}
\newcommand{\cexpn}[2]{\mathrm{\mathbb{E}}_{#1}\left(#2\right)}
\def\gnp{G_{n,p}}
\def\G{\mathcal{G}}
\def\lflr{\left\lfloor}
\def\rflr{\right\rfloor}
\def\lcl{\left\lceil}
\def\rcl{\right\rceil}

\newcommand{\qbinom}[2]{\binom{#1}{#2}_{\!q}}
\newcommand{\binomdim}[2]{\binom{#1}{#2}_{\!\dim}}

\newcommand{\grass}{\mathrm{Gr}}

\newcommand{\brackets}[1]{\left(#1\right)}
\def\sm{\setminus}
\newcommand{\Set}[1]{\{#1\}}
\newcommand{\set}[2]{\{#1\,:\;#2\}}
\newcommand{\krq}[2]{K^{(#1)}_{#2}}
\def\In{\subseteq}

\begin{abstract}
A subgraph of an edge-coloured graph is called rainbow if all its edges have distinct colours. Our main result implies that, given any optimal colouring of a sufficiently large complete graph $K_{2n}$, there exists a decomposition of $K_{2n}$ into isomorphic rainbow spanning trees. This settles conjectures of Brualdi--Hollingsworth (from 1996) and Constantine (from 2002) for large graphs.
\end{abstract}

\maketitle

\section{Introduction}

Given an edge-coloured graph~$G$, we say a subgraph $H$ is \defn{rainbow} if all the edges of $H$ have distinct colours. Moreover, we say that $H_1,\dots,H_t$ \defn{decompose $G$} if $H_1,\dots,H_t$ are edge-disjoint subgraphs of $G$ covering all the edges of~$G$.

In this paper, we address the problem of decomposing an optimally edge-coloured complete graph $K_{2n}$ into (isomorphic) rainbow spanning trees.
The study of rainbow decomposition problems can be traced back to the work of Euler, who investigated
for which $n$ one can find a pair of orthogonal Latin squares of order $n$.
That is, equivalently, for which $n$ does there exist an optimally edge-coloured $K_{n,n}$ which can be decomposed into rainbow perfect matchings? Euler gave a construction for all $n\not\equiv 2\mod{4}$ and conjectured that these are the only admissible values. His conjecture was disproved by Parker, Bose and Shrikhande who provided constructions for the missing values, except for $n=6$ (which corresponds to Euler's famous `36 officers problem', for which the non-existence had already been shown by Tarry in 1901).

On the other hand, given an \emph{arbitrary} optimally edge-coloured $K_{n,n}$, a decomposition into rainbow perfect matchings need not exist. In fact, there are examples of such colourings that do not admit a single rainbow perfect matching. (An important conjecture widely attributed to Ryser--Brualdi--Stein postulates that there always exists a rainbow matching of size $n-1$.)

Perfect matchings are, in some sense, very rigid objects, and it is natural to ask analogous questions for other types of subgraphs. In particular, several natural conjectures arose concerning decompositions into rainbow spanning trees.
Here, the most notable are the Brualdi--Hollingsworth conjecture,
Constantine's conjecture and the Kaneko--Kano--Suzuki conjecture.
Our main result implies the first two of these.

\subsection{Decompositions into rainbow spanning trees}

Note that if $K_{2n}$ is optimally edge-coloured, then the colour classes form a \emph{$1$-factorization}, that is, a decomposition of $K_{2n}$ into perfect matchings. We will here use the term $1$-factorization synonymously with an edge-colouring whose colour classes form a $1$-factorization. Note that if a $1$-factorization of $K_n$ exists, then $n$ is even. We now state the Brualdi--Hollingsworth conjecture.

\begin{conj}[Brualdi and Hollingsworth, \cite{BH:96}]\label{conj:BH}
For all $n>4$ and any $1$-factorization of~$K_n$, there exists a decomposition of $K_n$ into rainbow spanning trees.
\end{conj}
Note that the condition that $n>4$ is necessary.%
\COMMENT{Inserting $n>4$ here and in Conj 3 is not part of the original conjectures,
but it seems odd to omit this.}
Suppose a $1$-factorization of~$K_n$ is given. Clearly, there always exists \emph{one} rainbow spanning tree, for instance the star at any vertex.
Brualdi and Hollingsworth in their original paper~\cite{BH:96} showed that one can find two edge-disjoint rainbow spanning trees. Shortly afterwards, Krussel, Marshall and Verrall~\cite{KMV:00} were able to find three such trees. Later, Horn~\cite{horn:18} significantly improved on this by finding $\Omega(n)$ edge-disjoint rainbow spanning trees.
Very recently, Montgomery, Pokrovskiy, and Sudakov~\cite{MPS:18b} proved Conjecture~\ref{conj:BH} approximately by showing that one can guarantee $(1-o(1))n/2$ edge-disjoint rainbow spanning trees.

Several related problems have been studied, in two main directions. Firstly, we may wish to strengthen the conditions on the trees in the decomposition, most commonly by requiring the trees in the decomposition to be isomorphic. Secondly, we may wish to weaken the conditions on the colouring, most commonly by allowing non-optimal proper edge-colourings. (Naturally, those two directions can also be combined.)

A decomposition of an edge-coloured $K_n$ into isomorphic rainbow spanning trees is also known in the literature as a \defn{multicoloured tree parallelism (MTP)}.
It turns out that the problem of finding an MTP is non-trivial even if one is allowed to choose the $1$-factorization. Partial results were obtained by Constantine~\cite{constantine:02}.
Akbari, Alipour, Fu and~Lo~\cite{AAFL:06} then proved that for all $n\in \bN$ with $n>2$, there exists a $1$-factorization of $K_{2n}$ which admits an MTP.
Moreover, Constantine conjectured that in fact an MTP should exist for any given $1$-factorization, thus generalizing the Brualdi--Hollingsworth conjecture.

\begin{conj}[Constantine,~\cite{constantine:02,constantine:05}]\label{conj:constantine}
For all $n>4$, any $1$-factorization of $K_{n}$ admits a decomposition into isomorphic rainbow spanning trees.
\end{conj}
\COMMENT{
DO: The 2002 paper writes `We color the edges of $K_{2n}$ with $2n-1$
colors by assigning one color to each edge. The coloration is proper if whenever two edges that have one vertex in common carry
different colors.' So he clearly refers to $1$-factorizations in that paper.}

Unsurprisingly, for this conjecture, much less was known than for the Brualdi--Hollingsworth conjecture. In~\cite{FL:15} it was shown that three isomorphic rainbow spanning trees can be guaranteed. In~\cite{PS:18}, Pokrovskiy and Sudakov showed that one can find $10^{-6}n$ edge-disjoint rainbow spanning trees all isomorphic to a so-called \emph{$t$-spider} (which is even independent of the given $1$-factorization).
Montgomery, Pokrovskiy, and Sudakov~\cite{MPS:18b}, and independently, Kim, K\"uhn, Kupavskii and Osthus~\cite{KKKO:18}, proved a weak asymptotic version of the conjecture by showing that there are $(1-o(1))n/2$ edge-disjoint rainbow paths each of length $(1-o(1))n$.

We now discuss results on more general colourings. Intuitively it might seem that dealing with an optimal colouring is the hardest case, as having more colours should make finding rainbow subgraphs easier. However, non-optimal colourings seem genuinely harder to deal with than $1$-factorizations.
Kaneko, Kano and Suzuki proved that for any proper colouring of~$K_n$, there exist three edge-disjoint rainbow spanning trees, and also generalized the Brualdi--Hollingsworth conjecture as follows.
\begin{conj}[Kaneko, Kano, and Suzuki, 2002]\label{conj:KKS}
For all $n>4$, every properly edge-coloured $K_{n}$ contains $\lfloor n/2\rfloor$ edge-disjoint rainbow spanning trees.
\end{conj}
Note that any proper colouring is \defn{$n/2$-bounded}, that is, every colour appears on at most $n/2$ edges. Under the weaker assumption that the colouring is $n/2$-bounded, Akbari and Alipour~\cite{AA:07} showed that one can guarantee two edge-disjoint rainbow spanning trees, and this was significantly improved by Carraher, Hartke, and Horn~\cite{CHH:16} who showed that $\Omega(n/\log n)$ such trees exist. For proper colourings, a linear number of rainbow spanning trees was independently obtained by Pokrovskiy and Sudakov~\cite{PS:18} and by Balogh, Liu and Montgomery~\cite{BLM:18}, where in the former work, the trees are even isomorphic. Finally, the aforementioned result from \cite{MPS:18b} on Conjecture~\ref{conj:BH} also applies to proper colourings, thus proving Conjecture~\ref{conj:KKS} approximately.

We now state our main theorem, which implies the Brualdi--Hollingsworth conjecture and Constantine's conjecture for large~$n$. This is the first general exact rainbow decomposition result for spanning subgraphs, where each subgraph in the decomposition has to use all the colours.

\begin{theorem} \label{thm:main}
For all sufficiently large~$n$, there exists a tree $T$ on $n$ vertices such that for any $1$-factorization\COMMENT{this implicitly always implies that $n$ is even} of~$K_n$, there exists a decomposition into rainbow subgraphs each isomorphic to~$T$.
\end{theorem}

Note that whereas Constantine's conjecture says that given a $1$-factorization, one can decompose into isomorphic rainbow spanning trees, we actually show that one can use the same tree $T$ for \emph{any} $1$-factorization. This tree is made up of a path of length $n-o(n)$, with $o(n)$ short paths attached to it (see~Definition~\ref{def:tree}).
By modifying our proof slightly, we can even ensure that $\Delta(T)\le 3$. This is best possible in the sense that there exist $1$-factorizations which do not admit a single rainbow Hamilton path~\cite{MM:84}.

Our argument relies upon the fact that the colouring is a $1$-factorization. It would be very interesting to prove the result for more general colourings, in particular proper colourings.

It would also be interesting to investigate the $n/2$-bounded setting further. The best known bound is the one from~\cite{CHH:16} mentioned earlier, which provides $\Omega(n/\log n)$ edge-disjoint rainbow spanning trees. A natural question is to ask for the maximum number $k$ of such trees that can be guaranteed. It seems unlikely that a decomposition can be obtained, but it would be interesting to see whether $k=\Omega(n)$ is possible or not.
It is also natural to impose further local conditions on the colouring, e.g.~that the colouring is \defn{locally $\Delta$-bounded}, which means that the maximum degree of each colour class is at most~$\Delta$. For instance, in~\cite{KKKO:18} it is shown that for any $n/2$-bounded colouring which is locally $o(n)$-bounded, there exists an approximate decomposition into almost spanning rainbow cycles (and thus into almost spanning paths).

\subsection{Related problems}\label{subsec:related}

We now discuss some related results concerning rainbow decompositions. Let us first revisit the perfect matching case. As mentioned earlier, there exist proper optimal colourings of $K_{n,n}$ which do not contain a rainbow perfect matching. However, by imposing slightly stronger boundedness conditions on the colouring, one can obtain strong results. For example, Alon, Spencer and Tetali~\cite{AST:95} showed that if $n$ is a power of $2$ and the edge-colouring is $o(n)$-bounded (and not necessarily proper), there exists a decomposition into rainbow perfect matchings. Montgomery, Pokrovskiy and Sudakov~\cite{MPS:18b} showed that any proper edge-colouring of $K_{n,n}$, where at most $(1-o(1))n$ colours appear more than $(1-o(1))n$ times, contains $(1-o(1))n$ edge-disjoint rainbow perfect matchings. This implied a conjecture of Akbari and Alipour in a strong form (which was proved independently
by Keevash and Yepremyan~\cite{KY:18}) and a conjecture of Barat and Nagy approximately, both for large~$n$.
Kim, K\"uhn, Kupavskii and Osthus~\cite{KKKO:18} proved that for any $(1-o(1))n$-bounded and locally $o(n/\log^2 n)$-bounded edge-colouring of $K_{n,n}$, there exist $(1-o(1))n$ edge-disjoint rainbow perfect matchings.
The authors of both~\cite{KKKO:18} and \cite{MPS:18b} also obtain analogous
results (in their respective settings) on approximate decompositions of $K_n$ into rainbow Hamilton cycles. Furthermore,~\cite{KKKO:18} contains results for approximate decompositions of $K_n$ into rainbow $F$-factors (for any given graph~$F$).


A further tantalizing problem concerning rainbow tree decompositions is the following special case of Rota's basis conjecture. Let $T_1,\dots,T_{n-1}$ be spanning trees on a common vertex set of size~$n$, each monochromatic in a different colour. Then their union (allowing multiple edges) can be decomposed into $n-1$ rainbow spanning trees. The general version of Rota's conjecture concerns the rearrangement of bases of a matroid into disjoint transversal bases. Recently, Buci\'{c}, Kwan, Pokrovskiy and Sudakov~\cite{BKPS:18} showed that $(1/2-o(1))n$ disjoint transversal bases can be found.

\section{Notation}

Given a graph~$G$ with edge colouring $\phi \colon E(G)\to C$, we say a subgraph $H$ is \defn{$D$-rainbow} if $H$ is rainbow and $\phi(E(H))=D$. We refer to an edge $e=uv$ with colour $c$ as a \defn{$c$-edge}, and $v$ is a \defn{$c$-neighbour} of~$u$. For each colour $c$, $E_c(G)$ is the set of $c$-edges in~$G$.
For each vertex $v$ of $G$, we let $\partial_G(v)$ denote the set of all edges of $G$ incident to~$v$. For any $S\subseteq V(G)$, $N_G(S)$ is the common neighbourhood in $G$ of the vertices in $S$. For any $x\in V(G)$ and $U\subseteq V(G)$, $d_G(x,U)$ is the number of neighbours of $x$ in~$U$. We denote by $G-H$ the graph obtained from $G$ by deleting the edges of~$H$.

For a hypergraph $\cH$, let $\Delta^c(\cH)$ denote its maximum codegree, that is, the maximum number of edges containing any two fixed vertices.

Given a set $X$ and $p\in[0,1]$, a \defn{$p$-random subset} is a random subset $Y\In X$ which is obtained by including each element of $X$ independently with probability~$p$. If not otherwise stated, we always assume that such random subsets are independent. For instance, if we say that $Y$ is a $p$-random subset of $X$ and $Y'$ is a $p'$-random subset of $Y$, we implicitly assume that these random choices are made independently.
Similarly, if $G$ is a graph, then a \defn{$p$-random subgraph} is the random graph with vertex set $V(G)$ and a $p$-random subset of $E(G)$ as edge set.

On the other hand, we often split a random subset further into disjoint subsets. For instance, if $Y$ is a $(p+p')$-random subset of $X$, we might say that we split $Y$ into a $p$-random set $Y_1$ and a $p'$-random set $Y_2$, by which we mean that for each $y\in Y$ independently, we include $y$ in $Y_1$ with probability~$p/(p+p')$ and into $Y_2$ otherwise. Note that then $Y_1$ is indeed a $p$-random subset of $X$ and $Y_2$ is a $p'$-random subset of~$X$, but they are obviously not independent. To split into more sets, we use the following notation: By splitting $X$ randomly as
$$\begin{array}{cccccccccccccccccccccccccccccc}
  X   &=&  X_1  &\cupdot&	   \dots  &\cupdot&	  X_m\\
	1   &=&  p_1  &+& \dots &+&  p_m
\end{array}$$
we mean that for every element in~$X$ independently, we choose an index $i\in [m]$ according to the probability distribution $(p_i)_{i=1}^m$, and put this element into the corresponding set~$X_i$.

We say that a random event holds \defn{with high probability} if the probability that it holds tends to $1$ as $n$ tends to infinity (where $n$ is usually the number of vertices and the event depends on~$n$).

We write $[n]:=\Set{1,\dots,n}$. For $a,b,c\in \mathbb{R}$,
we write $a=b\pm c$ whenever $a\in [b-c,b+c]$.
For $a,b,c\in (0,1]$,
we write $a\ll b \ll c$ in our statements to mean that there are increasing functions $f,g:(0,1]\to (0,1]$
such that whenever $a\leq f(b)$ and $b \leq g(c)$,
then the subsequent result holds.


\section{Proof sketch}

Our proof is based on hypergraph matching results and new absorption techniques.
Suppose we are given a 1-factorization $\phi$ of the complete graph $K_n$ with colour set~$C$.
We build the $t:=n/2$ rainbow trees simultaneously, beginning with our absorbing structures and then gradually extending these structures to cover all the vertices and edges. For this, we further develop a recent `distributive' form of the absorption method: we form an absorption structure along with a reservoir, such that, given any subset
(of given size) from the reservoir we can distribute the elements of this subset among the different parts of the absorbing structure to always obtain a copy of the same tree. More precisely, we create a `global' reservoir of edges, as well as `local'  reservoirs of colours and vertices
(as explained below, `local' refers to the fact that there is one such reservoir for each tree, while the `global' reservoir is common to all trees). The structure of these
absorbers and the corresponding reservoirs is described in more detail in Section~\ref{sec:absorb}.

 Already, however, we can outline our proof strategy, as follows.

\begin{enumerate}
\item Create an edge absorption structure and a global edge reservoir.
\item For each tree, create a colour absorption structure and a colour reservoir.
\item For each tree, create a vertex absorption structure and a vertex reservoir.
\item Find $t=n/2$ edge-disjoint almost spanning rainbow paths $P_i$ covering most of the remaining vertices.
\item Link up the absorbers and paths to form $t$ rainbow forests~$F_i$ and thereby cover all non-reservoir vertices.
\item Cover non-reservoir edges by adding each such edge to one of the forests $F_i$.
\item Incorporate non-reservoir colours for each forest, by adding a suitable edge from the edge reservoir.
\item Absorb the uncovered reservoir vertices into each forest, using edges and colours from the reservoirs.
\item Absorb the uncovered reservoir colours into each forest, using the colour absorption structure.
 \item Absorb the uncovered edges from the global edge reservoir by distributing them among the forests to complete these forests into rainbow spanning trees~$T_i$.
\end{enumerate}
To find all of the structures we use, we apply results on matchings in certain auxiliary hypergraphs, as described in Section~\ref{sec:hypermatch}.
This allows the structures we find to look random-like, which in turn means that at each stage of the construction
of the trees $T_i$, the currently unused sets are also random-like. In particular, this
means that the leftover sets which need to be absorbed are sufficiently small and well-distributed
(again, the sets we track here are vertices and edges as well as colour sets).

The main difficulty in our proof lies in obtaining a decomposition into spanning trees.
The property that these trees are isomorphic (even to some $T$ fixed in advance)
can be achieved with only a little extra care. We comment more on this in Section~\ref{sec:iso}.
In Section~\ref{sec:tools}, we list the tools that we use in our proof.

The above strategy is implemented in Section~\ref{sec:main}, following the proof of several lemmas allowing some of these tasks. In Section~\ref{sec:approx}, we find a set of almost spanning rainbow paths. In Section~\ref{sec:colour}, we find our colour absorption structure. In Section~\ref{sec:edge}, we find our edge absorption structure. In Section~\ref{sec:connect}, we show how we will connect these structures together. In Section~\ref{sec:vertex}, we find suitable rainbow  matchings which we will use to absorb vertices.

\subsection{Designing absorbers}\label{sec:absorb}
The absorbing method has its roots in work by Erd\H{o}s, Gy\'arf\'as, and Pyber,
as well as Krivelevich, before its general codification by R\"odl, Ruci\'nski and Szemer\'edi. The key novelty in our work is to construct a `nested' absorbing structure for the edges, colours and vertices. As the edges of a tree define its colours and vertices, we start by building an edge absorption structure and an accompanying edge reservoir
(i.e.~the edges in the reservoir are those which can later be absorbed).

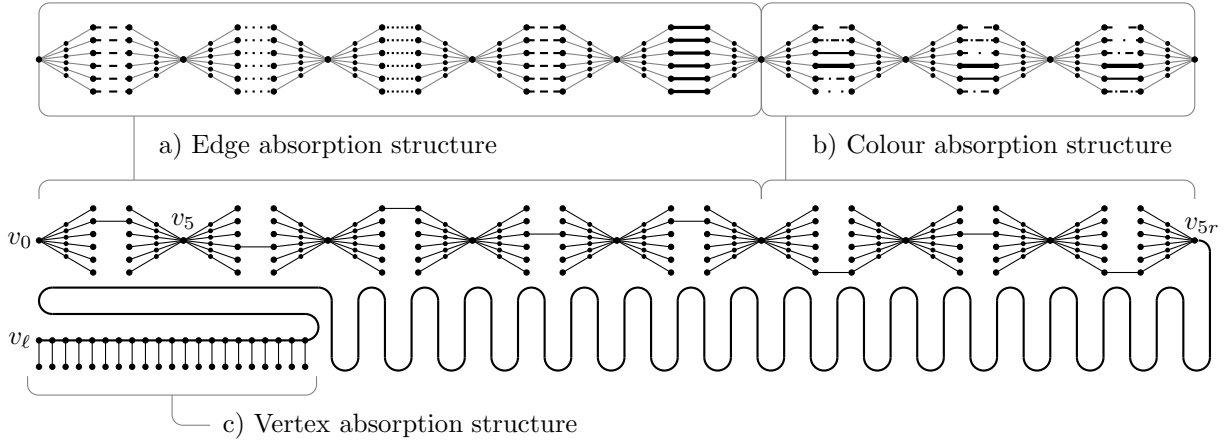
\begin{figure}[t]
\hspace{-0.3cm}
\begin{tikzpicture}[scale=1]
\def\vxrad{0.035cm}
\def\vxradoth{0.025cm}
\def\baseunit{0.475}
\def\yspacer{0.0425}
\def\spaceit{4*\baseunit}
\def\highery{5.6}
\def\yheight{0.75}
\def\yheightoth{4}
\def\textspacer{0.4}
\def\dropper{0.2}

\draw [rounded corners,gray] (2,{\highery-\yheight}) -- ({5*\spaceit+0.25},{\highery-\yheight}) -- ({5*\spaceit+0.25},{\highery+\yheight}) -- (0.25,{\highery+\yheight}) -- (0.25,{\highery-\yheight}) -- (2,{\highery-\yheight});
\draw [gray] (1.5,\yheightoth) -- (1.5,{\highery-\yheight});
\draw [rounded corners,gray] (0.25,{\yheightoth-0.25}) -- (0.25,\yheightoth) -- ({5*\spaceit+0.25},\yheightoth) -- ({5*\spaceit+0.25},{\yheightoth-0.25});

\draw [rounded corners,gray] (12,{\highery-\yheight}) -- ({8*\spaceit+0.25},{\highery-\yheight}) -- ({8*\spaceit+0.25},{\highery+\yheight}) -- ({5*\spaceit+0.25},{\highery+\yheight}) -- ({5*\spaceit+0.25},{\highery-\yheight}) -- (12,{\highery-\yheight});
\draw [gray] ({5.3*\spaceit},\yheightoth) -- ({5.3*\spaceit},{\highery-\yheight});
\draw [rounded corners,gray] ({8*\spaceit+0.25},{\yheightoth-0.25}) -- ({8*\spaceit+0.25},\yheightoth) -- ({5*\spaceit+0.25},\yheightoth) -- ({5*\spaceit+0.25},{\yheightoth-0.25});

\draw ({4.05},{\highery-\yheight-\textspacer}) node {\small a) Edge absorption structure};
\draw ({6.625*\spaceit+0.2},{\highery-\yheight-\textspacer}) node {\small b) Colour absorption structure};

\draw [rounded corners,gray]  (0.1,{1.35+0.25-\dropper}) -- (0.1,{1.35-\dropper}) -- (3.9,{1.35-\dropper}) -- (3.9,{1.35+0.25-\dropper});
\draw [rounded corners,gray] (2,{1.35-\dropper}) -- (2,{1.35-\dropper-0.4}) -- (2.5,{1.35-\dropper-0.4});
\draw (5,{1.35-\dropper-0.4}) node {\small c) Vertex absorption structure};

\draw (0,1.875) node {\small $v_\ell$};
\draw (0,3.2) node {\small $v_0$};
\draw ({\spaceit+0.25},3.45) node {\small $v_5$};
\draw (15.55,3.4) node {\small $v_{5r}$};

\def\colcyan{dashed}
\def\colmagenta{dotted}
\def\colteal{densely dotted}
\def\colEmerald{densely dashed}
\def\cololive{very thick}
\def\colred{loosely dotted}
\def\colblue{dashdotted}
\def\colgreen{ultra thick}
\def\colorange{thick}
\def\colblack{densely dashdotted}
\def\colpurple{loosely dashed}
\def\colviolet{loosely dashdotted}
\foreach \x/\matchingcolour in {0/{\colcyan},1/{\colmagenta},2/{\colteal},3/{\colEmerald},4/{\cololive}}
{
\foreach \y in {-5,-3,-1,1,3,5}
{
\draw [thick,\matchingcolour]  ({0.25+\baseunit*1.5+\x*\spaceit},{\highery+\y*2*\yspacer}) -- ({0.25+\baseunit*2.5+\x*\spaceit},{\highery+\y*2*\yspacer});
}
}

\foreach \x/\y/\colour in {5/-5/{\colred},5/-3/{\colblue},5/-1/{\colgreen},5/1/{\colorange},5/3/{\colblack},5/5/{\colpurple},
6/-5/{\colblue},6/-3/{\colorange},6/-1/{\colgreen},6/1/{\colviolet},6/3/{\colblack},6/5/{\colpurple},
7/-5/{\colblack},7/-3/{\colorange},7/-1/{\colgreen},7/1/{\colblue},7/3/{\colviolet},7/5/{\colpurple}}
{
\draw [thick,\colour]({0.25+\baseunit*1.5+\x*\spaceit},{\highery+\y*2*\yspacer}) -- ({0.25+\baseunit*2.5+\x*\spaceit},{\highery+\y*2*\yspacer});
}

\foreach \x in {0,...,7}
{
\foreach \y in {-5,-3,-1,1,3,5}
{
\draw [gray] ({0.25+\baseunit*0+\x*\spaceit},{\highery}) -- ({0.25+\baseunit*1.5+\x*\spaceit},{\highery+\y*2*\yspacer});
\draw [gray] ({0.25+\baseunit*4+\x*\spaceit},\highery) -- ({0.25+\baseunit*2.5+\x*\spaceit},{\highery+\y*2*\yspacer});

\draw [fill] ({0.25+\baseunit*0.75+\x*\spaceit},{\highery+\y*\yspacer}) circle [radius=\vxradoth];
\draw [fill] ({0.25+\baseunit*1.5+\x*\spaceit},{\highery+\y*2*\yspacer}) circle [radius=\vxrad];
\draw [fill] ({0.25+\baseunit*2.5+\x*\spaceit},{\highery+\y*2*\yspacer}) circle [radius=\vxrad];
\draw [fill] ({0.25+\baseunit*3.25+\x*\spaceit},{\highery+\y*\yspacer}) circle [radius=\vxradoth];
}
\draw [fill] ({0.25+\baseunit*0+\x*\spaceit},{\highery}) circle [radius=\vxrad];
\draw [fill] ({0.25+\baseunit*4+\x*\spaceit},\highery) circle [radius=\vxrad];
}

\foreach \x/\y in {0/3,1/-1,2/5,3/1,4/3,5/-5,6/1,7/-5}
{
\draw ({0.25+1.5*\baseunit+\x*\spaceit},{3.2+\y*\yspacer*2}) -- ({0.25+1.5*\baseunit+\baseunit+\x*\spaceit},{3.2+\y*\yspacer*2});
}

\foreach \x in {0,...,7}
{
\foreach \y in {-5,-3,-1,1,3,5}
{
\draw ({0.25+\baseunit*0+\x*\spaceit},{3.2}) -- ({0.25+\baseunit*1.5+\x*\spaceit},{3.2+\y*2*\yspacer});
\draw ({0.25+\baseunit*4+\x*\spaceit},3.2) -- ({0.25+\baseunit*2.5+\x*\spaceit},{3.2+\y*2*\yspacer});

\draw [fill] ({0.25+\baseunit*0.75+\x*\spaceit},{3.2+\y*\yspacer}) circle [radius=\vxradoth];
\draw [fill] ({0.25+\baseunit*1.5+\x*\spaceit},{3.2+\y*2*\yspacer}) circle [radius=\vxrad];
\draw [fill] ({0.25+\baseunit*2.5+\x*\spaceit},{3.2+\y*2*\yspacer}) circle [radius=\vxrad];
\draw [fill] ({0.25+\baseunit*3.25+\x*\spaceit},{3.2+\y*\yspacer}) circle [radius=\vxradoth];
}
\draw [fill] ({0.25+\baseunit*0+\x*\spaceit},{3.2}) circle [radius=\vxrad];
\draw [fill] ({0.25+\baseunit*4+\x*\spaceit},3.2) circle [radius=\vxrad];
}

\draw [thick] (15.5,3.2) arc (90:0:0.15);
\draw [thick] (15.65,3.05) -- (15.65,1.65);

\foreach \x in {-1,0,...,31}
{
\draw [thick] ({14.95-\x*0.35},1.65) -- ({14.95-\x*0.35},2.4);
}
\foreach \x in {-2,0,2,4,6,8,10,12,14,16,18,20,22,24,26,28,30}
{
\draw [thick] ({14.95-\x*0.35},1.65) arc (0:-180:0.175);
}
\foreach \x in {-1,1,3,5,7,9,11,13,15,17,19,21,23,25,27,29}
{
\draw [thick] ({14.95-\x*0.35},2.4) arc (0:180:0.175);
}

\draw [thick] (4.1,2.4) arc (0:90:0.175);
\draw [thick] ({4.1-0.175},{2.4+0.175}) -- ({0.25+0.175},{2.4+0.175});
\draw [thick] ({0.25+0.175},{2.4+0.175}) arc (90:270:0.175);
\draw [thick] ({4.1-0.35},{2.4+0.175-0.35}) -- ({0.25+0.175},{2.4+0.175-0.35});
\draw [thick] ({4.1-0.35},{2.4+0.175-0.35}) arc (90:-90:0.175);
\draw [thick] ({4.1-0.35},{2.4+0.175-0.7}) -- ({0.25},{2.4+0.175-0.7});

\foreach \x in {0,...,20}
{
\draw ({0.25+\x*0.175},1.875) -- ({0.25+\x*0.175},1.525);
\draw [fill] ({0.25+\x*0.175},1.875) circle [radius=\vxrad];
\draw [fill] ({0.25+\x*0.175},1.525) circle [radius=\vxrad];
}
\end{tikzpicture}
\vspace{-0.5cm}
\caption{The tree $T_{n;r,b}$ which we use for our decomposition (see Definition~\ref{def:tree}), with the edge absorption structure and colour absorption structure highlighted as a) and b) above, with the vertex absorption structure marked with c). Note that, whichever our choice of one edge from each of the monochromatic matchings in a) and one edge from each of the rainbow matchings in b), the resulting tree is the same.}\label{illustration}
\end{figure}

\smallskip
\noindent\textbf{Edge absorbers via monochromatic matchings.}
We create an edge absorption structure for a set of reservoir edges as follows, where $\eta$ is a small constant.
(Recall that $t=n/2$ is the number of trees in our decomposition.)
For each $i \in [t]$, we construct a rainbow forest $\tilde{F}_i$
(where we will have $\tilde{F}_i \subseteq T_i$) and
matchings $M_{i,c}$. Each $M_{i,c}$ will consist of edges of colour $c$, and $c$
ranges over all elements of some colour set $D_i'$, where $|D_i'| \sim 6\eta n$. The matchings may overlap but are edge-disjoint from $\tilde{F}_i$, and, for each matching $M_{i,c}$,
any one of its edges can be added to $\tilde{F}_i$ to obtain a rainbow tree.
More precisely, we have the following `local' edge absorption property for each
$i \in  [t]$:
\begin{itemize}
\item[(P)] If one edge $e_{i,c}$ is chosen from each matching $M_{i,c}$, then
$F^+_i:=\tilde{F}_i+\sum_{c \in D_i'} e_{i,c}$ is a rainbow tree with vertex set $V(\tilde{F}_i)$.
\end{itemize}
Note that since the $M_{i,c}$ are monochromatic, the colour set of $F_i^+$ does not depend on the choice of $e_{i,c}$.
See Figure~\ref{illustration} for our construction of such a subgraph $\tilde{F}_i$
and the matchings $M_{i,c}$. We think of $M_{i,c}$ as being (the essential part of) an absorber which is able to `absorb' exactly one of the edges it contains. The chosen edge $e_{i,c}$ is then added
to $\tilde{F}_i$ to become part of the tree $T_i$.

Since the $M_{i,c}$ will be small (of size 256) and monochromatic, the requirement that exactly one edge from each  $M_{i,c}$ is to be added to $\tilde{F}_i$  is very restrictive. However, by carefully choosing how edges appear in different
matchings $M_{i,c}$, we can combine these to create the following `global' edge absorption property for two suitable subgraphs $G_1$ and $G_2$ of $K_n$.
(Here $G_1$, $G_2$ and the forests $\tilde{F}_i$, $i\in [t]$, will be edge-disjoint.)
\begin{itemize}
\item[(Q)] For any subset $E^\ast\In E(G_1)$ which consists of precisely $\eta n$ edges of each colour $c\in C$, we can label $E^\ast\cup E(G_2)$ as $\{e_{i,c}:i\in [t],c\in D_i'\}$ so that $e_{i,c}\in M_{i,c}$ for each $i\in[t]$ and $c\in D_i'$.
\end{itemize}
Properties (P) and (Q) mean that, given any set $E^\ast \subseteq E(G_1)$ with the right number of edges of each colour, we can absorb these edges
(along with those in the `buffer set' $E(G_2)$)
into the forests $\tilde{F}_1,\ldots,\tilde{F}_t$ to obtain rainbow trees which span some pre-determined vertex set and colour set (these sets are different for each tree). In fact, the equidistribution condition on $E^\ast$ will be naturally satisfied as the resulting trees must contain exactly one edge of each colour.
Thus altogether, the local edge absorption structures give rise to a \emph{global edge reservoir} (namely $E(G_1)$), for which we can absorb a leftover edge set $E^*$ into the existing forests.

To choose the matchings $M_{i,c}\subseteq E(G_1\cup G_2)$, we
consider a set of auxiliary graphs (called `robustly matchable bipartite graphs'), introduced in~\cite{montgomery:18} and already a standard technique in the construction of absorbers.
As the name suggests, these graphs have the property that one can find a perfect matching even after the removal of an arbitrary set of vertices (of given size)
from the larger vertex class, $B$ say.
We will consider one such robustly matchable bipartite graph $H_c=H_c[A,B]$
for each colour $c$, where $B=E_c(G_1 \cup G_2)$.
The neighbourhood in $H_c$ of each vertex $a \in A$ will correspond to some matching $M_{i,c}$, where $i$ is such that $c \in D_i'$. Thus adjacencies in $H_c$ encode the possible
absorber matchings $M_{i,c}$ (and thus the possible trees) that a reservoir edge $e \in B$ can be assigned to.
A matching in $H_c$ saturating $E^*_c \cup E_c(G_2)$ (where $E^*_c$ is the set of $c$-edges in $E^*$) gives an assignment of these `leftover' edges of colour $c$ to absorbers and thus to the trees $T_i$.
Carrying this out for all $c \in C$ allows us to absorb all the leftover edges $E^*$ from
the edge reservoir $G_1$ and the buffer edges (i.e.~those in $G_2$).

The robustly matchable graphs are discussed in more detail in Section~\ref{sec:RMBGs} and the properties of the edge absorption structure are
described in Lemma~\ref{lem:monochromatic matchings}.

\smallskip
\noindent\textbf{Colour absorbers via rainbow matchings.} The above properties allow us to use part of the edge reservoir
$G_1$ to create separate \emph{colour absorbers} for each tree. This means that for the $i$th tree we
have a reservoir $C'_{i,1}$ of colours with the property that any `leftover' (i.e.~so far unused)
set of colours $C^* \subseteq C'_{i,1}$ of given size can be absorbed into the
$i$th rainbow forest so that the result is still a (larger) rainbow forest.

More precisely, for the $i$th tree, we find a rainbow forest $\tilde{F}_i'$ which is vertex- and colour-disjoint from $\tilde{F}_i$, along with
small rainbow matchings $M'_{i,1},\ldots,M'_{i,3s}$ which are edge-disjoint from $\tilde{F}'_i$, as well as colour sets $C'_{i,1}$ (the `colour reservoir') and $C'_{i,2}$
(the `buffer set'), such that the following `local colour absorption' property holds for each $i\in [t]$ (where $s=\eta n/768$, and $|C'_{i,1}|=|C'_{i,2}|=2s$).
\begin{enumerate}[label = (P$'$)]
\item Given any set $C^\ast\In C'_{i,1}$ of $s$ colours, we can choose one edge $f_{i,j}$  from each $M'_{i,j}$ so that $\tilde{F}'_i+f_{i,1}+\ldots+f_{i,3s}$ is a $(\phi(E(\tilde{F}'_i))\cup C^\ast\cup C'_{i,2})$-rainbow tree with vertex set
$V(\tilde{F}'_i)$.\label{thispprime}
\end{enumerate}
For each colour $c$ appearing on an edge in $M'_{i,j}$, we think of $M'_{i,j}$ as (part of) an
absorber which can absorb colour $c$ into the $i$th tree (and for each $c$, we will provide several of these absorbers). The edges of the $M'_{i,j}$ will lie in the edge reservoir
$G_1$. Crucially, this means that when absorbing a colour $c$,
it does not matter which edges/absorbers are actually involved in this colour absorption step -- we can absorb any unused ones later. This means that the colour absorption step is less delicate than the edge absorption step.
See Figure~\ref{illustration} for our construction of such an $\tilde{F}'_i$.
(In the main proof, we actually construct the forests $\tilde{F}_i$ and $\tilde{F}'_i$
simultaneously, and denote them  $\tilde{F}_i$.)

The matchings $M'_{i,j}$ will be small edge-disjoint rainbow matchings, where the colours of each matching $M'_{i,j}$ are chosen according to some auxiliary robustly matchable bipartite graph. We will consider one such auxiliary graph $H_i$ for each tree $T_i$, with the larger vertex
class consisting of the colour reservoir $C'_{i,1}$ together with the buffer set $C'_{i,2}$.
The edges of $H_i$ connect each colour $c$ to some indices $j \in [3s]$.
The colour set of $M'_{i,j}$ will consist of precisely those colours in $N_{H_i}(j)$.
For any set $C^* \subseteq C'_{i,1}$ of size $s$, a matching saturating
$C^* \cup C'_{i,2}$ absorbs all the `leftover' colours, as required.
The details are given in Lemma~\ref{cor:rainbow matchings}.


\smallskip
\noindent\textbf{Vertex absorbers.} We then use part of both the edge reservoir and the colour reservoir to create \emph{vertex absorbers}. This construction is relatively simple, and the resulting vertex reservoir consists of some vertices unused by the $i$th tree so far.
For each $i \in [t]$, we take a small random set $A_i$ of vertices and connect them into a rainbow \emph{vertex absorbing path}, while reserving a further random set of vertices $B_i$ that is slightly smaller than $A_i$.
When we reach Step (8), the set of uncovered vertices will be a subset of $B_i$
and contain almost all vertices of~$B_i$. (So one can view $B_i$ as a vertex reservoir.)
We will match those vertices in $B_i$ which are still uncovered
onto the vertex absorbing path.
 The randomness of $A_i$ and $B_i$ allows us to do this with a rainbow matching
 between $A_i$ and~$B_i$.

\smallskip
\noindent\textbf{Covering outside the reservoirs.} By construction, the edge and colour absorbing structures can only deal with edges/colours within the respective reservoirs. Thus, after we construct the $i$th forest $F_i$%
\COMMENT{actually, in the main proof this graph is more or less $F_i^*$, but
if we change to $F_i^*$ here, we might also want to change (6) too}
 which covers almost all the colours, we must extend it slightly so that it now uses every colour outside its reservoir, and that collectively the resulting forests use all the edges outside of the global edge reservoir.
We achieve this as follows:
To cover an edge $e$ outside the global edge reservoir (in Step (6)),
we include $e$ as an edge between $A_i$ and $B_i$
for some suitable~$i$.
Similarly, to cover a colour $c$ outside the $i$th colour reservoir (in Step (7)), we choose a suitable $c$-edge $e$ between $A_i$ and $B_i$, again from the edge reservoir.
We can carry this out in such a way that these edges form a relatively small $A_iB_i$-matching, thus enabling us to carry out the vertex absorption procedure
described above with only minor modifications.

\subsection{Almost-packing random subgraphs}\label{sec:hypermatch}
We will find the different structures for the strategy outlined above by defining
(for each of these structures) an
auxiliary hypergraph in which a large matching corresponds to the desired structure. The hypergraph will be roughly regular, with small codegrees, and thus the existence of this matching will follow from standard results
(see Theorem~\ref{thm:AY} in  Section~\ref{sec:nibble}). In each case, the hypergraph is defined in a similar way, but to give a concrete example we will sketch how to find $t=n/2$ almost-spanning rainbow paths in any optimally coloured $K_n$.
(Note that this construction as described below is already present in~\cite{KKKO:18}.
We repeat it informally here, as it forms a template for several more involved
applications in this paper.)

To simplify further, we note that by randomly reserving edges, colours and vertices, we can greedily connect a given set of long disjoint rainbow paths together via very short paths
(which use their own set of reserved edges, colours and vertices)
into a single rainbow path. Thus, to cut to the main part of the argument, let us suppose we want to find the following, where $\ell$ is a large constant, and $r\ell\leq (1-\eps)n$, for some small $\eps>0$.

\begin{itemize}
\item[Aim: ] To find in $K_n$, for each $i\in [t]$, a set $\cF_i$ of $r$ vertex-disjoint colour-disjoint rainbow cycles of length $\ell$, so that all the cycles in $\bigcup_{i\in [t]}\cF_i$ are edge-disjoint.
\end{itemize}

The key is to construct a hypergraph $\cH$ in which a large matching corresponds to the required cycles (where each matching edge directs us to include some cycle into some set $\cF_i$). We take vertices for $\cH$ as follows. We need all the cycles we find to be edge-disjoint, so each edge in $G$ will appear as a vertex of $\cH$. All the cycles in $\cF_i$ must be vertex-disjoint, so we wish to represent the vertices of $V=V(K_n)$ as vertices in~$\cH$. However, different cycles in different sets $\cF_i$ are permitted to share vertices. Thus, for each $i\in [t]$, we include a copy of $V$ by including the vertices in $\{i\}\times V$ in~$V(\cH)$. Similarly, we represent the colours for cycles by including $\{i\}\times C$ for each $i \in [t]$.
We define the hyperedges of $\cH$ as follows. For each rainbow cycle $F\subseteq K_n$ of length $\ell$ and $i\in [t]$, we include the hyperedge
\[
f_{i,F}:=E(F)\cup (\{i\}\times V(F))\cup (\{i\}\times \phi(E(F))).
\]
Suppose then we had a matching $\cM$ in $\cH$. Then, for each $i\in [t]$, let $\cF_i$ be the set of cycles $F$ with $f_{i,F}\in \cM$. Note that we have the following properties.
\begin{itemize}
\item If $F,L\in \cF_i$ are distinct, we have the following.
\begin{itemize}
	\item As $\{i\}\times V(F)\subseteq f_{i,F}$, $\{i\}\times V(L) \subseteq f_{i,L}$ and $f_{i,F},f_{i,L}\in \cM$, we have that $F$ and $L$ are vertex-disjoint.
\item As $\{i\}\times \phi(E(F))\subseteq f_{i,F}$, $\{i\}\times \phi(E(L))\subseteq f_{i,L}$ and $f_{i,F},f_{i,L}\in \cM$, we have that $F$ and $L$ are colour-disjoint.
\end{itemize}
\item For any $F\in \cF_i$ and $L\in \cF_j$ with $i\neq j$, we have $E(F)\subseteq f_{i,F}$, $E(L)\subseteq f_{j,L}$, and $f_{i,F},f_{j,L}\in \cM$, so $F$ and $L$ are edge-disjoint.
\end{itemize}
That is, each $\cF_i$ is a set of vertex- and colour-disjoint rainbow cycles, and the cycles in $\bigcup_{i\in [t]}\cF_i$ are edge-disjoint, as required in the above aim.

In the actual proof we will find the required structures within prescribed
(randomly chosen) vertex, edge and colour sets, with parameters carefully
chosen so that the construction uses almost all of the available sets each time.
Together, this has the advantage that the overall leftover after removing these
structures is also randomly distributed and sufficiently small so that it can be
absorbed.

\subsection{Isomorphic trees}\label{sec:iso}
The main achievement of our techniques is to find a decomposition into (any) spanning rainbow trees. However, by taking care at several points in our argument, the trees we construct can be kept isomorphic. The key point here is to observe that in Figure~\ref{illustration} the resulting structure from the absorber is the same regardless of which edges are used from the reservoir.

In fact, we not only find isomorphic trees, but we find copies of the same fixed tree, regardless of the $1$-factorization of $K_n$. We define this tree as follows (see Figure~\ref{illustration}).

\begin{defin}\label{def:tree}
Given $n,r,b\in \bN$ such that $\ell:=n-1020r-b-1> r + b$, we define the tree $T_{n;r,b}$ as follows:
Take a path $v_0 \dots  v_{\ell}$ of length~$\ell$. For all $k\in[r-1]$, add $510$ paths of length~$2$ to $v_{5k}$ (i.e.~$v_{5k}$ will become an endvertex of these $510$~paths), and add $255$ paths of length~$2$ to each of $v_0$ and $v_{5r}$. Moreover, take a set $B$ of $b$ new vertices and add a perfect matching between $B$ and $\Set{v_{\ell-b+1},\dots,v_\ell}$.\COMMENT{The number of vertices other than $v_0,\dots,v_{\ell}$ is $b+ (r-1)\cdot 510\cdot 2 + 4\cdot 255 = b+1020r$ }
\end{defin}

The set $B$ corresponds to the set $B_i$ in the vertex absorption structure.
For each $i \in [t]$ there will be an integer $r_i \le r$ so that
for each $k \in [r_i]$, the `middle' edge on the path between $v_{5(k-1)}$ and $v_{5k}$
will be an edge of some `absorber-matching' $M_{i,c}$ or $M'_{i,j}$.
Note that $|T_{n;r,b}|=n$ and $\Delta(T_{n;r,b}) \le 512$.
We will prove Theorem~\ref{thm:main} with $T=T_{n;r,b}$, where $r$ and $b$ are small but linear in~$n$. So $T$ contains an almost spanning path. After proving Theorem~\ref{thm:main} in Section~\ref{sec:main}, we describe how this construction can be slightly modified to achieve that $\Delta(T)=3$.

\section{Tools}\label{sec:tools}

\subsection{Hypergraph matchings}\label{sec:nibble}
We make frequent use of the existence of large matchings in almost regular hypergraphs with small codegrees (such matchings are constructed via semi-random nibble methods pioneered by R\"odl~\cite{rodl:85}). Moreover, we wish to have a matching which is `well-distributed' across a number of vertex subsets. To make this precise, we use the following definition.

\begin{defin} \label{gammaperfect}
Given a hypergraph $\cH$ and a collection $\cF$ of subsets of $V(\cH)$, we say a matching $\cM$ in $\cH$ is \defn{$(\gamma,\cF)$-perfect} if for each $F\in \cF$, at most $\gamma\cdot \max\Set{|F|,|V(\cH)|^{2/5}}$ vertices of $F$ are left uncovered by~$\cM$.
\end{defin}

Pippenger and Spencer~\cite{PS:89}  showed that in an almost regular hypergraph
$\cH$ with small codegrees there are many large edge-disjoint matchings. Alon and Yuster~\cite{AY:05} observed that by randomly splitting $V(\cH)$ into many parts,
and applying the Pippenger--Spencer theorem to each induced subhypergraph and then selecting a matching in each of these subhypergraphs at random, one can obtain an almost perfect matching $\cM$
of $\cH$ that is `well-distributed' in the sense of Definition~\ref{gammaperfect}. We will use the following consequence of Theorem 1.2 in~\cite{AY:05}.

\begin{theorem}[\cite{AY:05}] \label{thm:AY}
Suppose $1/n \ll \eps \ll \gamma,1/r$. Let $\cH$ be an $r$-uniform hypergraph on $n$ vertices such that for some $D\in \bN$, we have $d_{\cH}(x)=(1\pm \eps)D$ for all $x\in V(\cH)$, and $\Delta^c(\cH)\le D/\log^{9r} n$. Suppose that $\cF$ is a collection of subsets of $V(\cH)$ such that $|\cF|\le n^{\log n}$.
Then there exists a $(\gamma,\cF)$-perfect matching in~$\cH$.
\end{theorem}
\COMMENT{original \cite{AY:05}: for all $r\in \bN,C>1,\eps>0$ there exist $\mu,K>0$ such that for any $r$-graph $H$ on $n$ vertices with $\delta(H)\ge (1-\mu)\Delta(H)$, $g(H)> \max\Set{1/\mu,K(\ln n)^6}$ and $|\cF|\le C^{g(H)^{1/(3r-3)}}$ and $\min\set{|F|}{F\in \cF}\ge 5g(H)^{1/(3r-3)}\ln(|\cF|g(H))$, where $g(H):=\Delta(H)/\Delta^c(\cH)$, there exists a matching $\cM$ in $\cH$ such that $|F\sm V(\cM)|\le \eps |F|$ for all $F\in \cF$. Reduction: Apply with $C=2$, say, and $\eps=\gamma$, and choose $\eps$ such that $\frac{1-\eps}{1+\eps}\ge 1-\mu$. With $\Delta^c(\cH)\ge 1$ and $D\le n^{r-1}$, we have that $g(H)^{1/(3r-3)}\le n^{1/3}$. Moreover, by assumption, $g(H)^{1/(3r-3)}\ge \log^{3}n$. Thus, we could apply \cite{AY:05} if all sets in $\cF$ have size at least $n^{2/5}$. For all sets $F\in \cF$ which are too small, we add dummy vertices to bring them to size $n^{2/5}$. The result then carries over.}

We apply Theorem~\ref{thm:AY} to several different hypergraphs in our proof, each time checking the appropriate degree and codegree bounds. We comment here generally why our hypergraphs are almost regular with small codegree. Indeed, roughly speaking, in each hypergraph $\cH$ we define (see Section~\ref{sec:hypermatch}), estimating vertex degrees will correspond to counting the number of rainbow copies of a certain graph in $K_n$ with one fixed characteristic
(e.g.~one fixed vertex/edge/colour). The symmetry in our choice of random subsets and subgraphs will mean that for each characteristic, the vertex degrees in $\cH$ are roughly the same. Our choice of edge, colour and vertex probabilities then results in an almost regular hypergraph. (Here, it is also useful that we consider $1$-factorizations
rather than proper colourings.) Counting codegrees corresponds roughly to counting the number of copies of the same subgraph but with two characteristics fixed. This means that the codegrees are small in comparison to the degrees, giving the additional condition we need to apply Theorem~\ref{thm:AY}.

\subsection{Robustly matchable bipartite graphs}\label{sec:RMBGs}
As noted in Section~\ref{sec:absorb}, we use robustly matchable bipartite graphs as auxiliary graphs to tell us how to distribute edges during the absorbing steps. These graphs are defined as follows.

\begin{defin}
Given pairwise disjoint vertex sets $X,Y,Z$, an $\mathrm{RMBG}(X,Y,Z)$ 
is a bipartite graph $H$ with bipartition $(X,Y \cup Z)$ and the following crucial property: for any set $Y'\In Y$ with $|Y'|=|X|-|Z|$, the subgraph $H[X,Y'\cup Z]$ has a perfect matching.

We also refer to $H$ as an $\mathrm{RMBG}(|X|,|Y|,|Z|)$ with parts~$X,Y,Z$.
\end{defin}

Such graphs were introduced in~\cite{montgomery:18} in order to find given spanning trees in random graphs.

\begin{lemma}[{\cite[Lemma~10.7]{montgomery:18}}]\label{lem:rmbg original}
For all sufficiently large~$m$, there exists an $\mathrm{RMBG}(3m,2m,2m)$ with maximum degree at most~$100$.
\end{lemma}

We say that a bipartite graph $H$ with bipartition $(X,Y)$ is \defn{$(\ell,r)$-regular} if all the vertices in $X$ have degree $\ell$ and all the vertices in $Y$ have degree~$r$.
Using the Max-Flow-Min-Cut-theorem, it is straightforward to find a supergraph of an RMBG from Lemma~\ref{lem:rmbg original} which is appropriately regular.

\begin{cor} \label{cor:regular rmbg}
For all fixed $d\ge 59$ and sufficiently large~$m$, there exists a $(4d,3d)$-regular $\mathrm{RMBG}(3m,2m,2m)$.
\end{cor}

\proof
Let $H$ be an $\mathrm{RMBG}(3m,2m,2m)$ with parts $X$, $Y$ and $Z$ and maximum degree at most $100$, as in Lemma~\ref{lem:rmbg original}. Take new vertices $s,t$ and let $G$ be the directed graph obtained from the complete bipartite graph between $X$ and $Y\cup Z$ (with all edges directed towards $Y\cup Z$) by removing the edges of $H$ and adding all edges from $s$ to $X$ and from $Y\cup Z$ to~$t$.
An edge $sx$ receives capacity $4d-d_H(x)$, and an edge $yt$ receives capacity $3d-d_H(y)$. All edges in $G[X,Y\cup Z]$ receive capacity~$1$. We claim that $(\Set{s},V(G)\sm \Set{s})$ and $(\Set{t},V(G)\sm \Set{t})$ are minimal $(s,t)$-cuts. Indeed, first note that the capacity of these cuts is $12dm-e(H)$. Now, let $(S,T)$ be any $(s,t)$-cut. Let $S_1:=S\cap X$, $S_2:=S\cap (Y\cup Z)$, $T_1:=T\cap X$ and $T_2:=T\cap (Y\cup Z)$. The capacity $c(S,T)$ of the cut $(S,T)$ satisfies
\begin{align}
	c(S,T) &= \sum_{x\in T_1}(4d-d_H(x)) + |S_1||T_2| - e_H(S_1,T_2) + \sum_{y\in S_2}(3d-d_H(y))\label{hello}
	\\ &= 12dm+|S_1|(|T_2|-4d)+3d|S_2|-e(H)-e_H(T_1,S_2)\nonumber
	\\ &= 12dm+4d|T_1|+|T_2|(|S_1|-3d)-e(H)-e_H(T_1,S_2).\nonumber
\end{align}\COMMENT{where we use that $\sum_{x\in T_1}d_H(x) +  \sum_{y\in S_2}d_H(y) + e_H(S_1,T_2) = e(H)+e_H(T_1,S_2)$ }
Thus, if $|S_1|\ge 3d$ or $|T_2|\ge 4d$, then $c(S,T)\ge 12dm-e(H)$, as desired.\COMMENT{since $e_H(T_1,S_2)\le 100|S_2|,100|T_1|$} Therefore, assume that $|S_1|< 3d$ and $|T_2|< 4d$. Then, \eqref{hello} implies that $c(S,T)\ge (3m-3d)(4d-100) + (4m-4d)(3d-100) \ge 12md$, where the last inequality uses $d\ge 59$. This proves the claim.

By the Max-Flow-Min-Cut-theorem, there exists an $(s,t)$-flow in $G$ with value $12dm-e(H)$. This yields a subgraph $H'\In G\sm \Set{s,t}$ such that $d_{H'}(x)=4d-d_H(x)$ for all $x\in X$ and $d_{H'}(y)=3d-d_H(y)$ for all $y\in Y\cup Z$. Thus, $H\cup H'$ is the desired $(4d,3d)$-regular $\mathrm{RMBG}(3m,2m,2m)$.
\endproof

\subsection{Probabilistic tools}\label{sec:prob}
In order to show various properties of random subgraphs and subsets, we will use common concentration inequalities, as follows.

\begin{lemma}[see {\cite[{Corollary~2.3, Corollary~2.4 and Theorem 2.8}]{JLR:00}}] \label{lem:chernoff}
Let $X$ be the sum of $n$ independent Bernoulli random variables. Then the following hold.
\begin{enumerate}[label={\rm(\roman*)}]
\item For all $0\le\eps \le 3/2$, $\prob{|X - \expn{X}| \geq \eps\expn{X} } \leq 2\eul^{-\eps^2\expn{X}/3}$.\label{chernoff eps}
\item If $t\ge 7 \expn{X}$, then $\prob{X\ge t}\le \eul^{-t}$.\label{chernoff crude}
\end{enumerate}
\end{lemma}
Throughout, we will refer to \ref{chernoff eps} as `Chernoff's bound'. Often, we will use this in conjunction with an implicit union bound to show that several properties hold altogether with high probability.

\begin{fact}[cf.~{\cite[Lemma 8]{raman:90}}] \label{prop:jain}
Let $X_1, \ldots, X_n$ be Bernoulli random variables such that for all $i \in [n]$, we have $\prob{X_i = 1 \mid X_1, \ldots, X_{i-1}} \leq p$.
Let $B \sim Bin(n,p)$ and $X:=\sum_{i=1}^n X_i$. Then $\prob{X \geq t} \leq \prob{B \geq t}$ for all~$t\ge 0$.
\end{fact}

One important tool to prove concentration of our random variables is McDiarmid's inequality.
Let $X_1,\dots,X_m$ be independent random variables taking values in~$\cX$. Let $f\colon \cX^m\to \bR$ be a function of $X_1,\dots,X_m$ such that $$|f(x_1,\dots,x_{i-1},x_i,x_{i+1},\dots,x_m)-f(x_1,\dots,x_{i-1},x_i',x_{i+1},\dots,x_m)|\le c_i$$ for all $i\in [m]$, $x_i',x_1,\dots,x_m\in \cX$. If this holds, we say that $X_i$ \defn{affects} $f$ by at most~$c_i$.

\begin{theorem}[McDiarmid's inequality, see~{\cite[Lemma~1.2]{mcdiarmid:89}}] \label{thm:McDiarmid} Let $X_1,\ldots, X_m$, $c_1,\ldots,c_m$ and $f$ be as stated above. Then, for all $t>0$, $$\prob{|f(X_1,\ldots,X_m)-\expn{f(X_1,\ldots,X_m)}|\ge t} \le 2 \eul^{-\frac{2t^2}{\sum_{i=1}^m c_i^2}}.$$
\end{theorem}

For our purposes, we will have $\cX=\Set{0,1}$, and the $X_i$ will be indicator variables of certain events. We will often use different indicator variables (which are not necessarily independent) to compute~$\expn{f}$, and then use McDiarmid's inequality to prove concentration.

\section{Proof}

\subsection{Approximate decomposition}\label{sec:approx}
The main result in this Section is Lemma~\ref{lem:approximate dec},
which implies the existence of an approximate decomposition into rainbow almost spanning paths for any given $1$-factorization of $K_n$. As noted earlier, this result was already proved in~\cite{KKKO:18,MPS:18b}.

However, we need to strengthen the result somewhat -- in particular, we need to
constrain the paths to use given (randomly chosen) vertex and colour sets, and we
need the paths to be well-behaved towards given subsets of these sets
(the latter is encapsulated in the concept of `boundedness' defined below).

The proof of Lemma~\ref{lem:approximate dec} relies on ideas from Theorem~1.5 and Lemma~2.14 in~\cite{KKKO:18} (and simplifies some aspects of that argument).
As described in Section~\ref{sec:nibble}, the strategy is to first find
for each $i\in[t]$ an almost spanning collection of vertex-disjoint long rainbow cycles. Then we delete an edge from each such cycle and connect them into a long rainbow path via Lemma~\ref{lem:greedy connections}.

\begin{defin}[$m$-bounded]\label{boundeddefn}
Let $\phi$ be a $1$-factorization of the complete graph~$K_n$ with vertex set $V$ and colour set~$C$. Given a subgraph $G\In K_n$, vertex sets $\Set{V_i}_{i\in[t]}\In V$ and colour sets $\Set{C_i}_{i\in[t]}\In C$, we say that $(G,\Set{V_i}_{i\in[t]},\Set{C_i}_{i\in[t]})$ is \defn{$m$-bounded} if the following hold:
\begin{enumerate}[label=\rm{(B\arabic*)}]
\item for all $i\in[t]$, $|V_i|,|C_i|\le m$;\label{bounded:set sizes}
\item for all $v\in V$, we have $|\set{i\in[t]}{v\in V_i}|\le m$ and $d_G(v)\le m$;\label{bounded:vertex incidences}
\item for all $c\in C$, $|\set{i\in[t]}{c\in C_i}|\le m$ and $|E_c(G)|\le m$.\label{bounded:colour incidences}
\end{enumerate}
Here, we think of $(G,\Set{V_i}_{i\in[t]},\Set{C_i}_{i\in[t]})$ as being `leftovers' that we want to be `well-behaved' in the above sense.
\end{defin}

The following lemma allows us to embed rooted graphs in a rainbow fashion. We will often apply it to find the desired rainbow subgraphs that were missed by an application of Theorem~\ref{thm:AY}.

\begin{lemma} \label{lem:greedy connections}
Suppose $1/n\ll \gamma \ll 1/\Delta$ and $t\le n$. Let $G$ be a graph on $n$ vertices and $\phi\colon E(G)\to C$ a proper edge colouring. Moreover, let $V_1,\dots,V_t\In V(G)$ and $C_1,\dots,C_t\In C$ be such that for each $i\in[t]$ and any set $S\In V(G)$ with $|S|\le \Delta$, we have that
\begin{align}
	|\set{v\in N_G(S)\cap V_i}{\phi(uv)\in C_i \mbox{ for each }u\in S}| \ge \gamma^{1/3}n.\label{clock}
\end{align}

For each $i\in[t]$, let $H_i$ be a graph with $|V(H_i)|,|E(H_i)|\le \gamma n$ and $\Delta(H_i)\le \Delta$,\COMMENT{Could also use the degeneracy of $H_i$ rooted at $X_i$ here, which gives more efficient lemma, but at expense of another definition} and let $\Lambda_i\colon X_i \to V(G)$ be an injection,\COMMENT{most of the time $X_i$ will be a subset of $V(G)$ and $\Lambda_i$ the identity} where $X_i\In V(H_i)$ is independent. Assume that for all $u\in V(G)$, there are at most $\gamma n$ indices $i\in[t]$ for which $u\in \Ima(\Lambda_i)$.

Then, there exist embeddings $\psi_i\colon H_i\to G$, $i\in [t]$, such that,  for each $i\in[t]$,  $\psi_i(H_i)$ is rainbow with colours in~$C_i$, $\psi_i(x)=\Lambda_i(x)$ for all $x\in X_i$ and $\psi_i(x)\in V_i$ for all $x\in V(H_i)\sm X_i$, and such that $\psi_1(H_1),\dots,\psi_t(H_t)$ are edge-disjoint.
\end{lemma}

\proof
We find the embeddings $\psi_1,\dots,\psi_t$ successively and greedily. For $s\in[t]$ and a vertex $u\in V(G)$, let $r(u,s)$ be the number of indices $i\in[s]$ for which $u\in \Ima(\Lambda_i)$. By assumption, $r(u,s)\le \gamma n$.

Suppose that, for some $s\in[t]$, we have already found suitable embeddings $\psi_1,\dots,\psi_{s-1}$ such that, additionally,
\begin{align}
	\mbox{the degree of each $u\in V$ in $\bigcup_{i\in[s-1]}\psi_i(H_i)$ is at most $2\sqrt{\gamma} n + r(u,s-1)\Delta$.} \label{clever greedy degree}
\end{align}

Now, we find a suitable embedding $\psi_s$ such that \eqref{clever greedy degree} holds with $s$ replaced by~$s+1$.
Let $B$ be the set of all vertices whose degree in $\bigcup_{i\in[s-1]}\psi_i(H_i)$ is larger than $\sqrt{\gamma} n$. Since $\sum_{i\in[s-1]}|E(H_i)|\le \gamma n^2$, observe that $|B|\le 2\sqrt{\gamma}n$.

We can now greedily embed $H_s$ while avoiding the vertices in $B\setminus \Lambda_s(X_s)$. For all $x\in X_s$, define $\psi_s(x)=\Lambda_s(x)$. Order the remaining vertices of $H_s$ arbitrarily and embed them one by one into~$V(G)\sm B$ as follows. When we consider~$x\in V(H_s)\sm X_s$, let $S$ be the set of images of the neighbours of $x$ which have already been embedded. We would like to choose an image for $x$ from $Y:=\set{v\in N_G(S)\cap V_s}{\phi(uv)\in C_s \mbox{ for each }u\in S} \sm B$.

First, note that by~\eqref{clock} with $i=s$, and as $|B|\le 2\sqrt{\gamma}n$, we have $|Y|\ge \gamma^{1/3}n -2\sqrt{\gamma}n$.
At most $\gamma n$ vertices of $Y$ are blocked because they have already been chosen as an image for $H_s$, and at most $|S|\gamma n\le \Delta \gamma n$ vertices $v\in Y$ are blocked because $\set{\phi(uv)}{u\in S}$ contains a colour that has already been used. Moreover, invoking~\eqref{clever greedy degree}, at most $\Delta (2\sqrt{\gamma}+\Delta\gamma)n$ vertices $v$ are blocked because $uv\in E(\bigcup_{i\in[s-1]}\psi_i(H_i))$ for some $u\in S$. Hence, there exists a suitable image for~$x$. Thus, we can finish the embedding of $H_s$ in this way.

Clearly, $\psi_s(H_s)$ is rainbow with colours in~$C_s$, and edge-disjoint from $\psi_1(H_1),\dots,\psi_{s-1}(H_{s-1})$. Moreover, for any vertex $u\in V(G)$, the degree of $u$ in $\psi_s(H_s)$ is at most $\Delta$. In particular, if $u\notin B$, then the degree of $u$ in $\bigcup_{i\in[s]}\psi_i(H_i)$ is at most $\sqrt{\gamma} n+\Delta\le 2\sqrt{\gamma} n$. Moreover, if $u\in B$, then $u\in V(\psi_s(H_s))$ if and only if $u\in \Ima(\Lambda_s)$, in which case the degree of $u$ in $\psi_s(H_s)$ is at most $\Delta=(r(u,s)-r(u,s-1))\Delta$. Thus, \eqref{clever greedy degree} holds with $s$ replaced by~$s+1$, which completes the proof.
\endproof

\begin{lemma}\label{lem:approximate dec}
Suppose $1/n \ll \gamma,\kappa \ll p$ and let $q:=\beta:=p$ and $t:=n/2$. Let $\phi$ be a $1$-factorization of the complete graph~$K_n$ with vertex set $V$ and colour set~$C$. For every $i\in[t]$, let $V_i$ be a $p$-random subset of~$V$, and let $C_i$ be a $q$-random subset of~$C$. Moreover, let $G$ be a $\beta$-random subgraph of~$K_n$.

Then with high probability, there exist edge-disjoint rainbow paths $P_1,\dots,P_t$ in $G$ such that
\begin{enumerate}[label=\rm{(P\arabic*)}]
\item $V(P_i)\In V_i$, $\phi(E(P_i))\In C_i$;
\item $(G-\bigcup_{i\in[t]}P_i$, $\Set{V_i\sm V(P_i)}_{i\in[t]}$, $\Set{C_i \sm \phi(E(P_i))}_{i\in[t]})$ is $5\gamma n$-bounded;\label{approx dec bounded leftover}
\item for all $v\in V$, the number of $i\in[t]$ for which $v\in V(P_i)$ and the subpath from $v$ to one of the endvertices of $P_i$ has length at most $\kappa n$, is at most $21\kappa p^{-1} n$.\label{path ends}
\end{enumerate}
\end{lemma}

\proof
Choose new constants $\eps>0$ and $\ell\in \bN$ such that $1/n \ll \eps \ll 1/\ell \ll \gamma, \kappa\ll p$.

For each $i\in[t]$, split $V_i$ into a $p(1-\gamma)$-random subset $V_{i,1}$ and a $p\gamma$-random subset $V_{i,2}$, and split $C_i$ into a $q(1-\gamma)$-random subset $C_{i,1}$ and a $q\gamma$-random subset $C_{i,2}$. We also split $G$ into a $\beta(1-\gamma)$-random subgraph $G_1$ and a $\beta \gamma$-random subgraph~$G_2$.
We claim that with high probability, we have the following:
\begin{enumerate}[label=\rm{(\alph*)}]
\item for each $i\in[t]$ and $S\In V$ with $|S|\le 2$,\label{typicality condition} $$|\set{v\in N_{G_1}(S)\cap V_{i,1}}{\phi(uv)\in C_{i,1} \mbox{ for each }u\in S}| = (1\pm \eps)(\beta q(1-\gamma)^2)^{|S|}p(1-\gamma)n;$$
\item for each $i\in[t]$ and $S\In V$ with $|S|\le 2$,\label{greedy connection condition} $$|\set{v\in N_{G_2}(S)\cap V_{i,2}}{\phi(uv)\in C_{i,2} \mbox{ for each }u\in S}| \ge \gamma^6 n;$$
\item for each $e\in E(K_n)$, $|\set{i\in[t]}{e\In V_{i,1},\phi(e)\in C_{i,1}}| = (1\pm \eps)p^2q(1-\gamma)^3 n/2$;\label{approx edge incidences}
\item for each $i\in[t]$ and $c\in C$, we have that $|E_c(G_1[V_{i,1}])| = (1\pm \eps)p^2\beta(1-\gamma)^3 n/2$;\label{approx colour incidences}
\item for each $i\in[t]$, $|V_i|=(1\pm \eps)pn$;\label{approx vertex size}
\item $(G_2,\Set{V_{i,2}}_{i\in[t]},\Set{C_{i,2}}_{i\in[t]})$ is $2\gamma n$-bounded.\label{approx reserve bounded}
\end{enumerate}
Indeed, using Chernoff's bound, it is straightforward to check that \ref{approx edge incidences}--\ref{approx reserve bounded} hold with high probability. For \ref{typicality condition} and \ref{greedy connection condition}, we can apply McDiarmid's inequality, since of the at most $|S|n$ edges incident with $S$, each has an effect of at most~$1$, each vertex has an effect of at most~$1$, and each colour has an effect of at most~$|S|$.

Henceforth, assume that these random choices have been made and satisfy the above properties.

For all $i\in[t]$, let $\O_i$ be the collection of all rainbow cycles $R$ of length~$\ell$ in $G_1$ for which $V(R)\In V_{i,1}$ and $\phi(E(R))\In C_{i,1}$. (Note that the $\cO_i$'s are not necessarily disjoint.)
For $v\in V$, $c\in C$ and $e\in E(K_n)$, we let $\cO_i(v)$, $\cO_i(c)$ and $\cO_i(e)$ denote the set of all $R\in \cO_i$ with $v\in V(R)$, $c\in \phi(E(R))$ and $e\in E(R)$, respectively.
Using \ref{typicality condition}, we can now count that, for all $v\in V_{i,1}$,
\begin{align}
	|\cO_i(v)|&=\frac{1}{2}\cdot ((1\pm 2\eps)(1-\gamma)^3\beta q p n)^{\ell-2} \cdot (1\pm 2\eps)(1-\gamma)^5(\beta q)^2 pn\nonumber \\
	         &= (1\pm \sqrt{\eps})\frac{1}{2}(1-\gamma)^{3\ell-1}\beta^\ell q^\ell p^{\ell-1} n^{\ell-1}, \label{usefulvx}
\end{align}\COMMENT{$\frac{1}{2}$ for orientation of cycle, $2\eps$ instead of $\eps$ because can't use vertex or colour again}
 and, for all $e\in E(G_1[V_{i,1}])$ with $\phi(e)\in C_{i,1}$,
\begin{align}
	|\cO_i(e)| &= ((1\pm 2\eps)(1-\gamma)^3\beta q p n)^{\ell-3} \cdot (1\pm 2\eps)(1-\gamma)^5(\beta q)^2 pn \nonumber\\
	    &= (1\pm \sqrt{\eps})(1-\gamma)^{3\ell-4} \beta^{\ell-1} q^{\ell-1} p^{\ell-2} n^{\ell-2}.\label{usefuledge}
\end{align}

We define an auxiliary hypergraph $\cH$ as follows. The vertex set of $\cH$ consists of three parts. The first part is simply~$E(G_1)$. The second part is the set $\cV$ of all pairs $(i,v)$ with $i\in [t]$ and $v\in V_{i,1}$. The third part is the set $\cC$ of all pairs $(i,c)$ with $i\in [t]$ and $c\in C_{i,1}$.

Now, we define the edge set of~$\cH$.
For each $i\in[t]$ and $R\in \cO_i$, we add the hyperedge
\begin{align}
f_{i,R}:= E(R) \cup (\Set{i}\times V(R)) \cup (\Set{i}\times \phi(E(R))).\label{hyperedge}
\end{align}
Hence, $\cH$ is $3\ell$-uniform.

Clearly, using \eqref{usefulvx}, we have for each $(i,v)\in \cV$ that
\begin{align}\label{reg1}
	d_\cH((i,v))= |\cO_i(v)| = (1\pm \sqrt{\eps})\frac{1}{2}(1-\gamma)^{3\ell-1}\beta^\ell q^\ell p^{\ell-1} n^{\ell-1}.
	\end{align}
	Moreover, we have for each $e\in E(G_1)$ that
\begin{align}\label{reg2}
	d_\cH(e) = \sum_{i\in[t]\colon e\In V_{i,1},\phi(e)\in C_{i,1}}|\cO_i(e)|
	\overset{\ref{approx edge incidences},\eqref{usefuledge}}{=} (1\pm 2\sqrt{\eps}) \frac{1}{2} (1-\gamma)^{3\ell-1} \beta^{\ell-1} q^{\ell} p^{\ell} n^{\ell-1}
\end{align}
and for all $(i,c)\in \cC$ that
\begin{align}\label{reg3}
	d_\cH((i,c))= \sum_{e\in E_c(G_1[V_{i,1}])} |\cO_i(e)|
	\overset{\ref{approx colour incidences},\eqref{usefuledge}}{=} (1\pm 2\sqrt{\eps}) \frac{1}{2} (1-\gamma)^{3\ell-1} \beta^{\ell} q^{\ell-1} p^{\ell} n^{\ell-1}.
\end{align}
(Note that no hyperedge is counted more than once since each rainbow cycle contains at most one $c$-edge.)

\begin{NoHyper}
\begin{claim}\label{claim:aux codegree approx}
$\Delta^c(\cH)\le \ell^4 n^{\ell-2}$.
\end{claim}
\end{NoHyper}

\claimproof Recall, from \eqref{hyperedge}, that each hyperedge of $\cH$ is uniquely fixed by some $i\in [t]$ and $R\in \cO_i$.
Note that for a set $S$ of vertices and $i\in[t]$, the number of $R\in \cO_i$ with $S\In V(R)$ is at most $\ell^{|S|}n^{\ell-|S|}$.
This easily implies that codegrees of pairs in $E(G_1)\times \cV$, $E(G_1)\times \cC$ and $\cV\times \cV$ are at most $\ell^{2}n^{\ell-2}$.

Next, consider distinct $e,e'\in E(G_1)$. For each $i\in[t]$, the number of $R\in \cO_i$ with $e\cup e' \In V(R)$ is at most $\ell^{3}n^{\ell-3}$. Summing over all $i\in[t]$ yields the desired bound.

Now, take $(i,v)\in \cV$ and $(i,c)\in \cC$. Each $R\in \cO_i$ with $\{(i,v),(i,c)\}\subseteq f_{i,R}$ will contain some $c$-edge~$e$. We distinguish two cases for $e$. If $e$ is incident to~$v$, there is only one choice for~$e$, and then at most $\ell^2 n^{\ell-2}$ choices left. If $e$ is not incident to~$v$, there are at most $n/2$ choices for $e$, and then at most $\ell^3 n^{\ell-3}$ choices left. Thus, in total, there are at most $2\ell^3n^{\ell-2}$ choices.

Similarly, we check that the codegree of $(i,c)$ and $(i,c')$ for distinct $c,c'\in C$ is at most $\ell^4 n^{\ell-2}$. We have to choose a $c$-edge $e$ and a $c'$-edge~$e'$ and again distinguish two cases. If $e$ and $e'$ share a vertex~$v$, there are at most $n$ choices for~$v$ (which determines $e$ and~$e'$), and then at most $\ell^3 n^{\ell-3}$ choices left. If $e$ and $e'$ form a matching, there are at most $(n/2)^2$ ways to choose $e$ and~$e'$, and then at most $\ell^4 n^{\ell-4}$ choices left.
\endclaimproof

 For each $v\in V$, let $\cV_v$ be the set of all pairs $(i,v)$ with $i\in [t]$ and $v\in V_{i,1}$. For each colour $c\in C$, let $\cC_c$ be the set of all pairs $(i,c)$ with $i\in [t]$ and $c\in C_{i,1}$. Let $$\cF:=\set{\Set{i}\times V_{i,1},\Set{i}\times C_{i,1}}{i\in[t]} \cup \set{\cV_v,\partial_{G_1}(v)}{v\in V} \cup \set{\cC_c,E_c(G_1)}{c\in C}.$$
Using \eqref{reg1}, \eqref{reg2}, \eqref{reg3}, and Claim~\ref*{claim:aux codegree approx}, we now apply Theorem~\ref{thm:AY} to obtain a $(\gamma,\cF)$-perfect matching $\cM$ in~$\cH$. For $i\in[t]$, let $\cO_i'$ be the collection of all $R\in \cO_i$ for which $f_{i,R}\in \cM$.
For distinct $R,R'\in\cO_i$, $\{i\}\times V(R)\subseteq f_{i,R}$ and $\{i\}\times V(R')\subseteq f_{i,R'}$, and therefore, as $\cM$ is a matching, $R$ and $R'$ are vertex-disjoint. Similarly, $R$ and $R'$ are colour-disjoint. Thus, $\cO_i'\subseteq \cO_i$ is a collection of vertex-disjoint $\ell$-cycles in $G_1[V_{i,1}]$ whose union is rainbow with colours in~$C_{i,1}$. Moreover, as for each $i\in [t]$ and $R\in \cO_i$, $E(R)\subseteq f_{i,R}$, all these cycles are edge-disjoint.

For each $i\in [t]$, we will now randomly break each cycle in $\cO_i'$ into a path, before joining all these paths together into a single cycle.
For each $i\in[t]$ and all $R\in \cO_i'$, choose an edge $e_{i,R}\in E(R)$ uniformly at random. For each $i\in [t]$, let $X_i:=\bigcup_{R\in \cO_i'}e_{i,R}\In V_{i,1}$ and $D_i:=\set{\phi(e_{i,R})}{R\in \cO_i'}\subseteq C_{i,1}$. We claim that, with high probability, we have
\begin{align}
	|\set{i\in[t]}{v\in X_i}|, |\set{i\in[t]}{c\in D_i}| \le 7n/\ell\;\;\;\; \mbox{ for all }v\in V,c\in C.\label{random edge deletion}
\end{align}
Indeed, fix a vertex $v\in V$. We have $\prob{v\in X_i}\le 2/\ell$ for all $i\in[t]$, and those events are independent. Similarly, for a fixed colour $c\in C$, we have $\prob{c\in D_i}\le 1/\ell$ for all $i\in[t]$, and those events are independent too.
Thus, the claim follows with Lemma~\ref{lem:chernoff}\ref{chernoff crude} and a union bound.

Now, assume that~\eqref{random edge deletion} holds.
For each $i\in[t]$, let $H_i$ be the graph obtained as follows: Give every edge $\Set{e_{i,R}}_{R\in \cO_i'}$ an orientation, and (cyclically) enumerate these edges. Now, for each edge, add a path of length~$2$ between its head and the tail of the next edge, using a new vertex as the internal vertex each time. (So $H_i$ consists of the union of all these $|\cO_i'|$ paths of length~$2$, but does not contain the edges~$e_{i,R}$. In particular, $X_i$ is independent in $H_i$.) Observe that $|V(H_i)|,|E(H_i)|\le 3|\cO_i'| \le 3n/\ell$.

By~\ref{greedy connection condition} and~\eqref{random edge deletion}, we can apply Lemma~\ref{lem:greedy connections} (with $G_2,\Set{V_{i,2}}_{i\in[t]},\Set{C_{i,2}}_{i\in[t]}$ taking the place of $G$,$\Set{V_{i}}_{i\in[t]}$,$\Set{C_{i}}_{i\in[t]}$) to find
for each $i\in[t]$ an embedding $\psi_i\colon H_i\to G_2$ such that $\psi_i(H_i)$ is rainbow with colours in~$C_{i,2}$ and $\psi_i(x)=x$ for all $x\in X_i$ and $\psi_i(x)\in V_{i,2}$ for all $x\in V(H_i)\sm X_i$, and such that $\psi_1(H_1),\dots,\psi_t(H_t)$ are edge-disjoint.

For $i\in[t]$, let $$\tilde{R}_i:=\bigcup_{R\in \cO_i'}(R-e_{i,R}) \cup \psi_i(H_i).$$

We have that $\tilde{R}_1,\dots,\tilde{R}_t$ are edge-disjoint rainbow cycles in~$G$, where $\tilde{R}_i$ is rainbow with colours in~$C_i$ and $V(\tilde{R}_i)\In V_i$.
Moreover, by the definition of $\cF$, \ref{approx reserve bounded} and~\eqref{random edge deletion} and the fact that $|D_i|= |\cO_i'|\le n/\ell$ for all $i\in[t]$, we have that $(G-\bigcup_{i\in[t]}\tilde{R}_i$, $\Set{V_i\sm V(\tilde{R}_i)}_{i\in[t]}$, $\Set{C_i \sm \phi(E(\tilde{R}_i))}_{i\in[t]})$ is $4\gamma n$-bounded.

Finally, choose for each $i\in[t]$ an edge $e_i\in E(\tilde{R}_i)$ uniformly at random and let $P_i:=\tilde{R}_i-e_i$.
For a vertex $v\in V$, let $I_v$ be the set of indices $i\in[t]$ for which $v\in V(\tilde{R}_i)$ and the subpath from $v$ to one of the endvertices of $P_i$ has length at most $\kappa n$. Note that, for each $i\in[t]$, as $|V_i\setminus V(\tilde{R}_i)|\le 4\gamma n$, the cycle $\tilde{R}_i$ has length at least $pn/2$ by~\ref{approx vertex size}, implying $\prob{i\in I_v} \le \frac{3\kappa n}{pn/2}= 6\kappa p^{-1}$, and these events are independent. Thus, with Lemma~\ref{lem:chernoff}\ref{chernoff crude}, we conclude that~\ref{path ends} holds with high probability.\COMMENT{The random variable $|I_v|$ is the sum of independent Bernoulli variables and has expectation at most $3\kappa p^{-1}n$, so can apply Lemma~\ref{lem:chernoff}\ref{chernoff crude} with $t=21\kappa p^{-1}n$}
Similarly, we can deduce that with high probability, for every $v\in V$, the number of $i\in[t]$ for which $v$ is incident with $e_i$, is at most~$\log^2 n$, and for every $c\in C$, the number of $i\in[t]$ for which $\phi(e_i)=c$, is at most~$\log^2 n$.
This implies that \ref{approx dec bounded leftover} still holds with high probability.
\endproof

\subsection{Matchings for colour absorption}\label{sec:colour}

In this subsection, we find the rainbow matchings which form the
crucial part of the colour absorption structure (see Lemma~\ref{cor:rainbow matchings}). The following lemma prepares the ground for this.

\begin{lemma} \label{lem:rainbow matchings}
Suppose $1/n \ll \gamma\ll \eta\ll 1$ and let $p:=2\eta$, $q:=\eta/192$ and $t:=n/2$. Let $\phi$ be a $1$-factorization of~$K_n$ with vertex set~$V$ and colour set~$C$. For every $i\in[t]$, let $V_i$ be a $p(1+\gamma)$-random subset of~$V$, and let $C_i$ be a $q$-random subset of~$C$. Moreover, let $G$ be an $\eta(1+\gamma)$-random subgraph of~$K_n$. Then, with high probability, there exist edge-disjoint matchings $M_1,\dots,M_t$ in $G$ such that the following hold:
\begin{enumerate}[label=\rm{(\roman*)}]
\item $V(M_i) \In V_i$ for all $i\in[t]$;
\item for all $i\in[t]$, $M_i$ consists of $192$ $c$-edges for each $c\in C_i$;\label{rainbow matchings:colours}
\item for every vertex $v\in V$, the number of $i\in[t]$ for which $v$ is covered by $M_i$ is $(1\pm 3\gamma)pt$.
\end{enumerate}
\end{lemma}

Later on some edges of $M_i$ will be used to construct the $i$th tree $T_i$ of the decomposition.

\proof
Choose a new constant $\eps>0$ such that $1/n \ll \eps \ll \gamma \ll \eta\ll 1$.
For each $i\in[t]$, we randomly split $V_i$ into a $p$-random set $V_{i,1}$ and a $p\gamma$-random set $V_{i,2}$. Similarly, we split $G$ into an $\eta$-random subgraph $G_1$ and an $\eta\gamma$-random subgraph~$G_2$. For $c\in C$ and $i\in[t]$, let $Y_{c,i}$ denote the number of $c$-edges in~$G_2[V_{i,2}]$.

We define a (random) auxiliary hypergraph $\cH$ as follows. The vertex set of $\cH$ consists of three parts: The first part is simply~$E(G_1)$. The second part is the set $\cV$ of all pairs $(i,v)$ with $i\in [t]$ and $v\in V_{i,1}$. The third part of $V(\cH)$ is the set $\cC$ of all triples $(i,c,\ell)$ with $i\in [t]$, $c\in C_i$ and $\ell \in [192]$.
For all $e=uv\in E(K_n)$, $i\in[t]$ and $\ell\in [192]$, we add the hyperedge
\begin{align}
	f_{e,i,\ell}:=\Set{e,(i,u),(i,v),(i,\phi(e),\ell)}
\end{align}
if and only if $e\in E(G_1)$, $u,v\in V_{i,1}$ and $\phi(e)\in C_i$. Thus, $\cH$ is a $4$-uniform hypergraph.

\begin{NoHyper}
\begin{claim}\label{claim:aux reg 1}
With high probability, for each $e\in E(G_1)$, $d_{\cH}(e)=(1\pm \eps)192tp^2q$.
\end{claim}
\end{NoHyper}

\claimproof
Fix an edge $e=uv$ and assume $e\in E(G_1)$.\COMMENT{We condition here on the event that $e\in E(G_1)$, which does not affect the randomness/independence of the following events.} For $i\in[t]$, let $X_i$ be the indicator variable of the event that $u,v\in V_{i,1}$ and $\phi(e)\in C_i$. Note that $d_{\cH}(e)=192\sum_{i\in[t]}X_i$. Since $\prob{X_i=1}=p^2q$ for each $i$ and the $X_i$'s are independent, we can easily deduce from Chernoff's bound that the claim holds.
\endclaimproof

\begin{NoHyper}
\begin{claim}\label{claim:aux reg 2}
With high probability, for each $(i,v)\in \cV$, $d_{\cH}((i,v))=(1\pm \eps)192npq\eta$.
\end{claim}
\end{NoHyper}

\claimproof
Fix $i\in [t]$ and $v\in V$ and assume $(i,v)\in \cV$. For every vertex $u\neq v$, let $X_u$ be the indicator variable of the event that $u\in V_{i,1}$, $uv\in E(G_1)$ and $\phi(uv)\in C_i$. Note that $d_{\cH}((i,v))=192 \sum_{u\in V\sm\Set{v}} X_u$. Since $\prob{X_u=1} = p \eta q$ for each $u$ and the $X_u$'s are independent,\COMMENT{different $u$'s produce different edges and also different colours} we can easily deduce from Chernoff's bound that the claim holds.
\endclaimproof

\begin{NoHyper}
\begin{claim}\label{claim:aux reg 3}
With high probability, for each $(i,c,\ell)\in \cC$, $d_{\cH}((i,c,\ell))=(1\pm \eps) p^2 \eta n/2$, and, for each $c\in C$ and $i\in[t]$, $Y_{c,i}= (1\pm \eps)(p\gamma)^2 (\eta \gamma) n/2$.
\end{claim}
\end{NoHyper}

\claimproof
Fix $i\in [t]$, $c\in C$ and $\ell\in [192]$ and assume that $(i,c,\ell)\in \cC$. For every $c$-edge $e\in E(K_n)$, let $X_e$ be the indicator variable of the event that $e\in E(G_1)$ and $e\In V_{i,1}$. Note that $d_{\cH}((i,c,\ell))=\sum_{e\in E_c(K_n)} X_e$. Since $\prob{X_e=1} = p^2 \eta$ for each $e\in E(K_n)$ and the $X_e$'s are independent,\COMMENT{since the $c$-edges form a matching, every vertex choice only affects one edge} we can easily deduce from Chernoff's bound that the claim holds. (A similar argument works for~$Y_{c,i}$.)
\endclaimproof

\begin{NoHyper}
\begin{claim}\label{claim:basics}
With high probability, we have $|\set{i\in[t]}{v\in V_{i,1}}|=(1\pm \eps)pt$ and $|\set{i\in[t]}{v\in V_{i,2}}|=(1\pm \eps)p\gamma t$ for all $v\in V$.
\end{claim}
\end{NoHyper}

\claimproof
This is an easy consequence of Chernoff's bound.
\endclaimproof

\begin{NoHyper}
\begin{claim}\label{claim:aux codeg}
$\Delta^c(\cH)\le 192$.
\end{claim}
\end{NoHyper}

\claimproof
Clearly, the codegree of pairs in $E(G_1)\times E(G_1)$ and $\cC\times \cC$ is~$0$. Moreover, the codegree of any pair in $E(G_1)\times \cC$ and $\cV\times \cC$ is at most $1$, and the codegree of any pair in $E(G_1)\times \cV$ is at most~$192$. Finally, consider a pair in $\cV\times \cV$, say $(i,u)$ and $(i',v)$. If $i\neq i'$, then the codegree is $0$, so assume $i=i'$. If $uv\notin E(G_1)$, then the codegree is also zero, so assume otherwise and let $c$ be the colour of $uv$. Then the codegree is at most~$192$.
\endclaimproof

We now assume that the properties stated in Claims~\ref*{claim:aux reg 1}--\ref*{claim:aux codeg} are satisfied.
By our choice of $p,q,\eta$, we have that $d_{\cH}(x)=(1\pm \eps)192 pq\eta n$ for all $x\in V(\cH)$.
For every vertex $v\in V$, let $\cV_v$ be the set of all pairs $(i,v)\in \cV$ with $i\in[t]$ and $v\in V_{i,1}$. For every $c\in C$, let $\cC_c$ be the set of all $(i,c,\ell)$ with $c\in  C_i$ and $\ell\in[192]$.
Let $$\cF:=\set{\cV_v}{v\in V} \cup \set{\Set{i}\times C_i \times [192]}{i\in[t]} \cup \set{\cC_c}{c\in C}$$
Thus, we can apply Theorem~\ref{thm:AY} to find a $(\gamma^5,\cF)$-perfect matching~$\cM$ in~$\cH$.
For $i\in[t]$, let $M_i'$ be the set of all edges $e\in E(G_1)$ such that $f_{e,i,\ell}\in \cM$ for some $\ell\in[192]$.

Hence, by definition of $\cH$, we have that $M_1',\dots,M_t'$ are edge-disjoint matchings in~$G_1$, and, for each $i\in[t]$, we have $V(M_i')\In V_i$ and $M_i'$ consists of at most 192 edges with colour $c$, for each $c\in C_i$.

For each $i\in[t]$ and $c\in C_i$, let $r_{i,c}:=192-|M_i'\cap E_c(K_n)|$. Thus, $r_{i,c}$ is the number of $c$-edges that are missing in $M_i'$ in order to satisfy~\ref{rainbow matchings:colours}. Since $\cM$ is $(\gamma^5,\cF)$-perfect, we have for each $i\in [t]$ and $c\in C$ that
\begin{align}
	\sum_{c'\in C_i}r_{i,c'} \le \gamma^5|\{i\}\times C_i\times [192]|\le 192\gamma^5 n \;\;\mbox{ and } \;\;\sum_{i'\in [t]:c\in C_{i'}}r_{i',c} \le \gamma^5|\cC_c|  \le 192\gamma^5 n. \label{small leftover}
\end{align}
Moreover, for each vertex $v\in V$, the number of $i\in[t]$ for which $v\in V_{i,1}$ but $v$ is not covered by $M_i'$, is at most $\gamma^5 n$.
Since $|\set{i\in[t]}{v\in V_{i,1}}|=(1\pm \eps)pt$ by Claim~\ref*{claim:basics}, this implies that the number of $i\in[t]$ for which $v$ is covered by $M_i'$ is $(1\pm \gamma^4)pt$.

Now, we wish to find edge-disjoint matchings $M_1'',\dots, M_t''$ in $G_2$ such that, for each $i\in[t]$, $V(M_i'')\In V_{i,2}$ and $M_i''$ contains precisely $r_{i,c}$ $c$-edges for each $c\in C$.
This can be done in order greedily using Claim~\ref*{claim:aux reg 3} and~\eqref{small leftover}. Indeed, suppose we want to add $c$-edges to $M_i''$. By~\eqref{small leftover}, we added at most $192\gamma^5 n$ $c$-edges to previous matchings $M_j''$, $j<i$, and at most $192\gamma^5 n$ edges to~$M_i''$. Thus, at most $3\cdot 192\gamma^5 n$ $c$-edges are blocked (since every edge in $M_i''$ might block $2$ $c$-edges). Since $Y_{c,i}\ge \gamma^4 n$ by Claim~\ref*{claim:aux reg 3}, we can find $r_{i,c}$ suitable $c$-edges in~$G_2[V_{i,2}]$ and add them to $M_i''$.

Note that, by Claim~\ref*{claim:basics}, for every vertex $v\in V$, the number of $i\in[t]$ for which $v$ is covered by~$M_i''$, is at most $2\gamma pt$.
Finally, for each $i\in[t]$, let $M_i:=M_i'\cup M_i''$. It is easy to see that $M_1,\dots,M_t$ are the desired matchings.
\endproof

\begin{lemma} \label{cor:rainbow matchings}
Suppose $1/n \ll \gamma \ll \eta \ll 1$ and let $p:=2\eta$, $q:=\eta/192$ and $t:=n/2$. Suppose $s\in \bN$ with $s=(q/4- 2\gamma/5\pm \gamma^2) n$ and $0\le \alpha \le q/2-\gamma$. Let $\phi$ be a $1$-factorization of~$K_n$ with vertex set~$V$ and colour set~$C$. Let $G$ be an $\eta(1+\gamma)$-random subgraph of~$K_n$. For every $i\in[t]$, let $V_i$ be a $p(1+\gamma)$-random subset of~$V$, and let $C_{i,1}$, $C_{i,2}$ be disjoint $q/2$-random subsets of~$C$. Split $C_{i,1}$ further into an $\alpha$-random set $C_{i,1,1}$ and a $(q/2-\alpha)$-random subset $C_{i,1,2}$.

Then with high probability, for each $i\in[t]$, there exist $C_{i,1}',C_{i,2}'$ such that $C_{i,1,1} \In C_{i,1}'\In C_{i,1}$ and $ C_{i,2}'\In C_{i,2}$ and
vertex-disjoint rainbow matchings $\set{M_{i,j}}{j\in[3s]}$ in~$G[V_i]$, such that altogether the following hold:
\begin{enumerate}[label=\rm{(\roman*)}]
\item for each $i\in[t]$, $|C_{i,1}'|=|C_{i,2}'|=2s$;\label{rainbow matchings sizes}
\item for all $c\in C$, $|\set{i\in[t]}{c\in C_{i,2}\sm C_{i,2}'}|\le \sqrt{\gamma} n$;\label{rainbow matchings bounded}
\item for each $i\in[t]$, $M_i:=\bigcup_{j\in[3s]}M_{i,j}$ consists of $192$ $c$-edges for each $c\in C_{i,1}'\cup C_{i,2}'$; \label{rainbow matchings partition}
\item for each $i\in [t]$ and any subset $C_i^\ast\In C_{i,1}'$ of size~$s$, there exists $J_i\In M_i$ such that $J_i$ is $(C_i^\ast\cup C_{i,2}')$-rainbow and contains exactly one edge from each of $\set{M_{i,j}}{j\in[3s]}$; \label{colour absorber}
\item the matchings $M_1,\dots,M_t$ are edge-disjoint, and $|M_{i,j}|=256$ for all $(i,j)\in [t]\times [3s]$; \label{rainbow matchings partition 2}
\item for every vertex $v\in V$, the number of $i\in[t]$ for which $v$ is covered by $M_i$ is $(1\pm \sqrt{\gamma})pt$. \label{degree control vertices rainbow}
\end{enumerate}
\end{lemma}

Here, the crucial property is~\ref{colour absorber}, which will allow us to use some colours of $C_{i,1}'$ flexibly before assigning the remaining colours (i.e.~those in~$C_i^\ast$) together with the `buffer' $C_{i,2}'$ in such a way that each matching $\set{M_{i,j}}{i\in[t],j\in [3s]}$ contributes exactly one edge to $J_i$ which will be part of~$T_i$.

\proof
Clearly, we may assume that $\alpha=q/2-\gamma$.
We choose the random colour sets according to the following procedure: For each $i\in[t]$, let $C_i$ be a $q$-random subset of~$C$, and let $\tau_i\colon C\to [4]$ be a random function such that $\prob{\tau_i(c)=1}=\prob{\tau_i(c)=2}=1/2-\gamma/q$ and $\prob{\tau_i(c)=3}=\prob{\tau_i(c)=4}=\gamma/q$.
For $i\in[t]$ and $k\in[2]$, let $C_{i,k,1}=\set{c\in C_i}{\tau_i(c)=k}$ and $C_{i,k,2}=\set{c\in C_i}{\tau_i(c)=k+2}$, and let $C_{i,k}:=C_{i,k,1} \cup C_{i,k,2}$.
Then $C_{i,1},C_{i,2},C_{i,1,1},C_{i,1,2}$ are as in the statement.

Now, we first expose all random choices except the functions $\Set{\tau_i}_{i\in[t]}$.
By Lemma~\ref{lem:rainbow matchings}, with high probability, there exist edge-disjoint matchings $M_1',\dots,M_t'$ in~$G$ such that the following hold:
\begin{enumerate}[label=\rm{(\alph*)}]
\item $V(M_i') \In V_i$ for all $i\in[t]$;
\item for all $i\in[t]$, $M_i'$ consists of $192$ $c$-edges for each $c\in C_i$;
\item for every vertex $v\in V$, the number of $i\in[t]$ for which $v$ is covered by $M_i'$ is $(1\pm 3\gamma)pt$. \label{auxiliary matching degree}
\end{enumerate}

Henceforth, assume that these random choices have been made and satisfy the above properties. It remains to expose the functions~$\tau_i$.

With high probability, we have for all $i\in[t]$ and $k\in[2]$ that
\begin{align}
	|C_{i,k,1}|= (1\pm \gamma^2)(q/2-\gamma) n \quad \mbox{ and } \quad |C_{i,k,2}|= (1\pm \gamma^2)\gamma n. \label{colour set sizes}
\end{align}
With high probability, we also have for all $c\in C$ that
\begin{align}
	|\set{i\in [t]}{\tau_i(c)=4}|\le \sqrt{\gamma}n. \label{colour set incidences}
\end{align}

For $v\in V$, let us call $i\in[t]$ \defn{unreliable for $v$} if $v$ is covered by $M_i'$ via an edge whose colour is in $C_{i,1,2}\cup C_{i,2,2}$.
Then, also with high probability, for all $v\in V$,
\begin{align}
	\mbox{at most $2\gamma q^{-1}n$ indices $i\in[t]$ are unreliable for $v$.} \label{matching deletion control}
\end{align}\COMMENT{For each $i\in[t]$, the probability that it is unreliable is at most~$2\gamma/q$.}

From now on, assume that \eqref{colour set sizes}--\eqref{matching deletion control} hold.
For each $i\in[t]$ and $k\in[2]$, note that by \eqref{colour set sizes} we have $2s-2\gamma n/5 \le |C_{i,k,1}| \le 2s$ and thus, again by \eqref{colour set sizes}, we can choose $C_{i,k,2}'\In C_{i,k,2}$ such that $|C_{i,k,2}'|=2s-|C_{i,k,1}|$, and define $$C_{i,k}':= C_{i,k,1} \cup C_{i,k,2}'.$$ Then \ref{rainbow matchings sizes} clearly holds and \ref{rainbow matchings bounded} follows from~\eqref{colour set incidences}.
Moreover, let $$M_i:= \set{e\in M_i'}{\phi(e)\in C_{i,1}'\cup C_{i,2}'}$$ for each $i\in[t]$.
Observe that \ref{auxiliary matching degree} and~\eqref{matching deletion control} imply \ref{degree control vertices rainbow}.

We now use RMBG's to break each $M_i$ into small rainbow matchings. For each $i\in [t]$, let $H_i$ be a $(256,192)$-regular $\mathrm{RMBG}(3s,2s,2s)$ with parts~$[3s],C_{i,1}'$ and $C_{i,2}'$, which exist by Corollary~\ref{cor:regular rmbg}.
For each $i\in[t]$, partition $M_i$ into matchings $M_{i,1},\dots,M_{i,3s}$, such that, for each $j\in[3s]$, $M_{i,j}$ is an $N_{H_i}(j)$-rainbow matching. This can be done greedily since to do so we need precisely $192$ $c$-edges of each colour $c\in C_{i,1}'\cup C_{i,2}'$, which $M_i$ contains.
Clearly, $|M_{i,j}|=|N_{H_i}(j)|=256$, and thus \ref{rainbow matchings partition} and \ref{rainbow matchings partition 2} hold.

Finally, we check that the crucial property~\ref{colour absorber} holds.
Consider $i\in [t]$ and suppose $C_i^\ast\In C_{i,1}'$ has size~$s$. Since $H_i$ is an RMBG with parts $[3s]$, $C_{i,1}'$ and $C_{i,2}'$, there exists a perfect matching $\tau$ in $H_i$ between $[3s]$ and $C_i^\ast \cup C_{i,2}'$. Now, for each $j\in[3s]$, we select the $\tau(j)$-edge from $M_{i,j}$ and include it in $J_i$. (Here we view $\tau(j)$ as the colour matched to $j$ in the matching $\tau$, and we use that $M_{i,j}$ is $N_{H_i}(j)$-rainbow.) Clearly, $J_i$ is as desired.
\endproof

\subsection{Matchings for edge absorption}\label{sec:edge}
We now find the monochromatic matchings which form the crucial ingredients
for the edge absorption process.

\begin{lemma} \label{lem:monochromatic matchings}
Suppose $1/n \ll \eps \ll \gamma \ll \rhosub\ll 1$ and let $p:=3072\rhosub$, $q:=6\rhosub$, $t:=n/2$, and suppose $m\in \bN$ with $m=(\rhosub-\eps/5 \pm \eps^2) n$ and $0\le \alpha \le 4\rhosub-\eps$.
Let $\phi$ be a $1$-factorization of~$K_n$ with vertex set~$V$ and colour set~$C$.
Let $G_1',G_2'$ be edge-disjoint $4\rhosub$-random subgraphs of~$K_n$, and split $G_1'$ further into an $\alpha$-random subgraph $G_{1,1}'$ and a $(4\rhosub-\alpha)$-random subgraph~$G_{1,2}'$.
For each $i\in[t]$, let $V_i$ be a $p(1+\gamma)$-random subset of $V$ and let $D_i$ be a $q(1+\gamma)$-random subset of~$C$.

Then, with high probability, there exist $G_1$ and $G_2$ such that $G_{1,1}' \In G_1 \In G_1'$ and $G_2\In G_2'$ with $\Delta(G_2'-G_2)\le 2\eps n$ and,
for each $i\in[t]$, there exists $D_i'\In D_i$ of size $(1\pm 2\gamma)q n$ and vertex-disjoint matchings $\Set{M_{i,c}}_{c\in D_i'}$ in $(G_1\cup G_2)[V_i]$, where $M_{i,c}$ consists of $256$~$c$-edges, such that altogether the following hold:
\begin{enumerate}[label=\rm{(\roman*)}]
\item for each $c\in C$, $|E_c(G_1)|=|E_c(G_2)|=2m $ and $|\set{i\in[t]}{c\in D_i'}|=3m$; \label{colour regularity}
\item for any subset $E^\ast\In E(G_1)$ which consists of precisely $m$ edges of each colour $c\in C$, there exists a partition of $E^\ast\cup E(G_2)$ into sets $J_1,\dots,J_t$, such that, for each $i\in[t]$, $J_i$ contains exactly one edge from each of $\Set{M_{i,c}}_{c\in D_i'}$; \label{edge absorber}
\item every vertex $v\in V$ is covered by $(1\pm \gamma)pt$ of the matchings $\set{M_{i,c}}{i\in[t],c\in D_i'}$. \label{leftover control degree}
\end{enumerate}
\end{lemma}

\COMMENT{The second condition in \ref{colour regularity} also follows from the first and \ref{edge absorber}}

Here, the crucial property is~\ref{edge absorber}, which will allow us to use some edges of the global edge reservoir $G_1$ flexibly before assigning the remaining edges (i.e.~those in~$E^\ast$) together with the `buffer' $E(G_2)$ in such a way that each matching $\set{M_{i,c}}{i\in[t],c\in D_i'}$ contributes exactly one edge to $J_i$.
$J_i$ will then be assigned to the $i$th tree $T_i$.

\proof
We may clearly assume that $\alpha=4\rhosub-\eps$.
We also split $G_2'$ further into a $(4\rhosub -\eps)$-random subgraph $G_{2,1}'$ and an $\eps$-random subgraph~$G_{2,2}'$.
We first expose $G_{1,1}',G_{1,2}',G_{2,1}',G_{2,2}'$. Using Chernoff's bound, it is easy to see that, with high probability, we have for all $j\in[2]$, $c\in C$ and $v\in V$ that
\begin{align}
	|E_c(G_{j,1}')| = (1\pm \eps^2) (4\rhosub -\eps) n/2 \quad \mbox{ and }\quad  |E_c(G_{j,2}')| = (1\pm \eps^2) \eps n/2,  \label{edge splitting colour degrees} \\
	d_{G_{j,1}'}(v) = (1\pm \eps^2) (4\rhosub -\eps) n  \quad \mbox{ and }\quad  d_{G_{j,2}'}(v) = (1\pm \eps^2) \eps n. \label{edge splitting vertex degrees}
\end{align}
Henceforth, we assume that $G_{1,1}',G_{1,2}',G_{2,1}',G_{2,2}'$ are fixed with the above properties, and expose the other random sets.

By~\eqref{edge splitting colour degrees}, we have for $j\in[2]$ and $c\in C$ that $2m-\eps n/5 \le |E_c(G_{j,1}')| \le 2m$.\COMMENT{$2m - |E_c(G_{j,1}')| =  2(\rhosub-\eps/5 \pm \eps^2) n  - (1\pm \eps^2) (4\rhosub -\eps) n/2 = \eps n/10 \pm 3\eps^2 n$} Therefore, by \eqref{edge splitting colour degrees} again, for each $j\in[2]$, we can choose $G_{j,2}''\In G_{j,2}'$ such that $|E_c(G_{j,2}'')| = 2m - |E_c(G_{j,1}')|$ for all $c\in C$, and define $$G_j:=G _{j,1}' \cup G_{j,2}''.$$
Clearly, this choice of $G_1$ and $G_2$ satisfies the first part of~\ref{colour regularity}.
Moreover, from~\eqref{edge splitting vertex degrees}, we can infer that $\Delta(G_2'-G_2)\le 2\eps n$, as desired, and that
\begin{align}
	\mbox{$d_{G_1\cup G_2}(v)=(1 \pm \sqrt{\eps})8\rhosub n$ for all $v\in V$.} \label{final absorber graphs degree}
\end{align}

As indicated in the proof sketch, the key to obtaining \ref{edge absorber} is to use an RMBG for each colour which matches the $3m$~$c$-edges of $E^\ast \cup E(G_2)$ to $3m$~`absorbers'. Let $\hat{H}$ be a $(256,192)$-regular $\mathrm{RMBG}(3m,2m,2m)$ with parts~$[3m],\hat{Y_1},\hat{Y_2}$, which exists by Corollary~\ref{cor:regular rmbg}. We identify $\hat{Y_1}$ and $\hat{Y_2}$ with $E_c(G_1)$ and $E_c(G_2)$. We carry out this identification randomly in order to obtain a codegree condition in some hypergraph $\cH$ which we will define later. (This codegree condition will be needed when applying Theorem~\ref{thm:AY} to $\cH$.) For each colour $c\in C$, pick random bijections $\pi_{c,1}\colon E_c(G_1) \to \hat{Y_1}$ and $\pi_{c,2}\colon E_c(G_2) \to \hat{Y_2}$, all independently. Obtain a copy $H_c$ of $\hat{H}$ by identifying $E_c(G_1)$ with $\hat{Y_1}$ according to $\pi_{c,1}$ and $E_c(G_2)$ with $\hat{Y_2}$ according to~$\pi_{c,2}$.

For two vertices $v,v'$, we define $r_{v,v'}$ as the number of colours $c\in C$ for which $N_{H_c}(e)\cap N_{H_c}(e')\neq \emptyset$, where $e$ and $e'$ are the unique $c$-edges at $v$ and $v'$, respectively. (In particular, if $e$ or $e'$ is not contained in $E_c(G_1\cup G_2)$, then $c$ contributes $0$ to $r_{v,v'}$.)

\begin{NoHyper}
\begin{claim}\label{claim:random edge perm}
With positive probability, $r_{v,v'}\le 3\log{n}$ for all distinct vertices $v,v'\in V$.
\end{claim}
\end{NoHyper}

\claimproof
Fix two distinct vertices $v,v'\in V$. For $c\in C$, let $X_c$ be the indicator variable of the event that there exist $c$-edges $e,e'$ at $v$ and $v'$, respectively (which are unique if existent), and $N_{H_c}(e)\cap N_{H_c}(e')\neq \emptyset$. Fix $c\in C\sm\Set{\phi(vv')}$ and let $e,e'$ be as above. We claim that $\prob{X_c=1}\le 10^5 m^{-1}$.  Note that $e,e'$ are distinct. Let $k,k'\in\Set{1,2}$ be such that $e\in E(G_k)$ and $e'\in E(G_{k'})$.
Thus,
\begin{align*}
\prob{X_c=1} \le \sum_{j\in[3m]} \prob{\pi_{c,k}(e),\pi_{c,k'}(e')\in N_{\hat{H}}(j)} \le 3m \cdot \frac{256\cdot 255}{2m\cdot(2m-1)} \le 10^5 m^{-1}.
\end{align*}
Hence, $\expn{r_{v,v'}}\le \rhosub^{-2}$, and since the $X_c$'s are independent, Chernoff's bound implies that the probability that $r_{v,v'}>3\log{n}$ is smaller than $n^{-2}$. A union bound then implies the claim.
\endclaimproof

From now on, fix RMBG's $\Set{H_c}_{c\in C}$ for which the conclusion of Claim~\ref*{claim:random edge perm} holds. Let $\cA:=C\times [3m]$. For each $(c,j)\in \cA$, we define $A_{c,j}:=N_{H_c}(j)$. We refer to $A_{c,j}$ as an \defn{absorber} and will sometimes identify $A_{c,j}$ with $(c,j)\in \cA$. Note that $A_{c,j}$ is a matching consisting of $256$~$c$-edges.

By our choice of RMBG's, we have that for any two distinct vertices $v,v'\in V$,
\begin{align}
	\mbox{there are at most $192\cdot 3\log{n}$ absorbers $(c,j)\in \cA$ with $v,v'\in V(A_{c,j})$.} \label{vertex codegree}
\end{align}
\COMMENT{for each colour $c\in C$ with $N_{H_c}(e)\cap N_{H_c}(e')\neq \emptyset$, we have $|N_{H_c}(e)\cap N_{H_c}(e')|\le 192$}
 We will now assign to each absorber an index $i\in[t]$.
\COMMENT{In the final decomposition, exactly one of the $256$ edges of such an absorber will be used in the $i$-th tree.}
The assignment will be obtained as follows: We first define an auxiliary hypergraph~$\cH$, in which we will find an almost perfect matching that provides an almost complete assignment. For the remaining absorbers not yet assigned, we will greedily pick images from a reserve.

In order to set aside this `reserve', we randomly split $V_i$ and $D_i$ further as follows. For each $i\in[t]$, split $V_i$ into a $p$-random set $V_{i,1}$ and a $p\gamma$-random set $V_{i,2}$, and split $D_i$ into a $q$-random set $D_{i,1}$ and a $q\gamma$-random set $D_{i,2}$.

We can now define the (random) auxiliary hypergraph $\cH$ as follows. The vertex set of $\cH$ consists of three different parts: The first part is simply the set $\cA$ which represents all the absorbers. The second part is the set $\cV$ of all pairs $(i,v)$ with $i\in [t]$ and $v\in V_{i,1}$. The third part is the set $\cC$ of all pairs $(i,c)$ with $i\in [t]$ and $c\in D_{i,1}$.

Now, we define the edge set of~$\cH$. For every $i\in[t]$ and every absorber $(c,j)\in \cA$, we add the hyperedge
\begin{align}
f_{c,j,i}:=\Set{(c,j),(i,c)} \cup (\Set{i}\times V(A_{c,j}))
\end{align}
if and only if $c\in D_{i,1}$ and $V(A_{c,j}) \In V_{i,1}$.
Hence, $\cH$ is $514$-uniform. (Recall that $A_{c,j}$ is a matching consisting of $256$~$c$-edges.)

Moreover, for each absorber $(c,j)\in \cA$, we define the random set $Y_{c,j}$ of indices $i\in[t]$ for which $c\in D_{i,2}$ and $V(A_{c,j}) \In V_{i,2}$.
We aim to apply Theorem~\ref{thm:AY} to~$\cH$. For this, we first establish the following properties.

\begin{NoHyper}
\begin{claim}\label{claim:aux reg 4}
With high probability, for each $(c,j)\in \cA$, $d_{\cH}((c,j))=(1\pm \eps)tp^{512} q$ and $|Y_{c,j}|=(1\pm \eps)t(p\gamma)^{512} q\gamma$.
\end{claim}
\end{NoHyper}

\claimproof
Fix an absorber $(c,j)\in \cA$. For $i\in[t]$, let $X_i$ be the indicator variable of the event that $c\in D_{i,1}$ and $V(A_{c,j})\In V_{i,1}$ and let $Y_i$ be the indicator variable of the event that $c\in D_{i,2}$ and $V(A_{c,j})\In V_{i,2}$. Note that $d_{\cH}((c,j))=\sum_{i\in[t]}X_i$ and $|Y_{c,j}|=\sum_{i\in[t]}Y_i$. For each $i\in[t]$, we have that $\prob{X_i=1}=p^{512} q$ and $\prob{Y_i=1}=(p\gamma)^{512}q\gamma$. Thus, $\expn{d_{\cH}((c,j))}=tp^{512} q$ and $\expn{|Y_{c,j}|}=t(p\gamma)^{512} q\gamma$. Since the $X_i$'s are independent, and similarly, the $Y_i$'s are independent, we can deduce with Chernoff's bound that the claim holds.
\endclaimproof

\begin{NoHyper}
\begin{claim}\label{claim:aux reg 5}
With high probability, for each $(i,c)\in \cC$, $d_{\cH}((i,c))=(1\pm \eps)3m p^{512}$.
\end{claim}
\end{NoHyper}

\claimproof
Fix $(i,c)\in \cC$. For $j\in[3m]$, let $X_j$ be the indicator variable of the event that $V(A_{c,j})\In V_{i,1}$. Note that $d_{\cH}((i,c))=\sum_{j\in[3m]}X_j$. For each $j\in[3m]$ we have $\prob{X_j=1} = p^{512}$. Thus, $\expn{d_{\cH}((i,c))}=3m p^{512}$.

Moreover, $d_{\cH}((i,c))$ is determined by the independent random variables $\set{\ind{v\in V_{i,1}}}{v\in V}$. Since $\ind{v\in V_{i,1}}$ affects $d_{\cH}((i,c))$ by at most $192$,\COMMENT{$c$ is fixed, there is only one $c$-edge at $v$, and this edge is contained in $192$ absorbers~$(c,j)$.} the claim follows by an application of McDiarmid's inequality.
\endclaimproof

\begin{NoHyper}
\begin{claim}\label{claim:aux reg 6}
With high probability, for all $(i,v)\in \cV$, $d_{\cH}((i,v))=(1\pm 2\sqrt{\eps})1536\rhosub p^{511}q n$.
\end{claim}
\end{NoHyper}

\claimproof
Fix $(i,v)\in \cV$. For each edge $e$ at $v$ in $G_1\cup G_2$, say with colour $c$, $e$ has $192$ neighbours $j\in[3m]$ in $H_c$, and for each of those we have $f_{c,j,i}\in E(\cH)$ iff $c\in D_{i,1}$ and the $511$ other vertices of $A_{c,j}$ are contained in~$V_{i,1}$. Thus, $\expn{d_{\cH}((i,v))}= d_{G_1\cup G_2}(v) \cdot 192 \cdot p^{511} q = (1\pm \sqrt{\eps})\cdot 8\cdot 192\rhosub p^{511}q n$ by~\eqref{final absorber graphs degree}.

Moreover, $d_{\cH}((i,v))$ is determined by the independent random variables $\set{\ind{u\in V_{i,1}}}{u\in V\sm\Set{v}}\cup \set{\ind{c\in D_{i,1}}}{c\in C}$. The effect of $\ind{c\in D_{i,1}}$ on $d_{\cH}((i,v))$ is at most~$192$.\COMMENT{there is only one $c$-edge at~$v$, and this edge is contained in $192$ absorbers~$(c,j)$.} Moreover, for each $u\in V\sm\Set{v}$, by~\eqref{vertex codegree}, $\ind{u\in V_{i,1}}$ affects $d_{\cH}((i,v))$ by at most~$192\cdot 3\log{n}$.
The claim now follows from an application of McDiarmid's inequality.
\endclaimproof

\begin{NoHyper}
\begin{claim}\label{claim:aux codeg 2}
$\Delta^c(\cH)\le 192\cdot 3\log{n}$.
\end{claim}
\end{NoHyper}

\claimproof
Clearly, the codegree of pairs in $\cA\times \cA$ and $\cC\times \cC$ is~$0$. Moreover, the codegree of pairs in $\cA\times \cV$ and $\cA\times \cC$ is at most~$1$.\COMMENT{since $\cA$ fixes $c,j$ and $\cV$/$\cC$ fixes $i$}
It is also easy to see that the codegree of a pair in $\cV\times \cC$ is at most~$192$.\COMMENT{Say $(i,v)$ and $(i',c)$. If $i\neq i'$, then the codegree is~$0$, so assume $i=i'$. Let $e$ be the $c$-edge incident to~$v$ (if not existent, then codegree is~$0$.) Then $e$ is contained in $192$~absorbers.}

Finally, consider a pair in $\cV\times \cV$, say $(i,u)$ and $(i',v)$. If $i\neq i'$, then the codegree is~$0$, so assume $i=i'$. Crucially, by~\eqref{vertex codegree}, the codegree of $(i,u)$ and $(i,v)$ is at most $192\cdot 3\log{n}$.
\endclaimproof

We now assume that the properties stated in Claims~\ref*{claim:aux reg 4}--\ref*{claim:aux codeg 2} are satisfied. Using Chernoff's bound, we can assume that the following simple properties hold as well:
\begin{align}
	|D_{i,1}| = (1\pm \eps)q n \;\;\;\mbox{ and }\;\;\;|D_{i,2}| = (1\pm \eps)q\gamma n. \label{colour splitting sizes}
\end{align}

By our choice of $p,q,\rhosub,t,m$, we have that $d_{\cH}(x)=(1\pm 2\sqrt{\eps}) 3\rhosub p^{512} n$ for all $x\in V(\cH)$. In combination with Claim~\ref*{claim:aux codeg 2}, we can thus apply Theorem~\ref{thm:AY} to find an almost perfect matching in~$\cH$. In order to gain control over the leftover vertices in $\cH$, we define the following vertex sets. For each vertex $v\in V$, let $\cA_v$ be the set of all absorbers $(c,j)\in \cA$ for which $v\in V(A_{c,j})$. Note that
\begin{align}
	\mbox{$|\cA_v|=192d_{G_1\cup G_2}(v) \overset{\eqref{final absorber graphs degree}}{=}(1\pm \sqrt{\eps})\cdot8 \cdot 192\rhosub n =(1\pm \sqrt{\eps})pt$.} \label{vertex control set}
\end{align}
Define
$$\cF:=\set{\cA_v}{v\in V} \cup \set{\Set{c}\times [3m]}{c\in C} \cup \set{\Set{i}\times D_{i,1}}{i\in[t]}.$$

Now, apply Theorem~\ref{thm:AY} to find a
\begin{align}
\mbox{$(\gamma^{515},\cF)$-perfect matching~$\cM$ in~$\cH$.}\label{find matching}
\end{align}
Our goal is to define a map $\sigma\colon \cA\to [t]$.
Let $\cA'$ be the set of absorbers $(c,j)\in \cA$ which are not covered by~$\cM$.
For each $(c,j)\in \cA\sm \cA'$, the absorber $(c,j)$ is covered by a (unique) hyperedge $f_{c,j,i}\in\cM$, and we define $\sigma(c,j):=i$. For all uncovered absorbers, we now use the `reserve' sets $V_{i,2}$ and $D_{i,2}$ to pick suitable images.

For all $(c,j)\in \cA'$, we successively define $\sigma(c,j)$ as follows: when we consider $(c,j)\in \cA'$, let $\cA''$ be the set of all previously considered $(c',j')\in \cA'$ with $c'=c$ or $V(A_{c,j})\cap V(A_{c',j'})\neq \emptyset$. By~\eqref{find matching}, we have that $$|\cA''|\le \gamma^{515} \cdot |\Set{c}\times [3m]| + \sum_{v\in V(A_{c,j})} \gamma^{515}|\cA_v| \overset{\eqref{vertex control set}}{\le} \gamma^{515} \cdot 3m + 512 \gamma^{515}\cdot 2 pt <\gamma^{514}n/2 .$$ Recall from Claim~\ref*{claim:aux reg 4} that $|Y_{c,j}|\ge \gamma^{514}n$. Thus, there is $i\in Y_{c,j}\sm \sigma(\cA'')$ and we define $\sigma(c,j):=i$.

Altogether, we have found a map $\sigma\colon \cA\to [t]$, which we show has the following properties:
\begin{enumerate}[label=\rm{(\alph*)}]
\item $V(A_{c,j})\In V_{\sigma(c,j)}$ and $c\in D_{\sigma(c,j)}$ for all $(c,j)\in \cA$; \label{abs assignment slices}
\item $V(A_{c,j})\cap V(A_{c',j'})=\emptyset$ whenever $\sigma(c,j)=\sigma(c',j')$; \label{abs assignment disjoint}
\item for all $c\in C$ and $i\in[t]$, there is at most one $j\in[3m]$ with $\sigma(c,j)=i$. \label{abs assignment injective}
\end{enumerate}

Here, \ref{abs assignment slices} clearly holds by the definitions of $\cH$, $Y_{c,j}$ and~$\sigma$. To see \ref{abs assignment disjoint}, suppose $\sigma(c,j)=\sigma(c',j')=i$. If $(c,j),(c',j')\in \cA\sm \cA'$, then we have $V(A_{c,j})\cap V(A_{c',j'})=\emptyset$ since $\cM$ is a matching and as such covers every vertex $(i,v)\in \cV$ at most once. If $(c,j)\in \cA$ and $(c',j')\in \cA'$, then $V(A_{c,j})\In V_{i,1}$ and $V(A_{c',j'})\In V_{i,2}$. Finally, suppose $(c,j),(c',j')\in \cA'$ and assume that we defined $\sigma(c,j)$ after $\sigma(c',j')$. If $V(A_{c,j})\cap V(A_{c',j'})\neq \emptyset$, then $(c',j')\in \cA''$ (with notation as above) and hence $i\in\sigma(\cA'')$, a contradiction.

For \ref{abs assignment injective}, fix $c\in C$ and $i\in[t]$. Suppose $\sigma(c,j)=i$ for some $j\in[3m]$. We consider two cases. In the first case, we have $(c,j)\in \cA\sm \cA'$ and $f_{c,j,i}\in \cM$. In particular, there is at most one $j$ which satisfies this and we must have $c\in D_{i,1}$. In the second case, we must have $(c,j)\in \cA'$ and $c\in D_{i,2}$, and there can only be one $j$ which satisfies this by definition of $\cA''$ above. Since $D_{i,1}$ and $D_{i,2}$ are disjoint, \ref{abs assignment injective} follows.

\bigskip
Now, for every $c\in C$, define $\sigma_c:= \sigma(c,\cdot)$ and $I_c:=\Ima(\sigma_c)$. By~\ref{abs assignment injective}, $\sigma_c\colon [3m] \to I_c$ is a bijection.
For all $i\in[t]$, define $$D_i':=\set{c\in C}{i\in I_c}.$$
For all $c\in C$, we have $|\set{i\in[t]}{c\in D_i'}|=|I_c|=3m$, so the second part of~\ref{colour regularity} holds too.
Observe that if $c\in D_i'$, then $i\in I_c$ and hence there exists some $j\in[3m]$ for which $\sigma(c,j)=i$. By~\ref{abs assignment slices}, we have $c\in D_i$. Thus, $D_i'\In D_i$. In particular, we have $|D_i'|\le |D_i|\le (1+2\gamma)q n$ by~\eqref{colour splitting sizes}. Moreover, since $\Set{i}\times D_{i,1}\in \cF$, at least $(1-\gamma^{515})|D_{i,1}|$ elements of $\Set{i}\times D_{i,1}$ are covered by~$\cM$, which means that for at least $(1-\gamma^{515})|D_{i,1}|$ colours $c$ in~$D_{i,1}$, we have $i\in I_c$ and therefore $c\in D_i'$. Thus, $|D_i'|\ge (1-\gamma^{515})|D_{i,1}| \ge (1-2\gamma)q n$ by~\eqref{colour splitting sizes}. Hence, $|D_i'|=(1\pm 2\gamma)qn$, as required.

Furthermore, for all $i\in[t]$ and $c\in D_i'$, let $$M_{i,c} := A_{c,\sigma_c^{-1}(i)}.$$\COMMENT{$c\in D_i'$ implies $i\in I_c$ and hence $\sigma_c^{-1}(i)$ is defined}

Clearly, $M_{i,c}$ is a matching consisting of $256$ $c$-edges in $G_1\cup G_2$. Using~\ref{abs assignment slices}, we can also see that $V(M_{i,c})\In V_i$. Moreover, for fixed $i\in[t]$, all the matchings $\Set{M_{i,c}}_{c\in D_i'}$ are vertex-disjoint by~\ref{abs assignment disjoint}.
To check~\ref{leftover control degree}, consider any vertex $v\in V$. Clearly, the number of matchings $\set{M_{i,c}}{i\in[t],c\in D_i'}$ covering $v$ is at most~$|\cA_v|\le (1+\sqrt{\eps})pt$ by~\eqref{vertex control set}. Moreover, since $\cM$ covers all but at most $\gamma^{515}|\cA_v|$ absorbers in~$\cA_v$, we obtain a lower bound of $(1-\gamma^{515})|\cA_v|\ge (1-\gamma)pt$, as desired.

\bigskip
It remains to show the crucial property~\ref{edge absorber}. Suppose $E^\ast\In E(G_1)$ consists of precisely $m$ edges of each colour $c\in C$.
For each $c\in C$, let $E^\ast_c$ be the set of $c$-edges in~$E^\ast$.
Since $H_c$ is an RMBG with parts $[3m]$, $E_c(G_1)$ and $E_c(G_2)$, there exists a bijection $\tau_c\colon [3m] \to E^\ast_c \cup E_c(G_2)$ such that $\tau_c(j)\in N_{H_c}(j)$ for all $j\in[3m]$.

We can now define the desired partition of $E^\ast \cup E(G_2)$ as follows. Let $e\in E^\ast \cup E(G_2)$. Let $c$ be the colour of~$e$ and $j:=\tau_c^{-1}(e)$. Thus, we have $e\in N_{H_c}(j)=A_{c,j}=M_{i,c}$, where $i:=\sigma_c(j)$. Note that $i\in I_c$ and hence $c\in D_i'$. Assign $e$ to~$J_i$. Clearly, this defines a partition of $E^\ast \cup E(G_2)$ into $J_1,\dots,J_t$.
Consider $i\in[t]$. By construction, every edge $e\in J_i$ belongs to some $M_{i,c}$ with $c\in D_i'$. Moreover, for fixed $c\in D_i'$, only one edge of $M_{i,c}$ is included in~$J_i$ because $\sigma_c$ and $\tau_c$ are bijective.
\endproof

\subsection{Connecting lemma}\label{sec:connect}

The following lemma will be used to efficiently connect up the (edges from the) matchings produced by Lemmas~\ref{cor:rainbow matchings} and~\ref{lem:monochromatic matchings} of the trees~$T_i$.

Given a $k$-uniform matching~$\cR$, we say that a graph $F$ is an \defn{$\cR$-connector} if $F$ is obtained from the empty graph on $V(\cR)$ by adding, for every $R\in \cR$, new vertices $v_{R,1},\dots,v_{R,k+1}$, a perfect matching between $R$ and $\Set{v_{R,1},\dots,v_{R,k}}$ and all edges from $\Set{v_{R,1},\dots,v_{R,k}}$ to~$v_{R,k+1}$.

\begin{lemma}\label{lem:connecting}
Suppose $1/n \ll \eps \ll \gamma \ll p'\ll 1/k$ and let $p:=p'/k$ and $\beta:=q:=2p'$ and $t:=n/2$ and suppose $p''=(1\pm \eps)p'$.
Let $\phi$ be a $1$-factorization of the complete graph~$K_n$ with vertex set $V$ and colour set~$C$. Let $\tilde{G}$ be a $\beta(1+\gamma)$-random subgraph of~$K_n$. For every $i\in[t]$, let $U_i,\tilde{V_i}$ be disjoint subsets of~$V$ that are $p''$-random and $(p'+p)(1+\gamma)$-random, respectively, and let $\tilde{C_i}$ be a $q(1+\gamma)$-random subset of~$C$.

Then, the following holds with high probability:
Let $\cR$ be any $k$-uniform (multi-)hypergraph which is the union of $t$ matchings $\cR_1,\dots,\cR_t$ such that $V(\cR_i)\In U_i$ and $|\cR_i|=(1\pm \eps)p n$ for all $i\in[t]$, and such that $d_{\cR}(x) = (1\pm \eps)p't$ for all $x\in V$.
Then, for each $i\in[t]$, there exists an $\cR_i$-connector $\tilde{F}_i$ in $\tilde{G}[U_i\cup \tilde{V_i}]$ such that the following hold:
\begin{enumerate}[label=\rm{(\roman*)}]
\item $\tilde{F}_1,\dots,\tilde{F}_t$ are edge-disjoint;
\item for each $i\in[t]$, $\tilde{F}_i$ is rainbow with colours in~$\tilde{C_i}$;
\item $(\tilde{G}-\bigcup_{i\in[t]}\tilde{F}_i,\Set{(U_i\cup\tilde{V_i})\sm V(\tilde{F}_i)}_{i\in[t]},\Set{\tilde{C_i}\sm \phi(E(\tilde{F}_i))}_{i\in[t]})$ is $2\gamma n$-bounded.\label{bounded leftover connecting}
\end{enumerate}
\end{lemma}

In the proof, we will find most of the required connections via Theorem~\ref{thm:AY} (which allows us to do this `efficiently') and the remaining ones via Lemma~\ref{lem:greedy connections}.

\proof
Choose a new constant $\xi>0$ such that $\eps\ll \xi \ll \gamma$.
Split $\tilde{G}$ further into a $\beta$-random subgraph $G$ and a $\beta\gamma$-random subgraph $G'$. Moreover, for each $i\in[t]$, split $\tilde{V_i}$ into a $p'$-random subset $V_i$, a $p$-random subset~$W_i$ and a $(p'+p)\gamma$-random subset $V_i'$. Split $\tilde{C_i}$ into a $q$-random subset~$C_i$ and a $q\gamma$-random subset~$C_i'$.
We will now establish a few properties concerning the random sets which hold with high probability. From these properties, we can then (deterministically) find the desired connections for any admissible~$\cR$.

For $i\in[t]$, let $G_i$ be the spanning subgraph of $G$ with all $C_i$-edges, and let $G_i'$ be the spanning subgraph of $G'$ with all $C_i'$-edges.

For each edge $e\in E(K_n)$, let $\tilde{I}_{e,1}$ be the set of $i\in[t]$ for which $\phi(e)\in C_i$ and $e$ intersects both $U_i,V_i$, and let $I_{e,2}$ be the set of $i\in[t]$ for which $\phi(e)\in C_i$ and $e$ intersects both $V_i,W_i$.

For $i\in[t]$ and $c\in C$, let $\tilde{E}_{i,c,1}$ be the set of $c$-edges in $E_G(U_i,V_i)$, and let $E_{i,c,2}$ be the set of $c$-edges in $E_G(V_i,W_i)$.

We claim that the following hold with high probability:\COMMENT{Recall $p''=(1\pm \eps)p'$}
\begin{enumerate}[label=\rm{(\alph*)}]
\item for all $i\in[t]$, $|U_i|=(1\pm 3\eps)p'n$, $|V_i|=(1\pm \eps)p'n$, and $|W_i|=(1\pm \eps)pn$;\label{sizes}
\item for all $u\in V$, we have $|\set{i\in[t]}{u\in U_i}|=(1\pm 3\eps)p't$;\label{incidences}
\item for all $i\in[t]$ and $x\in V$, $d_{G_i}(x,U_i)= (1\pm 3\eps) p'q\beta n$, $d_{G_i}(x,V_i)= (1\pm \eps)p'q\beta n$ and $d_{G_i}(x,W_i)= (1\pm \eps)pq\beta n$;\label{degrees}
\item for all $i\in[t]$ and distinct $x,y\in V$, $|N_{G_i}(\Set{x,y})\cap V_i|=(1\pm \eps) p'q^2\beta^2  n$;\label{common neighbours}
\item for all $i\in [t]$ and $S\In V$ with $1\le |S|\le k$, we have $|N_{G_i'}(S)\cap V_i'|\ge \gamma^{2k+2} n$;\label{common neighbours reserve}
\item for all $e\in E(K_n)$, $|\tilde{I}_{e,1}|=(1\pm 2\eps)p'^2qn$ and $|I_{e,2}|=(1\pm \eps)p'pqn$;\label{edge degrees}
\item for all $i\in[t]$ and $c\in C$, $|\tilde{E}_{i,c,1}|=(1\pm 2\eps)p'^2 \beta n$ and $|E_{i,c,2}|=(1\pm \eps)p'p \beta n$;\label{colour degrees}
\item $(G',\Set{V_i'}_{i\in[t]},\Set{C_i'}_{i\in[t]})$ is $\gamma n$-bounded. \label{bounded:reserve}
\end{enumerate}
Indeed, \ref{sizes}, \ref{incidences}, \ref{degrees}, \ref{edge degrees}, \ref{colour degrees} and~\ref{bounded:reserve} follow easily from Chernoff's bound. For~\ref{common neighbours} and~\ref{common neighbours reserve}, we use McDiarmid's inequality, as follows. Consider $i\in[t]$ and distinct $x,y\in V$. Clearly, $\expn{|N_{G_i}(\Set{x,y})\cap V_i|}=p'\beta^2 q^2(n-2)$. Moreover, of the at most $2n$ edges incident with either $x$ or $y$, each has an effect of at most $1$. Each vertex has an effect of at most~$1$, and each colour has an effect of at most~$2$, and so McDiarmid's inequality applies. A similar argument works for~\ref{common neighbours reserve}.

\medskip
Now assume that \ref{sizes}--\ref{bounded:reserve} hold.
Let $\cR$ be given arbitrarily as in the lemma statement. Let $U_i':=V(\cR_i)$. By \ref{sizes} and since $|V(\cR_i)|=k|\cR_i|=(1\pm \eps)p'n$, we have that
\begin{align}
	|U_i\sm U_i'|\le 4\eps p' n \le \eps n. \label{bounded new size}
\end{align}
Moreover, for every vertex $u\in V$, it follows from~\ref{incidences} and since $d_{\cR}(u) = (1\pm \eps)p' t$ that
\begin{align}
	|\set{i\in[t]}{u\in U_i\sm U_i'}| \le 4\eps p' t \le \eps n.\label{bounded new incidences}
\end{align}
From \ref{degrees} and~\eqref{bounded new size} we infer that
\begin{align}
	\mbox{for all $i\in[t]$ and $x\in V$, $d_{G_i}(x,U_i')= (1\pm \sqrt{\eps}) p'q\beta n$.} \label{connecting back degree}
\end{align}
For an edge $e\in E(K_n)$, let $I_{e,1}$ be the set of $i\in \tilde{I}_{e,1}$ for which $e$ intersects~$U_i'$.
From~\ref{edge degrees} and~\eqref{bounded new incidences}, we deduce that
\begin{align}
	|I_{e,1}|=(1\pm \sqrt{\eps})p'^2qn. \label{new edge degrees}
\end{align}
For all $i\in[t]$ and $c\in C$, let $E_{i,c,1}$ be the set of $c$-edges in $E_G(U_i',V_i)$. By~\ref{colour degrees} and~\eqref{bounded new size}, we have
\begin{align}
	|E_{i,c,1}|=(1\pm \sqrt{\eps})p'^2 \beta n. \label{new colour degrees}
\end{align}

We now define an auxiliary hypergraph~$\cH$ whose vertex set is the union of five parts. The first part is simply $E(G)$. The second part is the set $\cR^\ast$ of all pairs $(i,R)$ such that $R\in \cR_i$. The third part is the set $\cV$ of all pairs $(i,v)$ with $v\in V_i$. The fourth part is the set $\cW$ of all pairs $(i,w)$ with $w\in W_i$. The fifth part is the set $\cC$ of all pairs $(i,c)$ with $c\in C_i$.

We now define the edge set of~$\cH$.
For disjoint $R$, $T$, $\Set{w}\In V$ and a bijection $\pi\colon T\to R$, let $S_{R,T,w,\pi}$ denote the graph on $R\cup T \cup \Set{w}$ with edge set $\set{\pi(v)v,vw}{v\in T}$. Note that $S_{R,T,w,\pi}$ is an $\{R\}$-connector.

For all $(i,R)\in \cR^\ast$, $T\In V_i$, $w\in W_i$ and bijections $\pi\colon T\to R$, we add the hyperedge
\begin{align*}
	f_{i,R,T,w,\pi}:=E(S_{R,T,w,\pi}) \cup \Set{(i,R)} \cup (\Set{i} \times T) \cup \Set{(i,w)} \cup (\Set{i}\times \phi(E(S_{R,T,w,\pi})))
\end{align*}
to $\cH$ if and only if $S_{R,T,w,\pi}$ is a rainbow subgraph of~$G_i$.
Note that $\cH$ is $(5k+2)$-uniform since $|T|=|R|=k$ and hence $E(S_{R,T,w,\pi})=2k$. We will apply Theorem~\ref{thm:AY} to~$\cH$. For this, we first check that $\cH$ is roughly regular.

For each $i\in[t]$ and $e\in E(G_i)$, let $$d_{i,e}:=|\set{f_{i,R,T,w,\pi}\in E(\cH)}{e\in E(S_{R,T,w,\pi})}|,$$
and let $d_{i,e}:=0$ for each $i\in[t]$ and $e\notin E(G_i)$.

\begin{NoHyper}
\begin{claim}\label{claim:aux edge index count}
For each $i\in[t]$ and $e\in E(G_i)$, we have $$d_{i,e}=
\begin{cases}
(1\pm \sqrt{\eps})pq\beta n (p'\beta^2 q^2 n)^{k-1} & \mbox{if $e$ intersects both $U_i'$ and $V_i$}, \\
(1\pm 3\sqrt{\eps})p'q\beta n (p'\beta^2 q^2 n)^{k-1} & \mbox{if $e$ intersects both $V_i$ and $W_i$}, \\
0 & \mbox{otherwise.}
\end{cases}$$
\end{claim}
\end{NoHyper}

\claimproof
First, assume $e=uv$ with $u\in U_i'$ and $v\in V_i$. There is a unique $R\in \cR_i$ with $u\in R$. By~\ref{degrees}, there are $(1\pm \eps)pq\beta n$ choices for~$w\in N_{G_i}(v)\cap W_i$. For each $u'\in R\sm \Set{u}$ in turn, by~\ref{common neighbours}, we have $(1\pm 2\eps) p'q^2\beta^2 n$ choices for $\pi^{-1}(u')\in N_{G_i}(u')\cap N_{G_i}(w)\cap V_i$ while avoiding previously chosen vertices and previously used colours. We deduce that
$d_{i,e}=(1\pm \sqrt{\eps})pq\beta n (p'\beta^2 q^2 n)^{k-1}$.

Next, assume $e=vw$ with $v\in V_i$ and $w\in W_i$. By~\eqref{connecting back degree}, there are $(1\pm \sqrt{\eps})p'q\beta n$ choices for $\pi(v)\in N_{G_i}(v)\cap U_i'$, which yields a unique $R\in \cR_i$ with $\pi(v)\in R$. Using~\ref{common neighbours} as above, we conclude that
$d_{i,e}=(1\pm 3\sqrt{\eps})p'q\beta n (p'\beta^2 q^2 n)^{k-1}$.

Clearly, in any other case, we have $d_{i,e}=0$.
\endclaimproof

We will use Claim~\ref*{claim:aux edge index count} below without explicit reference.

\begin{NoHyper}
\begin{claim}\label{claim:aux connect all degrees}
For all $x\in V(\cH)$, we have $d_{\cH}(x)=(1\pm \eps^{1/3})p p'^{k} \beta^{2k} q^{2k} n^{k+1}$.
\end{claim}
\end{NoHyper}

\claimproof
First, consider $e\in E(G)$. We have
\begin{eqnarray*}
	d_{\cH}(e) &=& \sum_{i\in [t]}d_{i,e} = |I_{e,1}|\cdot (1\pm \sqrt{\eps})pq\beta n (p'\beta^2 q^2 n)^{k-1} + |I_{e,2}|\cdot (1\pm 3\sqrt{\eps})p'q\beta n (p'\beta^2 q^2 n)^{k-1} \\
       &\overset{\ref{edge degrees},\eqref{new edge degrees}}{=}& (1\pm \eps^{1/3}) 2 p p'^{k+1} \beta^{2k-1} q^{2k} n^{k+1}=(1\pm \eps^{1/3})p p'^{k} \beta^{2k} q^{2k} n^{k+1},
\end{eqnarray*}
as $p'=\beta/2$.

Next, consider $(i,R)\in \cR^\ast$. By~\ref{sizes}, there are $(1\pm \eps)pn$ choices for $w\in W_i$. For each $u\in R$ in turn, by~\ref{common neighbours}, we have $(1\pm 2\eps) p'\beta^2 q^2 n$ choices for $\pi^{-1}(u)\in N_{G_i}(u)\cap N_{G_i}(w)\cap V_i$ while avoiding previously chosen vertices and previously used colours. We deduce that
$$d_{\cH}((i,R))=(1\pm \sqrt{\eps})pn(p'\beta^2 q^2 n)^{k} = (1\pm \sqrt{\eps})p p'^{k} \beta^{2k} q^{2k} n^{k+1}.$$

Now, consider $(i,v)\in \cV$. We have
\begin{eqnarray*}
	d_{\cH}((i,v)) &=& \sum_{w\in N_{G_i}(v)\cap W_i} d_{i,vw} \overset{\ref{degrees}}{=} (1\pm \eps)pq\beta n \cdot (1\pm 3\sqrt{\eps})p'q\beta n (p'\beta^2 q^2 n)^{k-1} \\
	            &=& (1\pm \eps^{1/3})p p'^{k} \beta^{2k} q^{2k} n^{k+1}.
\end{eqnarray*}

Next, consider $(i,w)\in \cW$. By assumption, we have $|\cR_i|=(1\pm \eps)pn$ choices for $R\in \cR_i$. For each $u\in R$ in turn, by~\ref{common neighbours}, we have $(1\pm 2\eps) p'\beta^2 q^2 n$ choices for $\pi^{-1}(u)\in N_{G_i}(u)\cap N_{G_i}(w)\cap V_i$ while avoiding previously chosen vertices and previously used colours. We deduce that
$$d_{\cH}(i,w)=(1\pm \eps^{1/3})pn(p'\beta^2 q^2 n)^{k}= (1\pm \eps^{1/3})p p'^{k} \beta^{2k} q^{2k} n^{k+1}.$$

Finally, consider $(i,c)\in \cC$. Note that $d_{\cH}((i,c))=\sum_{e\in E_c(G)}d_{i,e}$. Hence,
\begin{eqnarray*}
	d_{\cH}((i,c))&=& |E_{i,c,1}|\cdot (1\pm \sqrt{\eps})pq\beta n (p'\beta^2 q^2 n)^{k-1} + |E_{i,c,2}|\cdot (1\pm 3\sqrt{\eps})p'q\beta n (p'\beta^2 q^2 n)^{k-1}\\
	      &\overset{\ref{colour degrees},\eqref{new colour degrees}}{=}& (1\pm \eps^{1/3})((p'^2 \beta n)(pq\beta n) (p'\beta^2 q^2 n)^{k-1} + (p'p \beta n)(p'q\beta n) (p'\beta^2 q^2 n)^{k-1}) \\
				&=& (1\pm \eps^{1/3})2p p'^{k+1} \beta^{2k} q^{2k-1} n^{k+1}=(1\pm \eps^{1/3})p p'^{k} \beta^{2k} q^{2k} n^{k+1},
\end{eqnarray*}
since $p'=q/2$.
\endclaimproof

\begin{NoHyper}
\begin{claim}\label{claim:aux codeg 3}
$\Delta^c(\cH)\le 4k^2 n^k$.
\end{claim}
\end{NoHyper}

\claimproof
Clearly, the codegrees of pairs in $\cR^\ast \times \cR^\ast$ and $\cW \times \cW$ are~$0$.
Moreover, by Claim~\ref*{claim:aux edge index count}, we have $d_{i,e}\le n^k$ for all $i\in[t]$ and $e\in E(G_i)$. This implies that the codegrees of pairs in $E(G)\times \cR^\ast$, $E(G)\times \cV$, $E(G)\times \cW$ and $E(G)\times \cC$ are at most $n^k$, as required.
It is also easy to see that the codegrees of pairs in $\cR^\ast \times \cV$, $\cR^\ast \times \cW$, $\cV\times \cV$ and $\cV\times \cW$ are at most $n^k$, since for fixed~$i$, we always have at most $|\cR_i|\le n$ choices for $R$ and at most $n$ choices for each remaining vertex.

Consider distinct $e,e'\in E(G)$. There are $t\le n$ choices for~$i$. If $e=vw$ and $e'=v'w$ with $v,v'\in V_i$ and $w\in W_i$, then there are at most $n$ choices for $R$ and at most $n^{k-2}$ choices for $T\sm \Set{v,v'}$. Otherwise, we may assume that $e=uv$ for $u\in U_i',v\in V_i$ and $e'$ is incident to a vertex $x\in (V_i\cup W_i)\sm \Set{v}$. Now $u$ determines $R$ and there are at most $n^{k-1}$ choices for $(T\cup \Set{w})\sm \Set{v,x}$. Altogether, we conclude that the codegree of $e,e'$ is at most $n^k$.

Next, consider $(i,c),(i,c')\in \cC$ with $c\neq c'$. We have to provide an upper bound for the number of $R,T,w,\pi$ for which $f_{i,R,T,w,\pi}\in E(\cH)$ and $S_{R,T,w,\pi}$ contains a $c$-edge $e$ and a $c'$-edge $e'$. To count these possibilities, we distinguish some cases regarding how $e,e'$ intersect $U_i',V_i,W_i$. First, assume that $e,e'\in E_G(U_i',V_i)$. In this case, there are at most $n$ choices for~$R$ and then at most $k(k-1)$ choices for $e,e'$. Moreover, since $e,e'$ must form a matching, two vertices of $T$ are determined. This leaves at most $n^{k-1}$ choices for the remaining vertices, which yields a total of $k^2 n^k$ choices in this case.
Next, assume that $e,e'\in E_G(V_i,W_i)$. In this case, there are at most $n$ choices for $w$, which then determines $e$ and~$e'$ and thus two vertices from~$T$. There are at most $n$ choices for $R$ and at most $n^{k-2}$ choices for the remaining vertices of~$T$.
Finally, assume that $e\in E_G(U_i',V_i),e'\in E_G(V_i,W_i)$. We divide this case into two subcases. First, assume that $e,e'$ share their endpoint $v$ in~$V_i$. Then we have at most $n$ choices for~$v$, which determines $e,e'$, which in turn determines $R$ and $w$, and leaves at most $n^{k-1}$ choices for the vertices in $T\sm\Set{v}$. On the other hand, if $e,e'$ form a matching, then we have at most $n^2$ choices for $e,e'$, which determines $R$ and $w$ as before and leaves at most $n^{k-2}$ choices for the remaining vertices in~$T$.
Thus, altogether, the codegree of $(i,c),(i,c')$ is at most $k^2 n^k + n^k + 2n^k\le 4k^2 n^k$.

Next, consider $(i,R)\in \cR^\ast$ and $(i,c)\in \cC$. We have to choose a $c$-edge~$e$. If $e\in E_G(U_i',V_i)$, there are at most $k$ choices for $e$, which also fixes one vertex of~$T$, and leaves at most $n^{k}$ choices for the remaining vertices. If $e\in E_G(V_i,W_i)$, then there are at most $n$ choices for~$e$, which fixes $w$ and one vertex from $T$, and leaves at most $n^{k-1}$ choices for the remaining vertices. Thus, $(i,R)$ and $(i,c)$ have codegree at most $(k+1)n^k$.

Finally, consider $(i,x)\in \cV\cup \cW$ and $(i,c)\in\cC$. We have to choose a $c$-edge~$e$. If $e$ is incident with~$x$, then there is only one choice for~$e$. This either fixes $R$, in which case at most $n^{k}$ choices are left for the remaining vertices, or it fixes another vertex from $T\cup \Set{w}$, in which case there are at most $n$ choices for $R$ and at most $n^{k-1}$ choices for the remaining vertices. If $e$ is not incident with~$x$, then there are at most $n$ choices for~$e$. However, this either fixes $R$ and leaves at most $n^{k-1}$ choices for the remaining vertices, or it fixes two more vertices, which leaves at most $n$ choices for $R$ and at most $n^{k-2}$ choices for the remaining vertices. Thus, $(i,x)$ and $(i,c)$ have codegree at most $2n^k$.
\endclaimproof

For $v\in V$, let $\cR_{v}$ be the set of all $(i,R)\in \cR^\ast$ with $v\in R$, and let $\cV_v$ be the set of all pairs $(i,v)$ with $i\in [t]$ and $v\in V_i\cup W_i$. For a colour $c\in C$, let $\cC_c$ be the set of all pairs $(i,c)$ with $i\in [t]$ and $c\in C_i$. Let $$\cF:=\set{\Set{i}\times \cR_i,\Set{i}\times (V_i\cup W_i),\Set{i}\times C_i}{i\in[t]} \cup \set{\cR_v,\cV_v,\partial_G(v)}{v\in V} \cup \set{\cC_c,E_c(G)}{c\in C}.$$

Now, apply Theorem~\ref{thm:AY} to obtain a $(\xi,\cF)$-perfect matching $\cM$ in~$\cH$. For each $i\in[t]$, let $\cR_i'$ be the set of all $R\in \cR_i$ with $(i,R)\notin V(\cM)$. Since $\Set{i}\times \cR_i\in \cF$, we have $|\cR_i'|\le \xi|\cR_i| \le \xi n$. For each $R \in \cR_i\sm \cR_i'$, there is a unique edge $f_{i,R,T,w,\pi}$ in $\cM$ which covers~$(i,R)$. Let $S_{i,R}:=S_{R,T,w,\pi}$ and define $F_i:=\bigcup_{R\in \cR_i\sm \cR_i'}S_{i,R}$. By construction of~$\cH$, $F_1,\dots,F_t$ are edge-disjoint subgraphs of~$G$, and, for each $i\in [t]$, we have that $F_i$ is a rainbow $(\cR_i\sm \cR_i')$-connector with colours in~$C_i$, and $V(F_i)\cap V(\cR_i')=\emptyset$ and $V(F_i)\sm U_i' \In V_i\cup W_i$.
Moreover, observe that
\begin{align}
\mbox{$(G-\bigcup_{i\in[t]}F_i,\Set{(V_i\cup W_i)\sm V(F_i)}_{i\in[t]},\Set{C_i\sm \phi(E(F_i))}_{i\in[t]})$ is $\xi n$-bounded.}\label{bounded from PS}
\end{align}
Indeed, \ref{bounded:set sizes} holds since $\Set{i}\times (V_i\cup W_i),\Set{i}\times C_i\in \cF$ for every $i\in [t]$. Similarly, \ref{bounded:vertex incidences} holds since $\cV_v,\partial_G(v)\in \cF$ for every $v\in V$. Finally, \ref{bounded:colour incidences} holds since $\cC_c,E_c(G)\in \cF$ for every $c\in C$.

We will find the missing connectors using Lemma~\ref{lem:greedy connections}.
Let $H_i$ be an $\cR_i'$-connector. Clearly, $V(\cR_i')$ is an independent set in~$H_i$, $\Delta(H_i)\le k$ and $|V(H_i)|,|E(H_i)|\le (2k+1)|\cR_i'| \le 3k\xi n$. Moreover, for every vertex $x\in V$, the number of $i\in[t]$ for which $x\in V(\cR_i')$, is at most $|\cR_x|\le \xi n$.

Using~\ref{common neighbours reserve}, we can thus apply Lemma~\ref{lem:greedy connections} (with $G',\Set{V_i'}_{i\in[t]},\Set{C_i'}_{i\in[t]}$ playing the roles of $G,\Set{V_i}_{i\in[t]},\Set{C_i}_{i\in[t]}$) to find for each $i\in[t]$, an embedding $\psi_i\colon H_i\to G'$ such that $\psi_i(H_i)$ is rainbow with colours in~$C_i'$ and $\psi_i(x)=x$ for all $x\in V(\cR_i')$ and $\psi_i(x)\in V_i'$ for all $x\in V(H_i)\sm V(\cR_i')$, and such that $\psi_1(H_1),\dots,\psi_t(H_t)$ are edge-disjoint.

Finally, let $\tilde{F}_i:=F_i\cup \psi_i(H_i)$. Clearly, $\tilde{F}_i$ is a rainbow $\cR_i$-connector in $\tilde{G}[U_i\cup \tilde{V_i}]$ with colours in $\tilde{C_i}$, and $\tilde{F}_1,\dots,\tilde{F}_t$ are edge-disjoint.
Moreover, \ref{bounded leftover connecting} follows from \ref{bounded:reserve}, \eqref{bounded from PS}, \eqref{bounded new size} and \eqref{bounded new incidences}.
\endproof

\subsection{Rainbow perfect matchings}\label{sec:vertex}

Given a bipartite graph $G$ with vertex classes $V_1,V_2$, we say that $G$ is \defn{$(\eps,d)$-quasirandom} if for all $j\in[2]$ and distinct $v,v'\in V_j$, we have $d_G(v,V_{3-j})=(1\pm \eps) d|V_{3-j}|$ and $|N_G(\Set{v,v'})\cap V_{3-j}| = (1\pm \eps) d^2|V_{3-j}|$.

The next lemma follows easily from a result of Coulson and Perarnau~{\cite[{Lemma~6}]{CP:ta}}. Indeed, it is well known that the above notion of quasirandomness implies super-regularity, which is sufficient for the existence of a perfect matching. Moreover, using regularity, it is straightforward to count the number of `switchable edges', as required in the general statement in~\cite{CP:ta}.
\COMMENT{
[Coulson and Perarnau~{\cite[{Lemma~6}]{CP:ta}}]
Suppose $1/n\ll \mu \ll \gamma$. Let $G$ be a bipartite graph with bipartition $(A,B)$ such that $|A|=|B|=n$. Assume that $G$ has at least one perfect matching, and for every perfect matching $M$ of $G$ and for every $e\in M$ there are at least $\gamma n^2$ edges in $G$ that are $(e,M)$-switchable. Then, given any $\mu n$-bounded system of conflicts for $G$, a uniformly chosen perfect matching of $G$ is conflict-free with probability at least $\eul^{-\mu^{1/2} n}$.
}\COMMENT{It's easy to count the number of switchings in super-regular pairs. It is well known that super-regularity is implied by the above quasirandomness.}

\begin{lemma} \label{lem:switchings}
Suppose $1/n\ll \eps \ll d$. Let $G$ be a bipartite graph with vertex classes $A,B$ such that $|A|=|B|=n$ and $G$ is $(\eps,d)$-quasirandom. Then, given any edge-colouring of $G$ where each colour appears at most $\eps n$~times, there exists a rainbow perfect matching of~$G$.
\end{lemma}

We now use Lemma~\ref{lem:switchings} to obtain several edge-disjoint rainbow perfect matchings.

\begin{lemma} \label{lem:matching routine}
Suppose $1/n \ll \mu \ll d$ and let $t\le n$. Let $V$ be a vertex set of size $n$ and assume that $U_1,\dots,U_t$ are subsets of $V$\COMMENT{think of $2\mu$-random subsets} such that $|U_i|\ge \mu n$ for all $i\in[t]$, $|U_i\cap U_j|\le 5\mu^2 n$ for all distinct $i,j\in[t]$, and, for every $v\in V$, the number of $i\in[t]$ for which $v\in U_i$ is at most $3\mu t$.
For each $i\in[t]$, suppose $U_i$ is partitioned into equal-sized sets $A_i$ and~$B_i$, and $G_i$ is a $(\mu^{1/3},d)$-quasirandom bipartite graph with vertex classes $A_i,B_i$. Assume that $G_i$ is edge-coloured and each colour appears at most $2\mu^2 n$~times in~$G_i$.

Then there exist edge-disjoint $M_1,\dots,M_t$ such that $M_i$ is a rainbow perfect matching of $G_i$ for each $i\in[t]$.
\end{lemma}

We find $M_1,\dots,M_t$ using a randomised greedy algorithm.
\proof
Let $r:= \lceil 105 \mu^{3/2}n \rceil$. Suppose that we have already found $M_1,\dots,M_{s-1}$ for some $s\in[t]$. We now define $M_s$ as follows. Let $H_{s-1}:=\bigcup_{i=1}^{s-1}M_i$ and let $G_s':=G_s-H_{s-1}$. If $\Delta(H_{s-1}[U_s]) \le \mu^{3/2}n$, then $G_s'$ is $(\mu^{1/4},d)$-quasirandom. Thus, by Lemma~\ref{lem:switchings} used repeatedly, we can find edge-disjoint rainbow perfect matchings $M_{s,1},\dots,M_{s,r}$ of~$G_s'$. Otherwise, if $\Delta(H_{s-1}[U_s]) > \mu^{3/2}n$, let $M_{s,1},\dots,M_{s,r}$ be empty graphs on $U_s$. In either case, pick $j\in[r]$ uniformly at random and let $M_s:=M_{s,j}$. The lemma clearly follows if the following holds with positive probability:
\begin{align}
	\Delta(H_{s-1}[U_s]) \le \mu^{3/2}n \mbox{ for all }s\in[t].\label{eqn:Krmkey}
\end{align}
For $s \in [t]$ and $u \in U_s$, let $J^{s,u}$ be the set of indices $i\in [s-1]$ such that $u \in U_i$, so that $|J^{s,u}|\leq 3\mu t$, and for $i\in J^{s,u}$, let $Y^{s,u}_i$ be the indicator variable of the event that $uu'\in E(M_i)$ for some $u'\in U_s$.
Observe that $$d_{H_{s-1}[U_s]}(u)=\sum_{i\in J^{s,u}}Y^{s,u}_i.$$

Now, fix $s \in [t]$ and $u \in U_s$.
Crucially, for any $i\in J^{s,u}$, since $|U_s\cap U_i|\le 5\mu^2 n$, at most $5\mu^2 n$ of the matchings $M_{i,1},\dots,M_{i,r}$ that we picked in $G'_i$ contain an edge incident to~$u$ in~$G_s$ (regardless of the previous choices).
Let $i_1, \dots, i_{|J^{s,u}|}$ be the enumeration of $J^{s,u}$ in increasing order. By the above, for all $\ell \in [|J^{s,u}|] $, we have
\begin{align*}
	\prob{ Y^{s,u}_{i_\ell} = 1 \mid Y^{s,u}_{i_1}, \dots, Y^{s,u}_{i_{\ell-1}}} \le
	\frac{5\mu^2 n}{ r } \le \frac{\mu^{1/2}}{21}.
\end{align*}
Let $B \sim Bin( |J^{s,u}|  , \mu^{1/2}/21 )$. Since $ |J^{s,u}|  \le 3 \mu n$, we have $\expn{B}\le \mu^{3/2}n/7$.
Using Fact~\ref{prop:jain} and Lemma~\ref{lem:chernoff}\ref{chernoff crude}, we infer that
\begin{align*}
\prob{ \sum_{i \in J^{s,u}} Y^{s,u}_i  >  \mu^{3/2} n }
\le \prob{ B >  \mu^{3/2} n } \le e^{-\mu^{3/2} n}.
\end{align*}
Finally, a union bound implies that~\eqref{eqn:Krmkey} holds with high probability.
\endproof

\subsection{Proof of Theorem~\ref{thm:main}}\label{sec:main}

We are now ready to prove our main theorem.

\lateproof{Theorem~\ref{thm:main}}
Choose new constants $\eps ,\gamma ,\xi, \mu, \eta>0$ such that $$1/n \ll \eps \ll \gamma \ll \xi \ll  \mu \ll \eta \ll 1,$$
and let
$$t:=n/2 \hspace{15pt}   r:=\lceil (\eta /256+ 6\rhosub +3\gamma) n \rceil         \hspace{15pt}    b:= \lceil (\mu-\xi^{1/3}) n\rceil.$$
Let $\phi$ be a $1$-factorization of~$K_n$ with vertex set~$V$ and colour set~$C$. We will obtain a decomposition into $t$ rainbow copies of $T_{n;r,b}$ (cf.~Definition~\ref{def:tree}). Hence, $r$ and $b$ are essentially determined by $\eta$ and $\mu$, respectively, and $\eps,\gamma,\xi$ are best thought of as error parameters.

In order to apply the lemmas that we have proven in this section without interference, we will split $E(K_n)$, $V$ and $C$ into random subsets each reserved for the application of the relevant lemma. For convenience, we now define the relevant constants in one place (where the letters $p,q,\beta$ represent vertex, colour and edge probabilities, respectively).
\begin{center}
$p_{\rm{rb}}:=2\eta$    \hspace{15pt}  $q_{\rm{rb}} :=\eta/192$ \hspace{15pt}     $p_{\rm{mc}}:=3072\rhosub$  \hspace{15pt}  $q_{\rm{mc}}:=6\rhosub$  \\

\vspace{5pt}

 $m:=\lceil (\rhosub-\eps/5) n \rceil$    \hspace{15pt}      $s:=\lceil (q_{\rm{rb}}/4- 2\gamma/5) n\rceil$

\vspace{5pt}

$\tilde{p}':=p_{\rm{rb}} + p_{\rm{mc}}$  \hspace{15pt}  $\tilde{p}:=\eta/256+6\rhosub$  \hspace{15pt}   $\tilde{\beta}:=\tilde{q}:=2\tilde{p}'$\\

\vspace{5pt}

$p_{\circ}:= 1  - \tilde{p}'(1+\gamma) - (\tilde{p}'+\tilde{p})(1+\xi) - 2\mu$ \\

\vspace{5pt}

$q_{\circ,1}:= 1  - q_{\rm{rb}} - q_{\rm{mc}}(1+\gamma) - \tilde{q}(1+\xi)  - \mu$   \hspace{15pt}   $q_{\circ,2}:= p_{\circ} - q_{\circ,1} $    \\

\vspace{5pt}

$\beta_{\circ,1}:= 1-8\rhosub- \tilde{\beta}(1+\xi)-\mu$  \hspace{15pt}   $\beta_{\circ,2} :=  p_{\circ} -\beta_{\circ,1} $   \\

\vspace{5pt}

$q_\vartriangle:= (q_{\rm{rb}}/2 - q_{\circ,2})/3$  \hspace{15pt}  $\beta_\vartriangle := (4\rhosub-\eta(1+\gamma)-\beta_{\circ,2})/3 $.
\end{center}

Note that, as $\gamma\ll\xi$, $q_{\circ,2} = \eta/768-\mu \pm \xi$\COMMENT{$q_{\circ,2} = p_{\circ} - q_{\circ,1} = (1  - \tilde{p}'(1+\gamma) - (\tilde{p}'+\tilde{p})(1+\xi) - 2\mu) - (1  - q_{\rm{rb}} - q_{\rm{mc}}(1+\gamma) - \tilde{q}(1+\xi)  - \mu) = - 2\tilde{p}' - \tilde{p} + q_{\rm{rb}} + q_{\rm{mc}} +\tilde{q} -\mu \pm \xi = \eta/768-\mu \pm \xi$ since $2\tilde{p}'=\tilde{q}$ and $- \tilde{p} + q_{\rm{rb}} + q_{\rm{mc}}= -\eta/256-6\rhosub +\eta/192+6\rhosub = \eta/768$} and $\beta_{\circ,2} = 2\rhosub -\eta/256 -\mu \pm \xi$\COMMENT{$\beta_{\circ,2} = p_{\circ} - \beta_{\circ,1} = (1  - \tilde{p}'(1+\gamma) - (\tilde{p}'+\tilde{p})(1+\xi) - 2\mu) - (1-8\rhosub- \tilde{\beta}(1+\xi)-\mu) = - 2\tilde{p}' - \tilde{p} + 8\rhosub +\tilde{\beta} -\mu \pm \xi = 2\rhosub -\eta/256 -\mu \pm \xi$ since $2\tilde{p}'=\tilde{\beta}$ and $-\tilde{p}+8\rhosub = 2\rhosub-\eta/256$} and hence $q_\vartriangle\ge \eta/2304$\COMMENT{$q_{\rm{rb}}/2 - q_{\circ,2} \ge \eta/384 - \eta/768 = \eta/768$} and $\beta_\vartriangle \ge \rhosub/3$.

\medskip

\begin{NoHyper}
\begin{step}
Random splitting
\end{step}
\end{NoHyper}

\noindent\textbf{Split vertices.} For each $i\in[t]$, we split $V$ randomly as follows:
$$\begin{array}{cccccccccccccccccccccccccccccc}
  V   &=&  U_i  &\cupdot&	   \tilde{V}_i 	&\cupdot&   V^{\circ}_i   &\cupdot&  A_i    &\cupdot&   B_i\\
	1   &=&  \tilde{p}'(1+\gamma)  &+& (\tilde{p}'+\tilde{p})(1+\xi) &+& p_{\circ} &+& \mu  &+&  \mu.
\end{array}$$
We split $U_i$ and $B_i$ further as follows:

\begin{minipage}[b]{0.45\textwidth}
$$\begin{array}{cccccccccccccccccccccccccccccc}
  U_i   &=&  V^{\rm{rb}}_i &\cupdot& V^{\rm{mc}}_i \\
	\tilde{p}'(1+\gamma)  &=&  p_{\rm{rb}}(1+\gamma) &+&  p_{\rm{mc}}(1+\gamma)
\end{array}\;\;\;\;\;\;\;\text{ and }$$
\end{minipage}
\hfill
\begin{minipage}[b]{0.45\textwidth}
$$\begin{array}{cccccccccccccccccccccccccccccc}
  B_i   &=&  B_{i,1} &\cupdot& B_{i,2} \\
	\mu  &=&  \mu/2 &+&  \mu/2.
\end{array}$$
\end{minipage}

\medskip
\noindent\textbf{Split colours.} Moreover, for each $i\in[t]$, we split $C$ randomly as follows:
$$\begin{array}{cccccccccccccccccccccccccccccc}
  C   &=&  C_{i,1} &\cupdot& C_{i,2}  &\cupdot&	   D_i    &\cupdot&	 \tilde{C}_i 	&\cupdot&   C^\bullet_i   &\cupdot&  C^{\circ,1}_i\\
	1   &=&  q_{\rm{rb}}/2 &+&  q_{\rm{rb}}/2  &+& q_{\rm{mc}}(1+\gamma) &+& \tilde{q}(1+\xi) &+& \mu  &+& q_{\circ,1}.
\end{array}$$
We split $C_{i,1}$ further as follows:
$$\begin{array}{cccccccccccccccccccccccccccccc}
  C_{i,1}   &=& C_{i,1}^\vartriangle &\cupdot& C_{i,2}^\vartriangle &\cupdot& C_{i,3}^\vartriangle &\cupdot&  C^{\circ,2}_i\\
	q_{\rm{rb}}/2   &=&  q_{\vartriangle} &+& q_{\vartriangle} &+& q_{\vartriangle} &+& q_{\circ,2}.
\end{array}$$
Let $C^{\circ}_i:=C^{\circ,1}_i \cup C^{\circ,2}_i$. Hence, $C^{\circ}_i$ is a $p_{\circ}$-random set. Moreover, let $C_i^\vartriangle:=C_{i,1}^\vartriangle\cup C_{i,2}^\vartriangle \cup C_{i,3}^\vartriangle$.

\medskip
\noindent\textbf{Split edges.} We split $K_n$ randomly as follows:
$$\begin{array}{cccccccccccccccccccccccccccccc}
  K_n &=&   G_1 &\cupdot& G_2 &\cupdot& \tilde{G} &\cupdot&  G^\bullet  &\cupdot& G^{\circ,1} \\
	1      &=&  4\rhosub &+&  4\rhosub  &+& \tilde{\beta}(1+\xi)  &+& \mu &+& \beta_{\circ,1}.
\end{array}$$
Split $G_1$ further as follows:
$$\begin{array}{cccccccccccccccccccccccccccccccc}
  G_1    &=&   G^{\rm{rb}} &\cupdot& G^\vartriangle_1 &\cupdot& G^\vartriangle_2 &\cupdot& G^\vartriangle_3 &\cupdot&G^{\circ,2} \\
	4\rhosub   &=&  \eta(1+\gamma) &+&  \beta_\vartriangle  &+&  \beta_\vartriangle  &+&   \beta_\vartriangle  &+&  \beta_{\circ,2}.
\end{array}$$
Let $G^{\circ}:=G^{\circ,1} \cup G^{\circ,2}$. Thus, $G^{\circ}$ is a $p_{\circ}$-random subgraph. Moreover, let $G^\vartriangle:=G^\vartriangle_1 \cup G^\vartriangle_2  \cup G^\vartriangle_3 $.

\bigskip

\noindent\textbf{Create the edge reservoir.}
By Lemma~\ref{lem:monochromatic matchings} (with $G_1,G_2,G^{\rm{rb}} \cup G^\vartriangle,G^{\circ,2},\Set{V^{\rm{mc}}_i,D_i}_{i\in[t]}$ in place of $G_1',G_2',G_{1,1}',G_{1,2}',\Set{V_i,D_i}_{i\in[t]}$), with high probability,\COMMENT{$G^{\rm{rb}} \cup G^\vartriangle$ is admissible choice for $G_{1,1}'$ since the `gap' probability is $\beta_{\circ,2} \ge \rhosub$} there exist $G_1',G_2'$ such that $$G^{\rm{rb}} \cup G^\vartriangle \In G_1' \In G_1\mbox{ and }G_2'\In G_2$$ with $\Delta(G_2-G_2')\le 2\eps n$ and,
for each $i\in[t]$, there exists $D_i'\In D_i$ of size $(1\pm 2\gamma)q_{\rm{mc}} n$ and vertex-disjoint matchings $\Set{M_{i,c}'}_{c\in D_i'}$ in $(G_1'\cup G_2')[V^{\rm{mc}}_i]$, where $M_{i,c}'$ consists of $256$~$c$-edges, such that altogether the following hold:
\begin{enumerate}[label=\rm{(M\arabic*)}]
\item $|E_c(G_1')|=|E_c(G_2')|=2m $ and $|\set{i\in[t]}{c\in D_i'}|=3m$ for all $c\in C$; \label{colour regularity app}
\item for any subset $E^\ast\In E(G_1')$ which consists of precisely $m$ edges of each colour $c\in C$, there exists a partition of $E^\ast \cup E(G_2')$ into sets $J_1',\dots,J_t'$, such that for each $i\in[t]$, $J_i'$ contains exactly one edge from each of $\Set{M_{i,c}'}_{c\in D_i'}$; \label{edge absorber app}
\item every vertex $v\in V$ is covered by $(1\pm \gamma)p_{\rm{mc}}t$ of the matchings $\set{M_{i,c}'}{i\in[t],c\in D_i'}$.\label{leftover control degree app}
\end{enumerate}

\noindent\textbf{Create colour reservoirs.}
We apply Lemma~\ref{cor:rainbow matchings} (with $G^{\rm{rb}},\Set{V^{\rm{rb}}_i,C_{i,1},C_{i,2},C_i^\vartriangle,C^{\circ,2}_i}_{i\in[t]}$ in place of $G,\Set{V_i,C_{i,1},C_{i,2},C_{i,1,1},C_{i,1,2}}_{i\in[t]}$) to see that with high probability,\COMMENT{$C_i^\vartriangle$ is admissible choice for $C_{i,1,1}$ since the `gap' probability is $q_{\circ,2}\ge \eta/800$} for each $i\in[t]$, there exist $$C_i^\vartriangle \In C_{i,1}' \In C_{i,1}\mbox{ and } C_{i,2}'\In C_{i,2}$$ and
vertex-disjoint rainbow matchings $\set{M_{i,j}}{j\in[3s]}$ in $G^{\rm{rb}}[V^{\rm{rb}}_i]$, such that altogether the following hold:
\begin{enumerate}[label=\rm{(R\arabic*)}]
\item for each $i\in[t]$, $|C_{i,1}'|=|C_{i,2}'|=2s$;\label{rainbow matchings set sizes}
\item for all $c\in C$, $|\set{i\in[t]}{c\in C_{i,2}\sm C_{i,2}'}|\le \sqrt{\gamma} n$;\label{rainbow matchings bounded app}
\item for each $i\in[t]$, $M_i:=\bigcup_{j\in[3s]}M_{i,j}$ consists of $192$ $c$-edges for each $c\in C_{i,1}'\cup C_{i,2}'$; \label{rainbow matchings partition app}
\item for each $i\in [t]$ and any subset $C_i^\ast \In C_{i,1}'$ of size~$s$, there exists $J_i\In M_i$ such that $J_i$ is $(C_i^\ast \cup C_{i,2}')$-rainbow and contains exactly one edge from each of $\set{M_{i,j}}{j\in[3s]}$; \label{colour absorber app}
\item the matchings $\set{M_{i,j}}{(i,j)\in [t]\times [3s]}$ are edge-disjoint, and $|M_{i,j}|=256$ for all $(i,j)\in [t]\times [3s]$; \label{rainbow matchings partition 2 app}
\item for every vertex $v\in V$, the number of $i\in[t]$ for which $v$ is covered by $M_i$ is $(1\pm \sqrt{\gamma})p_{\rm{rb}}t$. \label{degree control vertices rainbow app}
\end{enumerate}

\noindent\textbf{Create short paths for the vertex absorption.}
By Lemma~\ref{lem:approximate dec} (with $\eps/5,2\xi^{1/3},\mu$ playing the roles of $\gamma,\kappa,p$), with high probability, there exist edge-disjoint rainbow paths $Q_1,\dots,Q_t$ in $G^\bullet$ such that
\begin{enumerate}[label=\rm{(Q\arabic*)}]
\item for each $i\in[t]$, we have $V(Q_i)\In A_i$ and $\phi(E(Q_i))\In C^\bullet_i$;
\item $(G^\bullet-\bigcup_{i\in[t]}Q_i$, $\Set{A_i\sm V(Q_i)}_{i\in[t]}$, $\Set{C^\bullet_i \sm \phi(E(Q_i))}_{i\in[t]})$ is $\eps n$-bounded;\label{bounded leftover approximate app 2}
\item for each $v\in V$, the number of $i\in[t]$ for which $v\in V(Q_i)$ and the subpath from $v$ to one of the endvertices of $Q_i$ has length at most $2\xi^{1/3} n$, is at most $\xi^{1/4} n$.\label{path ends Q}
\end{enumerate}

\smallskip

\noindent\textbf{Properties for vertex absorption and covering non-reservoir edges/colours/vertices.}
For $i\in[t]$, let $G^{(i)}$ be the subgraph of $G^\vartriangle_3[A_i,B_i]$ containing precisely the $C_{i,3}^\vartriangle$-edges.
In addition to the above, with high probability, the following hold:
\begin{enumerate}[label=\rm{(A\arabic*)}]
\item for all $i\in[t]$ and $c\in C$, there are at most $2\mu^2 n$ $c$-edges between $A_i$ and $B_i$;\COMMENT{there are $n/2$ $c$-edges each of which is between $A_i$ and $B_i$ with probability $2\mu^2$}\label{colour boundedness}
\item for all distinct $i,i'\in [t]$, we have $|(A_i\cup B_i)\cap (A_{i'}\cup B_{i'})| \le 5\mu^2 n$;\COMMENT{expectation is $4\mu^2 n$}\label{jain:intersections}
\item for every $v\in V$, the number of $i\in[t]$ for which $v\in A_i\cup B_i$, is at most~$3\mu t$;\COMMENT{expectation is $2\mu n/2$}\label{jain:incidences}
\item for all $i\in[t]$, $|A_i|,|B_i|=(1\pm \eps)\mu n$ and $G^{(i)}$ is $(\eps,\beta_\vartriangle q_\vartriangle)$-quasirandom;\label{super regular}
\item for all $i\in[t]$ and all $S\In V$ with $|S|\le 512$, we have that $$|\set{v\in N_{G^\vartriangle_1}(S)\cap B_{i,1}}{\phi(uv)\in C_{i,1}^\vartriangle \mbox{ for each }u\in S}| \ge \mu^2 n;$$\COMMENT{expected number is $(\beta_\vartriangle q_\vartriangle)^{|S|}\mu n/2 \ge 2\mu^2 n$}\label{greedy connection assumption}
\item for all $e\in E(K_n)$, the number of $i\in[t]$ for which $e$ intersects $A_i$ and $B_{i,2}$, and $\phi(e)\in C_{i,2}^\vartriangle$, is at least~$q_\vartriangle \mu^2 n/3$;\COMMENT{for fixed $i$, the probability is $2(\mu^2/2)q_\vartriangle$}\label{edge cover down}
\item for all $i\in[t]$ and $c\in C$, the number of $c$-edges in $E_{G^\vartriangle_2}(A_i,B_{i,2})$, is at least~$\beta_\vartriangle \mu^2 n/3$;\COMMENT{for fixed edge, the probability is $2(\mu^2/2)\beta_\vartriangle$, $n/2$ $c$-edges, all independent since they form a matching}\label{colour cover down}
\item for all $i\in[t]$, $|C_{i,2}|=(1\pm \eps)q_{\rm{rb}}n/2$ and $|D_i|=(1\pm 2\gamma)q_{\rm{mc}}n$;\label{misc:colour set size}
\item for all $c\in C$, $|E_c(G_2)|=(1\pm \eps)2\rhosub n$ and $|\set{i\in[t]}{c\in D_i}|=(1\pm \sqrt{\gamma})q_{\rm{mc}}t$.\label{misc:colour set incidences}
\end{enumerate}
Here, to deal with the sizes of the common neighbourhoods in \ref{super regular} and~\ref{greedy connection assumption}, we use McDiarmid's inequality. For all other claims, Chernoff's bound suffices.

\smallskip

\noindent\textbf{Find almost-spanning paths.} We now find an approximate decomposition of $G^\circ$ (and thus~$K_n$) into $t=n/2$ almost spanning rainbow paths. By Lemma~\ref{lem:approximate dec} (with $\eps/5,\eps^2,p_{\circ}$ playing the roles of $\gamma,\kappa,p$), with high probability, there exist edge-disjoint rainbow paths $P_1,\dots,P_t$ in $G^{\circ}$ such that
\begin{enumerate}[label=\rm{(P\arabic*)}]
\item for each $i\in[t]$, we have $V(P_i)\In V^{\circ}_i$ and $\phi(E(P_i))\In C^{\circ}_i$;
\item $(G^{\circ}-\bigcup_{i\in[t]}P_i$, $\Set{V^{\circ}_i\sm V(P_i)}_{i\in[t]}$, $\Set{C^{\circ}_i \sm \phi(E(P_i))}_{i\in[t]})$ is $\eps n$-bounded;\label{bounded leftover approximate app 1}
\item every vertex $v\in V$ is an endvertex of at most $\eps n$ paths.\label{path ends P}
\end{enumerate}

\smallskip

\noindent\textbf{Establish connection properties.} By Lemma~\ref{lem:connecting} (with $\gamma^{1/3},\xi,\tilde{p}',\tilde{p},512$ playing the roles of $\eps,\gamma,p',p,k$), with high probability, the following is true:
\begin{enumerate}[label=\rm{(C)}]
\item Let $\cR$ be any $512$-uniform (multi-)hypergraph which is the union of $t$ matchings $\cR_1,\dots,\cR_t$ such that $V(\cR_i)\In U_i$ and $|\cR_i|=(1\pm \gamma^{1/3})\tilde{p} n$ for all $i\in[t]$, and such that $d_{\cR}(x) = (1\pm \gamma^{1/3})\tilde{p}'t$ for all $x\in V$.
Then, for each $i\in[t]$, there exists an $\cR_i$-connector $\tilde{F}_i$ in $\tilde{G}[U_i\cup \tilde{V}_i]$ such that the following hold:
\begin{enumerate}[label=\rm{(C\arabic*)}]
\item $\tilde{F}_1,\dots,\tilde{F}_t$ are edge-disjoint;
\item for each $i\in[t]$, $\tilde{F}_i$ is rainbow with colours in~$\tilde{C}_i$;
\item $(\tilde{G}-\bigcup_{i\in[t]}\tilde{F}_i,\Set{(U_i\cup\tilde{V_i})\sm V(\tilde{F}_i)}_{i\in[t]},\Set{\tilde{C_i}\sm \phi(E(\tilde{F}_i))}_{i\in[t]})$ is $2\xi n$-bounded.
\end{enumerate}
\label{connecting}
\end{enumerate}

\bigskip
{\bf Henceforth, we assume that all random choices have been made and satisfy the above properties.}

\begin{NoHyper}
\begin{step}\label{step:connecting}
Connecting the pieces
\end{step}
\end{NoHyper}

For each $i\in[t]$, we now aim to connect the matchings $\set{M_{i,j}}{j\in[3s]}$ and $\set{M_{i,c}'}{c\in D_i'}$. For this, we will define a $512$-uniform matching $\cR_i$, which consists of $256$ vertices of one matching and $256$ vertices of the next matching, and then apply~\ref{connecting}.

To make this more precise, for each $i\in [t]$, let
\begin{align}
\mbox{$r_i:=3s+|D_i'|$ and $\cM_i:=\set{M_{i,j}}{j\in[3s]} \cup \set{M_{i,c}'}{c\in D_i'}$.}\label{def number absorbers}
\end{align}
So $r_i=|\cM_i|$. Note that $r_i= \tilde{p} n \pm 3\gamma n$ and hence $0\le r-r_i \le 7\gamma n$. Also note that since $V^{\rm{rb}}_i$ and $V^{\rm{mc}}_i$ are disjoint, all the matchings in $\cM_i$ are vertex-disjoint, and recall that each matching consists of $256$~edges.
For each $i\in[t]$, find two distinct $M_i^-,M_i^+\in \cM_i$ such that altogether,
\begin{align}
	\mbox{each vertex $v\in V$ is contained in $V(M_i^-\cup M_i^+)$ for at most $\eta^{-2}$ indices $i\in[t]$.}\label{start end matchings}
\end{align}
This can clearly be done greedily.\COMMENT{When choosing $M_i^-,M_i^+$, the number of vertices which are already used at least $\eta^{-2}/2$ times is at most $\eta^{2}1024n$, and thus at most that many matchings of $\cM_i$ are blocked. Since $|\cM_i| \ge 3s \ge \eta n/800$, there are two available matchings that we can choose.}
Now, for each $i\in[t]$, choose an arbitrary bijection $\sigma_i\colon [r_i]\to \cM_i$ such that $\sigma_i(1)=M_i^-$ and $\sigma_i(r_i)=M_i^+$, and partition for each matching $M\in \cM_i$ the vertices $V(M)$ arbitrarily into a `tail set' $T(M)$ and a `head set' $H(M)$ such that $M$ is a perfect matching between $T(M)$ and~$H(M)$.
Define $$\cR_i:=\set{H(\sigma_i(k)) \cup T(\sigma_i(k+1))}{k\in[r_i-1]}.$$
Hence, $\cR_i$ is a $512$-uniform matching in~$U_i$.
Note that $$|\cR_i|=r_i-1 = (1\pm \sqrt{\gamma})\tilde{p} n.$$
Let $\cR:=\cR_1\cup \dots \cup \cR_t$. By \ref{degree control vertices rainbow app}, \ref{leftover control degree app} and~\eqref{start end matchings}, we have that $d_{\cR}(x) = (1\pm 2\sqrt{\gamma})(p_{\rm{rb}} + p_{\rm{mc}}) t = (1\pm 2\sqrt{\gamma})\tilde{p}' t$ for all $x\in V$.
Hence, applying~\ref{connecting}, for each $i\in[t]$, there exists an $\cR_i$-connector $\tilde{F}_i$ in $\tilde{G}[U_i\cup \tilde{V}_i]$ such that altogether the following hold:
\begin{enumerate}[label=\rm{(C\arabic*$'$)}]
\item $\tilde{F}_1,\dots,\tilde{F}_t$ are edge-disjoint;
\item for each $i\in[t]$, $\tilde{F}_i$ is rainbow with colours in~$\tilde{C}_i$;
\item $(\tilde{G}-\bigcup_{i\in[t]}\tilde{F}_i,\Set{(U_i\cup\tilde{V_i})\sm V(\tilde{F}_i)}_{i\in[t]},\Set{\tilde{C_i}\sm \phi(E(\tilde{F}_i))}_{i\in[t]})$ is $2\xi n$-bounded.\label{bounded leftover connecting app}
\end{enumerate}

For each $i\in[t]$, let $$F_i := P_i \cup Q_i \cup \tilde{F}_i.$$
We will eventually have $F_i\In T_i$ for all $i\in[t]$. Note that $F_1,\dots,F_t$ are edge-disjoint rainbow forests in $K_n-(G^{\rm{rb}} \cup G^\vartriangle \cup G_2)$. Moreover, for each $i\in[t]$, $V(F_i)\In V\sm B_i$ and $\phi(E(F_i))\In C\sm (C_i^\vartriangle \cup C_{i,2} \cup D_i)$.
Let
\begin{align}
	\hat{V}_i &:= V\sm (V(F_i)\cup B_i), \\
	\hat{C}_i &:= C\sm (\phi(E(F_i)) \cup C_{i,1}'\cup C_{i,2}' \cup D_i'),\label{leftover colours} \\
	\hat{G} &:= K_n  - \bigcup_{i\in[t]}F_i - (G_1'\cup G_2').\label{leftover graph}
\end{align}

We think of the above as leftover sets. The following claim asserts that this leftover is well-behaved.

\begin{NoHyper}
\begin{claim}\label{claim:bounded leftover}
$(\hat{G}, \Set{\hat{V}_i}_{i\in[t]}, \Set{\hat{C}_i}_{i\in[t]})$ is $\sqrt{\xi} n$-bounded.
\end{claim}
\end{NoHyper}

\claimproof
Observe that
\begin{align*}
	\hat{V_i} &=((U_i\cup \tilde{V}_i )\sm V(\tilde{F}_i)) \cup (V^{\circ}_i \sm V(P_i)) \cup (A_i \sm V(Q_i)),\\
	\hat{C_i} &\In ((C_{i,2}\cup D_i)\sm (C_{i,2}'\cup D_i')) \cup (\tilde{C}_i \sm \phi(E(\tilde{F}_i))) \cup (C^{\circ}_i \sm \phi(E(P_i))) \cup (C^\bullet_i \sm \phi(E(Q_i))),\\
	\hat{G} &\In (G_2-G_2') \cup \bigg(\tilde{G}-\bigcup_{i\in[t]}\tilde{F}_i\bigg) \cup \bigg(G^{\circ}-\bigcup_{i\in[t]}P_i\bigg) \cup \bigg(G^\bullet-\bigcup_{i\in[t]}Q_i\bigg).
\end{align*}
\COMMENT{equality holds with $C^{\circ}_i \sm \phi(E(P_i))$ replaced with $C^{\circ}_i \sm (\phi(E(P_i))\cup C_{i,1}')$ and $G^{\circ}-\bigcup_{i\in[t]}P_i$ replaced with $G^{\circ}-(G_1'\cup \bigcup_{i\in[t]}P_i)$}
Recall that $\Delta(G_2-G_2')\le 2\eps n$ and $|D_i'|=(1\pm 2\gamma)q_{\rm{mc}} n$ for all $i\in[t]$. Thus, \ref{rainbow matchings set sizes} and \ref{misc:colour set size} imply that, for all $i\in[t]$,
 $$|(C_{i,2}\cup D_i)\sm (C_{i,2}'\cup D_i')|\le \gamma^{1/3} n$$
 and \ref{rainbow matchings bounded app}, \ref{colour regularity app}  and~\ref{misc:colour set incidences} imply that,  for all $c\in C$, $$|E_c(G_2-G_2')|,|\set{i\in[t]}{c\in (C_{i,2}\cup D_i)\sm (C_{i,2}'\cup D_i')}| \le \gamma^{1/3} n.$$

Hence, the claim follows together with \ref{bounded leftover approximate app 1}, \ref{bounded leftover approximate app 2} and~\ref{bounded leftover connecting app}.
\endclaimproof

We now use Lemma~\ref{lem:greedy connections} to join the pieces of each $F_i$ together. Moreover, since the sets $\cR_i$ have different sizes, we artificially add some structure that will ensure that ultimately, all trees are isomorphic to~$T$ (cf.~$(\dagger)$ below). In this process we can cover all remaining vertices outside the vertex reservoir~$B_i$.

For $i\in[t]$, let $v_i^-,v_i^+$ be the endvertices of~$P_i$, and let $w_i^-,w_i^+$ be the endvertices of~$Q_i$.
Let $$X_i:=T(M_i^-)\cup H(M_i^+) \cup \Set{v_i^-,v_i^+,w_i^-} \cup \hat{V}_i.$$
We now define a graph $H_i$ in which $X_i$ is independent and all other vertices are new vertices. Take new vertices $z_0,z_{r_i},\dots,z_r$. For each $x\in T(M_i^-)$, add a path of length~$2$ between $x$ and~$z_0$, and for each $x\in H(M_i^+)$, add a path of length~$2$ between $x$ and~$z_{r_i}$. For each $k\in\Set{r_i,\dots,r-1}$, add a path of length~$5$ between $z_k$ and $z_{k+1}$. For each $k\in\Set{r_i+1,\dots,r-1}$, add $510$ further paths of length~$2$ onto~$z_k$ (so $z_k$ will be an endvertex of those paths of length~$2$), and add $255$ paths of length~$2$ onto each of $z_{r_i}$ and~$z_r$. Connect $z_r$ and $v_i^-$ by a path which contains $\hat{V}_i$ such that $\hat{V}_i$ is an independent set.
Finally, add a path of length~$2$ between $v_i^+$ and~$w_i^-$. Clearly, $\Delta(H_i)\le 512$.
Note that since $|\hat{V}_i|\le \sqrt{\xi} n$, $r-r_i\le 7\gamma n$ and $|B_i|-b = \xi^{1/3}n\pm \eps n$ by~\ref{super regular}, we can choose $H_i$ in such a way that
\begin{align}
	 |V(H_i)\sm X_i|=|B_i|-b \mbox{ and }|V(H_i)|,|E(H_i)|\le  2\xi^{1/3}n.\label{leftover vertices precise}
\end{align}
Also, for every $v\in V$, by~\eqref{start end matchings}, \ref{path ends P}, \ref{path ends Q} and Claim~\ref*{claim:bounded leftover}, the number of $i\in[t]$ for which $v\in X_i$, is at most~$2\xi^{1/4}n$.

By~\ref{greedy connection assumption}, we can now apply Lemma~\ref{lem:greedy connections} (with $G^\vartriangle_1,\Set{C_{i,1}^\vartriangle}_{i\in[t]},\Set{B_{i,1}}_{i\in[t]}$ taking the place of $G$,$\Set{C_i}_{i\in[t]}$,$\Set{V_i}_{i\in[t]}$) to obtain, for each $i\in[t]$, an embedding $\psi_i\colon H_i\to G^\vartriangle_1$ such that $\psi_i(H_i)$ is rainbow with colours in~$C_{i,1}^\vartriangle$, $\psi_i(x)=x$ for all $x\in X_i$, $\psi_i(x)\in B_{i,1}$ for all $x\in V(H_i)\sm X_i$ and such that $\psi_1(H_1),\dots,\psi_t(H_t)$ are edge-disjoint.

For each $i\in[t]$, let $$F_i^\ast:=F_i\cup \psi_i(H_i)\mbox{ and }B_i^\ast:=V \sm V(F_i^\ast)\In B_i.$$
Note that $F_i^\ast$ is rainbow as $C_{i,1}^\vartriangle\In C_{i}^\vartriangle$.
By~\eqref{leftover vertices precise}, we have that $|B_i^\ast|=b$. Moreover, $B_{i,2}\In B_i^\ast$. Let $A_i^\ast$ be the set of the last $b$~vertices on~$Q_i$, containing~$w_i^+$, so that $A_i^\ast\In A_i$.
From~\ref{bounded leftover approximate app 2}, \ref{path ends Q} and~\ref{super regular}, we deduce that, for each $i\in[t]$,
\begin{align}
	\mbox{$|A_i\sm A_i^\ast|\le |A_i|-b \le 2\xi^{1/3}n$,}\label{final A set size}
\end{align}
and, for each $v\in V$,
\begin{align}
	\mbox{$|\set{i\in[t]}{v\in A_i\sm A_i^\ast}| \le 2\xi^{1/4} n$.}\label{final A set incidences}
\end{align}

Crucially, observe that for each $i\in[t]$,

\medskip

\noindent
\begin{minipage}{0.05\textwidth}
$(\dagger)$
	\end{minipage}
	\hfill
	\begin{minipage}{0.95\textwidth} any graph obtained from $F_i^\ast$ by adding a perfect matching between $A_i^\ast$ and $B_i^\ast$ and exactly one edge from each of the matchings in~$\cM_i$, is isomorphic to~$T_{n;r,b}$.
\end{minipage}

\medskip

In particular,
\begin{align}
	|E(F_i^\ast)|=n-1-b-r_i.\label{missing edge count}
\end{align}

\begin{NoHyper}
\begin{step}
Final absorption
\end{step}
\end{NoHyper}

We will find the perfect matchings between $A_i^\ast$ and $B_i^\ast$ using Lemma~\ref{lem:matching routine}, and then select exactly one edge from each of the matchings in~$\cM_i$ using \ref{colour absorber app} and~\ref{edge absorber app}. For the last step to work, we need to ensure that all leftover colours are in $C_{i,1}'\cup C_{i,2}' \cup D_i'$ and all leftover edges are in $G_1'\cup G_2'$. Thus, prior to applying Lemma~\ref{lem:matching routine}, we greedily deal with the colours in $\hat{C}_i$ and the edges in $\hat{G}$.

\smallskip

\noindent\textbf{Cover the remaining non-reservoir edges.} First, find a partition of $E(\hat{G})$ into rainbow matchings $\hat{M}_1,\dots,\hat{M}_t$ such that $\hat{M}_i\In E_{\hat{G}}(A_i^\ast,B_{i,2})$, $\phi(\hat{M}_i)\In C_{i,2}^\vartriangle$ and $|\hat{M_i}|\le \xi^{1/3}n$. This can be done greedily. Indeed, suppose we want to assign to $e\in E(\hat{G})$ an index $i\in[t]$. Let $e=xy$ and $c:=\phi(e)$. By~\ref{edge cover down} and~\eqref{final A set incidences}, the number of $i\in [t]$ for which $e\in E_{\hat{G}}(A_i^\ast,B_{i,2})$ and $c\in C_{i,2}^\vartriangle$, is at least~$\mu^3 n$.
By Claim~\ref*{claim:bounded leftover}, we have that $|E_c(\hat{G})|, d_{\hat{G}}(x),d_{\hat{G}}(y)\le \sqrt{\xi} n$ and $|E(\hat{G})|\le \sqrt{\xi} n^2$. Thus, there exists a suitable $i\in[t]$ such that no other $c$-edge of $\hat{G}$ has been assigned to~$i$, and $\hat{M}_i$ does not yet cover $x$ or~$y$ and contains at most $\xi^{1/3}n/2$ edges so far. Finally, by~\eqref{leftover graph} we have that
\begin{align}
	\mbox{$E(G_2') \In E(G_2'\cup G^{\rm{rb}} \cup G^\vartriangle_2 \cup G^\vartriangle_3) \In E\left(K_n-\bigcup_{i\in[t]}(F_i^\ast\cup \hat{M}_i)\right) \In E(G_1'\cup G_2')$.}\label{leftover graph end}
\end{align}

\smallskip

\noindent\textbf{Cover the remaining non-reservoir colours.}
Next, find edge-disjoint matchings $\hat{M}_1',\dots,\hat{M}_t'$ in $G^\vartriangle_2$ such that, for each $i\in [t]$, $V(\hat{M}_i')\cap V(\hat{M}_i)=\emptyset$, $\hat{M}_i'\In E_{G^\vartriangle_2}(A_i^\ast,B_{i,2})$, and $\hat{M}_i'$ consists of exactly one $c$-edge for each $c\in \hat{C}_i$. (Hence, $|\hat{M}_i'|=|\hat{C}_i|$.)
This can also be done greedily. Indeed, suppose we want to add a $c$-edge to~$\hat{M}_i'$. By~\ref{colour cover down} and~\eqref{final A set size}, there are at least $\mu^3 n$ $c$-edges in $E_{G^\vartriangle_2}(A_i^\ast,B_{i,2})$. By Claim~\ref*{claim:bounded leftover}, we have that $|\hat{C}_i|\le \sqrt{\xi} n$ and $|\set{i\in[t]}{c\in \hat{C}_i}|\le \sqrt{\xi} n$. Also recall that $|\hat{M}_i|\le \xi^{1/3}n$.
Hence, there exists a suitable $c$-edge which has not been used by another matching~$\hat{M}_{i'}'$, and whose endvertices are not covered by $\hat{M}_i$ or yet by $\hat{M}_i'$. Hence, for each $i\in[t]$, we have by~\eqref{leftover colours} that
\begin{align}
	\mbox{$C_{i,2}' \cup D_i' \In C_{i,2}' \cup D_i' \cup C_{i,3}^\vartriangle \In C\sm \phi(E(F_i^\ast)\cup \hat{M}_i\cup \hat{M}_i')\In C_{i,1}'\cup C_{i,2}' \cup D_i'$.}\label{leftover colours end}
\end{align}

\smallskip

\noindent\textbf{Absorb the uncovered vertices.}
We now extend $F_i^\ast\cup \hat{M}_i\cup \hat{M}_i'$ into a spanning forest by adding a matching~$M^\vartriangle_i$.
For each $i\in[t]$, let $A_i':=A_i^\ast\sm V(\hat{M}_i\cup \hat{M}_i')$ and $B_i':=B_i^\ast\sm V(\hat{M}_i\cup \hat{M}_i')$. We aim to apply Lemma~\ref{lem:matching routine} (with $\mu,\beta_\vartriangle q_\vartriangle,\Set{A_i',B_i',G^{(i)}[A_i',B_i']}_{i\in[t]}$ playing the roles of $\mu,d,\Set{A_i,B_i,G_i}_{i\in[t]}$).\COMMENT{This is where we use $\mu\ll \eta$} (Recall that $G^{(i)}$ was defined just before~\ref{colour boundedness}.)
Clearly, $|A_i'|=|B_i'|=b-|\hat{M}_i\cup \hat{M}_i'|\ge (\mu- 3\xi^{1/3})n$. In particular, by~\ref{super regular}, $|A_i\sm A_i'|,|B_i\sm B_i'| \le 4\xi^{1/3}n$ and thus $G^{(i)}[A_i',B_i']$ is $(\mu^{1/3},\beta_\vartriangle q_\vartriangle)$-quasirandom.\COMMENT{$\mu^{1/3}$ could be replaced with something much better, this is how it's applied}
Finally, since $A_i'\In A_i$ and $B_i'\In B_i$, the remaining conditions for Lemma~\ref{lem:matching routine} follow immediately from \ref{colour boundedness}, \ref{jain:intersections} and~\ref{jain:incidences}.
Therefore, we can find edge-disjoint $M^\vartriangle_1,\dots,M^\vartriangle_t$ such that $M^\vartriangle_i$ is a rainbow perfect matching of $G^{(i)}[A_i',B_i']$  for each $i\in[t]$.

\smallskip

\noindent\textbf{Absorb the uncovered colours.}
Now, for each $i\in[t]$, let
\begin{align}
	\mbox{$M_i^\ast := M_i^\vartriangle \cup \hat{M}_i\cup \hat{M}_i'$\;\;\;\; and \;\;\;\;$C_i^\ast:= C_{i,1}'\sm \phi(E(F_i^\ast) \cup M_i^\ast)$.}\label{brush}
\end{align}
Note that by~\eqref{leftover colours end}, $M_i^\ast$ is a rainbow perfect matching between $A_i^\ast$ and~$B_i^\ast$. Similarly, by~\eqref{leftover colours end} and~\eqref{brush}, $F_i^\ast \cup M_i^\ast$ is rainbow and $\phi(E(F_i^\ast) \cup M_i^\ast) = C\sm (C_i^\ast \cup C_{i,2}' \cup D_i')$. Also, $\bigcup_{i\in[t]}F_i^\ast$ is edge-disjoint from $\bigcup_{i\in[t]}M_i^\ast$ (since $G^\vartriangle_1 \cupdot G^\vartriangle_2  \cupdot G^\vartriangle_3=G^\vartriangle$). Since $F_i^\ast \cup M_i^\ast$ has $(n-1)-r_i$ edges by~\eqref{missing edge count}, we deduce that $|C_i^\ast \cup C_{i,2}' \cup D_i'|=r_i$, implying that $|C_i^\ast|=s$ by~\eqref{def number absorbers} and \ref{rainbow matchings set sizes}. Therefore, using~\ref{colour absorber app}, there exists $J_i\In M_i$ such that $J_i$ is $(C_i^\ast \cup C_{i,2}')$-rainbow and contains exactly one edge from each of $\set{M_{i,j}}{j\in[3s]}$. Note that the $J_i$ are edge-disjoint from each other by~\ref{rainbow matchings partition 2 app} and also edge-disjoint from $\bigcup_{i\in[t]}(F_i^\ast \cup M_i^\ast)$ by~\eqref{leftover graph end}.

\smallskip

\noindent\textbf{Absorb the uncovered edges.}
Finally, let $$E^\ast:=E(G_1') \mathbin{\big\backslash} \bigcup_{i\in[t]}(F_i^\ast \cup M_i^\ast \cup J_i).$$
We claim that $E^\ast$ contains precisely $m$ $c$-edges for every $c\in C$. Note that, by~\eqref{leftover graph end}, we have $\bigcup_{i\in[t]}(F_i^\ast \cup M_i^\ast \cup J_i) = K_n - (E^\ast \cup E(G_2'))$. Moreover, for each $i\in[t]$, $F_i^\ast \cup M_i^\ast \cup J_i$ is $(C\sm D_i')$-rainbow, implying that for each $c\in C$ the number of $c$-edges in $\bigcup_{i\in[t]}(F_i^\ast \cup M_i^\ast \cup J_i)$ is $|\set{i\in[t]}{c\in C\sm D_i'}|=t-3m$ by~\ref{colour regularity app}. Thus, the number of $c$-edges in $E^\ast \cup E(G_2')$ is $3m$, which implies the claim, using \ref{colour regularity app} again.

Thus, by~\ref{edge absorber app}, there exists a partition of $E^\ast \cup E(G_2')$ into sets $J_1',\dots,J_t'$, such that for each $i\in[t]$, $J_i'$ contains exactly one edge from each of $\Set{M_{i,c}'}_{c\in D_i'}$. In particular, $J_i'$ is $D_i'$-rainbow.

Let $T_i:=F_i^\ast \cup M_i^\ast \cup J_i \cup J_i'$. By~$(\dagger)$, $T_i$ is a rainbow spanning tree isomorphic to~$T_{n;r,b}$, and $T_1,\dots,T_t$ decompose $K_n$, as desired.
\endproof

Finally, we briefly mention how the proof can be adapted to prove Theorem~\ref{thm:main} with $\Delta(T)=3$. The only necessary change is in how we connect the matchings in $\cM_i$ by using Lemma~\ref{lem:connecting}. Suppose that in Step~\ref*{step:connecting} in the proof of Theorem~\ref{thm:main}, we want to connect the `head set' $H(M)$ with the `tail set' $T(M')$ for two consecutive $M,M'\in \cM_i$. In the current proof, we find a vertex $w$ and internally disjoint paths of length~$2$ from $w$ to each vertex in $H(M)\cup T(M')$. Instead, we could also connect $H(M)\cup T(M')$ as follows: let $B$ be a binary tree with root $b$ and leaves $H(M)$, and let $B'$ be a binary tree with root $b'$ and leaves $T(M')$, and such that $V(B)\cap V(B')=\emptyset$. (Recall that $|H(M)|=|T(M')|=2^8$.) Let $R$ be the graph obtained from $B\cup B'$ by adding a path of length~$2$ between $b$ and $b'$, and then subdividing every edge once. Clearly, $\Delta(R)\le 3$, and this construction ensures that still, the tree $T_i$ is always the same, independent of which edge of $M$ is ultimately selected for~$T_i$.
To find all the required connections~$R$, one could still employ Lemma~\ref{lem:connecting}, here repeatedly, with $k=2$.
However, this necessitates to split $V$, $C$ and $E(K_n)$ into even more subsets, so for clarity, we omitted this from the proof.

%
%

\providecommand{\bysame}{\leavevmode\hbox to3em{\hrulefill}\thinspace}
\providecommand{\MR}{\relax\ifhmode\unskip\space\fi MR }
\providecommand{\MRhref}[2]{%
  \href{http://www.ams.org/mathscinet-getitem?mr=#1}{#2}
}
\providecommand{\href}[2]{#2}

\end{document}